\documentclass[a4paper,11pt]{article}
\usepackage{amssymb}
\usepackage{stackengine}
\stackMath
\usepackage{latexsym}
\usepackage{amsmath}
\usepackage{amsthm}
\usepackage{amsfonts}
\usepackage{amssymb}

\usepackage{pstricks}

\usepackage{graphicx}

\newtheorem{Th}{Theorem}
\newtheorem{Prop}{Proposition}

\newtheorem{Lm}{Lemma}[section]
\newtheorem{Lma}{Lemma}[section]

\newcommand{\be}{\begin{equation}}
\newcommand{\ee}{\end{equation}}
\newcommand{\bes}{\begin{equation*}}
\newcommand{\ees}{\end{equation*}}

\newcommand{\R}{\mathbb{R}}
\newcommand{\N}{\mathbb{N}}
\newcommand{\C}{\mathbb{C}}
\newcommand{\Z}{\mathbb{Z}}

\newcommand{\reset}{\setcounter{equation}{0}\setcounter{Th}{0}\setcounter{Prop}{0}\setcounter{Co}{0}
\setcounter{Lm}{0}\setcounter{Rm}{0}}

\def\La{\Lambda}
\def\La{\Lambda}
\def\ti{\tilde}
\def\lf{\left}
\def\rg{\right}

\def\al{\alpha}
\def\la{\lambda}

\def\ep{\varepsilon}

\def\ds{\displaystyle}
\def\ov{\overline}
\def\Om{\Omega}
\def\om{\omega}
\def\p{\partial}

\def\Xint#1{\mathchoice
{\XXint\displaystyle\textstyle{#1}}
{\XXint\textstyle\scriptstyle{#1}}
{\XXint\scriptstyle\scriptscriptstyle{#1}}
{\XXint\scriptscriptstyle\scriptscriptstyle{#1}}
\!\int}
\def\XXint#1#2#3{{\setbox0=\hbox{$#1{#2#3}{\int}$ }
\vcenter{\hbox{$#2#3$ }}\kern-.6\wd0}}

\def\dashint{\Xint-}

\makeatletter
\newcommand{\tpitchfork}{%
  \vbox{
    \baselineskip\z@skip
    \lineskip-.52ex
    \lineskiplimit\maxdimen
    \m@th
    \ialign{##\crcr\hidewidth\smash{$-$}\hidewidth\crcr$\pitchfork$\crcr}
  }%
}
\makeatother

\begin{document}

\title{Morse Index Stability for Critical Points to Conformally Invariant Lagrangians }

\author{Francesca Da Lio, Matilde Gianocca  and Tristan Rivi\`ere\footnote{Department of Mathematics, ETH Zentrum,
CH-8093 Z\"urich, Switzerland.}}

\date{ }
\maketitle

{\bf Abstract :}{\it We prove the upper-semi-continuity of the Morse index plus nullity of critical points to general conformally invariant Lagrangians in dimension 2 under weak convergence.  Precisely we establish that the sum of the  Morse indices and the nullity of an arbitrary sequence of  weakly converging critical points to a general conformally invariant Lagrangians of maps from an arbitrary closed surface into an arbitrary closed smooth manifold passes to the limit in the following sense : it is asymptotically bounded from above by the sum of the Morse indices plus the nullity of the weak limit and the bubbles, while it was well known that the sum of the Morse index of the weak limit with the Morse indices of the bubbles is asymptotically bounded from above by the Morse indices of the weakly converging sequence.  The main result is then extended to the case of sequences of maps from sequences of domains degenerating to a punctured Riemann surface assuming that the lengths of the images by the maps of the collars  associated to this degeneration stay below some critical length.}


\section{Introduction}
\subsection{General Framework and Main Result}
Let $(\Sigma,h)$ be an arbitrary smooth closed and oriented 2-dimensional Riemannian manifold  and let $({\mathcal{N}}^n,g)$ be an arbitrary closed smooth Riemannian manifold of arbitrary dimension $n$ that we can assume to be isometrically embedded in some Euclidian Space ${\R}^m$ thanks to Nash theorem.

We consider the general form of conformally invariant Lagrangians of maps from $(\Sigma,h)$ into $({\mathcal{N}}^n,g)$ given by Gr\"uter\footnote{A classical result of Gr\"uter \cite{Gru} asserts that any strictly elliptic lagrangian of quadratic growth from $\Sigma$ into ${\mathcal{N}}^n$ is given  by (\ref{0.0001}) for some metric $g$ in ${\mathcal{N}}^n$.}
\be
\label{0.0001}
\mathfrak{L}(u):=\frac{1}{2}\int_\Sigma |du|^2_h\ dvol_h+u^\ast\al
\ee
where $\al$ is an arbitrary 2-form of ${\mathcal{N}}^n$. Observe that ${\frak L}(u)$ is unchanged by replacing $h$ by any conformally equivalent metric on $\Sigma$. Hence ${\frak L}$ can be considered as a lagrangian of maps from a given {\bf Riemann Surface} into $({\mathcal{N}}^n,g)$ and very often in the literature $h$ is chosen to have constant $0,+1, -1$ Gauss curvature thanks to the uniformization theorem. We are interested in critical points to $\frak{L}$ within the non linear Sobolev Space
\[
W^{1,2}(\Sigma, {\mathcal{N}}^n):=\lf\{u\in W^{1,2}(\Sigma,{\R}^m)\ ;\ u(x)\in {\mathcal{N}}^n\quad\mbox{a.e.}\rg\}.
\]
Let $H$ be the section of $TN\otimes\wedge^2T^\ast N$ such that
\be
\label{0.01}
\forall \,X,Y,Z\in T_z{\mathcal{N}}^n\quad\ d\al(X,Y,Z):= \lf<X, H_z(Y,Z)\rg>
\ee
where $<\cdot,\cdot>$ denotes the scalar multiplication in ${\R}^m$. The corresponding Euler-Lagrange equation is given by
\be
\label{0.1}
\Delta_hu+H_u(du\wedge du)_h= -\,{\mathbb I}_u(du,du)_h\quad\mbox{ in } {\mathcal D}'(\Sigma)\ .
\ee
where $\Delta_h$ is denoting the positive Laplace Beltrami operator and ${\mathbb I}_z(X,Y)$ is the second fundamental form of the embedding ${\mathcal{N}}^n\hookrightarrow {\R}^m$ at the point $z\in {\mathcal{N}}^n$ and taken over a pair of vectors $X,Y\in T_z{\mathcal{N}}^n$. The non linearity $ {\mathbb I}_u(du,du)_h$ is  equal to $\mbox{Tr}_h({\mathbb I}_u(du\otimes du))$. Using local coordinates in $\Sigma$ it is given by $${\mathbb I}_u(du,du)_h:=\sum_{i,j=1}^2 h^{ij}\,{\mathbb I}_u(\p_{x_i}u,\p_{x_j}u)\ .$$
and we define for any maps $f,g$ from $\Sigma$ into ${\R}^m$ using local coordinates
\[
H(df\wedge dg)=\lf[H_u(\p_{x_1}f,\p_{x_2}g)-H_u(\p_{x_2} f,\p_{x_1}g)\rg]\, dx_1\wedge dx_2\ ,
\]
 and\footnote{Observe that
 \[
 H_u(\p_{x_1}f,\p_{x_2}g)-H_u(\p_{x_2} f,\p_{x_1}g)=H_u(\p_{x_1}f,\p_{x_2}g)+H_u(\p_{x_1} g,\p_{x_2}f)
 \]which gives in particular $H_u(df\wedge dg)=H_u(dg\wedge df)$.} 
 \[
  H(df\wedge dg)_h\ dvol_h:=H(df\wedge dg)\ .
\]
If one denotes $P_z$ the $m\times m$ real symmetric matrix\footnote{We shall denote $S_m({\R})$ the space of real $m\times m$ real symmetric matrices.} given by the orthogonal projection of ${\R}^m$ onto $T_z{\mathcal{N}}^n$, for any map $u\in W^{1,2}(\Sigma, {\mathcal{N}}^n)$, the map $P_{u(x)}$ is clearly in $W^{1,2}(\Sigma,S_m({\R}))$ and one proves (see \cite{Hel}) that (\ref{0.1}) is equivalent to
\be
\label{0.2}
P_u\lf[\Delta_h u+H_u(du\wedge du)_h\rg]=0\quad\mbox{ in } {\mathcal D}'(\Sigma)\ .
\ee
Solutions to (\ref{0.2}) are known to be smooth \cite{Riv1}. Any solution to (\ref{0.2}) which is in addition conformal corresponds to a possibly branched immersion of $\Sigma$ with  {\it  mean curvature vector} $H(\vec{e}_1,\vec{e}_2)$ where $(\vec{e}_1,\vec{e}_2)$ is an orthonormal basis of the tangent space in $T{\mathcal{N}}^n$ to $u_\ast T\Sigma$.  For that reason, equation (\ref{0.2}) is called   {\it Prescribed Mean Curvature Equation} into ${\mathcal{N}}^n$. Solutions to the {\it Prescribed Mean Curvature Equation} into manifolds that we will simply call PMC solutions into ${\mathcal{N}}^n$  are relatively common mathematical objects from geometric calculus of variations. They are playing a special role in particular in the differential geometry of surfaces, Teichm\"uller space theory and uniformization as well as in mathematical models from general relativity or condensed matter physics.

Beside the value of it's energy of course there are two other important numbers attached to every harmonic map : the {\it Morse index} and the {\it nullity}. They are respectively the dimension of the largest space on which the second derivative of $\mathfrak{L}$ is negative and the dimension of it's kernel. Precisely, the second derivative of $\mathfrak{L}$ at a PMC map $u$ is defined on the space of infinitesimal variations of $u$ which are nothing but the sections of the pull-back by $u$ of the tangent bundle to ${\mathcal{N}}^n$ :
\be
\label{0.2-a}
V_u:=\Gamma(u^{-1}T{\mathcal{N}}^n):=\lf\{{w}\in W^{1,2}(\Sigma, {\mathcal{N}}^n)\ ;\ P_{u(x)}{w}=w \mbox{ for a. e. }x\in \Sigma\rg\}\ .
\ee 
The second derivative of
\[
E(u):=\frac{1}{2}\int_\Sigma |du|^2_h\ dvol_h
\]
at $u$ is a quadratic form on $V_u$ given by 
\be
\label{0.3}
D^2E_u(w)=\int_\Sigma |dw|_h^2\ dvol_h-\int_\Sigma \lf<{\mathbb I}_u(du,du)_h   , {\mathbb I}_u(w,w)\rg>\ dvol_h
\ee
where $<\cdot,\cdot>$ denotes the scalar multiplication in ${\R}^m$. In the sphere case for instance when ${\mathcal{N}}^n=S^n\subset {\R}^{n+1}$ the computations give
\be
\label{0.3-a}
D^2E_u(w)=\int_\Sigma|dw|^2_h- |du|^2_h\,|w|^2\, dvol_h\ .
\ee
In the general case, for $\om$ arbitrary, 
\be
\label{0.3-b}
\begin{array}{l}
\ds Q_u(w):=D^2\frak{L}_u(w)=\int_\Sigma |dw|^2_h\ dvol_h-\int_\Sigma \lf<{\mathbb I}_u(du,du)_h   , {\mathbb I}_u(w,w)\rg>\ dvol_h\\[5mm]
\ds\quad+\int_\Sigma w\cdot \lf[2\,H_u(dw\wedge du)_h+\nabla_w H_u(du\wedge du)_h\rg]\ dvol_h\ .
\end{array}
\ee
We define the {\it Morse index} of $u$ to be
\[
\mbox{Ind}_{\frak{L}}(u):=\sup\lf\{\mbox{dim}(W) \ ;\ W\mbox{ is a sub vector-space of }V\mbox{ s.t. }\lf.Q_u\rg|_{W}<0\rg\}
\]
and the {\it Nullity}
\[
\mbox{Null}_{\frak{L}}(u):=\mbox{dim}\,Ker (B_u)
\]
where $B_u$ is the bilinear form associated to $Q_u$.

The present paper is aiming at studying  the behaviour of $\mbox{Ind}_L(u_k)+\mbox{Null}_L(u_k)$  for  sequences of harmonic maps with uniformly bounded energy.

\medskip

A classical result in {\it concentration compactness theory} (see for instance \cite{LaRi1} theorem 3) asserts that modulo extraction of a subsequence there exists a limiting harmonic map $u_\infty$ from $\Sigma$ into ${\mathcal{N}}^n$, there exist $p^1\ldots p^Q\in \Sigma$, there exist
$x^1_k\ldots x^Q_k\in \Sigma$  such that $x^j_k\rightarrow p^j$ for any $j=1\ldots Q$, there exist $\delta^1_k\ldots \delta^Q_k$ positive numbers converging to zero and there exist
$v_\infty^1\ldots v_\infty^Q$ harmonic maps from ${\C}$ into ${\mathcal{N}}^n$ such that
\be
\label{0.4}
\lim_{k\rightarrow 0}\lf\|\nabla\lf(u_k-u_\infty- \sum_{j=1}^Q v^j_k\rg)\rg\|_{L^2(\Sigma)}=0\ ,
\ee
where, using normal coordinates for the metric $h$ around each $x^j_k$, we denote $v_k^j(\delta^j_k y+x_k^j):=v^j_\infty(y)$. The harmonic maps $v^j_\infty$ are called the {\it bubbles} and the convergence
(\ref{0.4}) which is proper to conformally invariant variational problems is called {\it bubble tree convergence} of $u_k$ towards $(u_\infty, v_\infty^1\ldots v_\infty^Q)$.

Our first main result in the present work is the following
\begin{Th}
\label{th-morse-stability}
Let $u_k$ be a sequence of solutions to the prescribed mean curvature equation (\ref{0.2}) from a closed oriented surface $(\Sigma,h)$ into a closed arbitrary $C^2$ riemannian manifold ${\mathcal{N}}^n$ for $H\in C^1$. Assume $u_k$ bubble tree converges towards $(u_\infty, v_\infty^1\ldots v_\infty^Q)$. Then, for $k$ large enough one has
\be
\label{0.5}
\mbox{Ind}_{\frak{L}}(u_k)+\mbox{Null}_{\frak{L}}(u_k)\le \mbox{Ind}_{\frak{L}}(u_\infty)+\mbox{Null}_{\frak{L}}(u_\infty)+\sum_{j=1}^Q\mbox{Ind}_{\frak{L}}(v_\infty^j)+\mbox{Null}_{\frak{L}}(v_\infty^j)\ .
\ee 
\hfill $\Box$
\end{Th}
Observe that, in view of the explicit form (\ref{0.3-b}) of $D^2{\frak{L}}$, the regularity assumptions we are making on ${\mathcal{N}}^n$ and $\al$ are the most general ones for ensuring the {\bf continuity of $D^2\frak{L}$}. Thanks to the method we are using, we don't require more regularity on $D^2\frak{L}$. 

Theorem~\ref{th-morse-stability} has been first established in \cite{yi-0} and \cite{yi-1} in the case of harmonic maps into $C^3$ manifolds based on pointwise estimates established in \cite{QiTi} (see more comment on these estimates below).

Recall that classical results (see \cite{MoRe} or more recently \cite{KaSt} for $\frak{L}=E$) imply that for $k$ large enough one has the lower-semi-continuity of the Morse index in the following sense
\be
\label{0.6}
\mbox{Ind}_{\frak{L}}(u_\infty)+\sum_{j=1}^Q\mbox{Ind}_{\frak{L}}(v_\infty^j)\le \mbox{Ind}_{\frak{L}}(u_k)\ .
\ee
The lower semi-continuity of the Morse index (\ref{0.6}) is also true while considering harmonic maps approximations and holds in general for limits of sequences solving some relaxed min-max problem
converging to possible non zero Morse index harmonic maps. Inequality (\ref{0.6}) is a rather general and ``robust'' inequality. It holds in various frameworks like for instance  in the viscosity method for minimal surfaces, \cite{Riv2}.

The upper-semi-continuity of the Morse Index plus nullity is a way more delicate to obtain than the lower semi-continuity of the Morse index and does not always hold in fact.
It has been established in two remarkable works  on minimal surface theory respectively by O. Chodosh and C. Mantoulidis for sequences of critical points to the Allen-Cahn Functional  (see \cite{ChoMa} Theorem 1.9)
and by F.C. Marques and A. Neves for limits of Almgren-Pitts minmax procedures (see \cite{MaNe} Theorem 1.3). In both cases a main assumption is that the limiting minimal surface has multiplicity one. It has been disproved when the limit has multiplicity higher than one (see \cite{dPKWY} and \cite{ChoMa} Example 5.2). 

 Inequality (\ref{0.5}) is not expected to hold for instance for general Palais-Smale harmonic map sequences or for sequences of harmonic maps from degenerating Riemann Surfaces. Indeed, in both cases
 {\it Energy quantization} (\ref{0.4}) does not necessarily holds (see \cite{Par}, \cite{Zhu}) and the neck energy could contribute to the negative part of the spectrum of the Hessian $D^2E$. More seriously, even when (\ref{0.4}) holds for some Sacks-Uhlenbeck type $\al-$harmonic maps, it is not excluded that the necks are given by geodesics with  non zero lengths (see  \cite{Moo} ). In such a case it could be that these geodesics contribute to the negative part of the spectrum of the Hessian $D^2E$.
 Our analysis throughout the present work is very much based on the {\it $C^0-$no neck property} for weakly converging critical points of $\frak L$.  This {\it $C^0$-no neck property} property is saying that there is no distance between the bubbles $v^j_\infty$ and the main part $u_\infty$, roughly speaking the bubbles $v_\infty^j$ are ``$L^\infty-$glued'' directly to the weak limit $u_\infty$. 
  It was first proved in the special case of harmonic maps into manifolds in a remarkable paper by J. Qing and G. Tian \cite{QiTi}. The main idea in this paper was to obtain pointwise estimates of the gradient of the harmonic maps in the neck region by bounding  the difference between the harmonic map and geodesic maps into ${\mathcal{N}}^n$ (which are nothing but axially symmetric harmonic maps from the neck regions into ${\mathcal{N}}^n$) thanks to a very refined {\it three circle theorem} by L.Simon \cite{Sim}. In a later work \cite{LiWa}, F.H. Lin and C.Y. Wang found a different strategy for proving the {\it $C^0$-no neck property} under the same assumptions than the ones made in \cite{QiTi}. They used mostly the fact that the angular energy of harmonic maps on surrounding circles in the neck regions satisfy a remarkable ordinary differential inequality, as it had been previously observed by T.Parker in \cite{Par}. Both  {\it $C^0-$no neck property} proved in \cite{QiTi} and \cite{LiWa} are valid for sequences of harmonic maps into  manifolds which are at least $C^3$ regular.
  
  \medskip
 
 A more general situation where strong convergence of sequences of critical points of $\frak L$ fails, arises when instead of considering $u_k$ from a fixed Riemann Surface $(\Sigma,h)$ one considers sequences of critical points $u_k$ of $\frak L$ from a sequence of Riemann surfaces $(\Sigma,h_k)$ into $({\mathcal{N}}^n,g)$ and when the sequence of Constant Gauss Curvature metric $h_k$ is degenerating in the Moduli space of $\Sigma$. The corresponding bubbling analysis in this situation has been carefully studied by M. Zhu in \cite{Zhu} for ${\frak L}=E$. According to the general moduli space compactification theory (see for instance \cite{Hum}) the degeneration of the Riemann surface $(\Sigma,h_k)$ is equivalent to the convergence to zero of the length of the shortest geodesics. One can extract a subsequence such that there exists a fixed number of pairwise disjoint closed embedded geodesics $(L_k^j)_{j=1\ldots P}$  of lengths $(l_k^j)_{j=1\ldots P}$ all tending to zero and conformal charts $\psi^j_k$ from  annuli ${\mathcal C}_k^j$ of degenerating conformal classes containing $L_k^j$ and called the {\bf collars} into $B_1(0)\setminus B_{\delta_k^j}(0)$ where $\delta_k^j:=\exp(-1/{l_k^j})$. Away from these collars the bubble tree analysis (\ref{0.4}) holds locally with respect to the {\bf nodal Riemann surface} obtained by collapsing the geodesics $L_k^j$ to points called {\bf punctures} (see \cite{Zhu} in the harmonic map case with $C^3$ targets and \cite{LaRi1} for the general case).  In order to simplify the presentation we assume that no bubble is formed in the collars that is to say, in the chart $\psi_k^j$,
 \[
 \lim_{\eta\rightarrow 0}\,\lim_{k\rightarrow+\infty}\sup_{\rho\in [\eta^{-1}\,\delta_k^j,\eta]}\int_{B_{2\rho}(0)\setminus B_\rho(0)}|\nabla u_k|^2\ dx^2=0\ .
 \]
 The {\bf limiting  average length } of  the images by $u_k$ of the collar ${\mathcal C}_k^j$ is given  by
 \[
\La^j:=\lim_{\eta\rightarrow 0}\,\limsup_{k\rightarrow +\infty}\int_{\eta^{-1}\,\delta_k^j}^\eta\frac{1}{2\pi}\int_{0}^{2\pi}\lf|\frac{d{u}_k}{dr}\rg|\ dr \ ,
 \]
 where the coordinates $(r,\theta)$ are referring to the charts $\psi_k^j$. Our second main result in this work is the following theorem

\begin{Th}
\label{th-morse-stability-deg}
Let $u_k$ be a sequence of solutions to the prescribed mean curvature equation (\ref{0.2}) from a sequence of closed  Riemann surfaces $(\Sigma,h_k)$ of genus larger than one into a closed arbitrary $C^2$ riemannian manifold ${\mathcal{N}}^n$ for $H\in C^1$. Assume $(\Sigma,h_k)$ is degenerating  to a punctured surface with a fixed number of collars $({\mathcal C}_k^j)_{j=1\ldots P}$ at the limit.  Assume $u_k$ bubble tree converges away from the collars towards $(u_\infty, v_\infty^1\ldots v_\infty^Q)$ and that no bubble is formed in the collars. There exists a constant $\La^\ast>0$ depending only on $m$, ${\mathcal{N}}^n$ and $H$ such that, if the  { limiting average  length } $\Lambda^j$ of  the images by $u_k$ of the collar ${\mathcal C}_k^j$  are less than $\La^\ast$ then, for $k$ large enough one has
\be
\label{0.5-b}
\mbox{Ind}_{\frak{L}}(u_k)+\mbox{Null}_{\frak{L}}(u_k)\le \mbox{Ind}_{\frak{L}}(u_\infty)+\mbox{Null}_{\frak{L}}(u_\infty)+\sum_{j=1}^Q\mbox{Ind}_{\frak{L}}(v_\infty^j)+\mbox{Null}_{\frak{L}}(v_\infty^j)\ .
\ee 
\hfill $\Box$
\end{Th}
The result can be understood in the following sense. As explained in \cite{Zhu}, in the harmonic map case $E={\frak L}$, the images of the collars converge towards portions of geodesics of ${\mathcal{N}}^n$. It is a classical result that below a critical length all geodesic arcs are stable and then the corresponding collar does not  contribute to the negative part of the spectrum. This is shedding some light on the reason why the upper-semi-continuity of the dimension of the negative part of the spectrum obtained for sequences of critical points of ${\frak L}$ from a fixed surface - i.e. (\ref{0.5}) - does extend to sequences of domains under the assumption that $\Lambda^j<\Lambda^\ast$ for any $j$. We believe that this assumption is optimal since geodesics beyond some length have negative Morse index and the collars could then contribute to the negative part of the spectrum, which could eventually break (\ref{0.5-b}).

 \medskip
 
\noindent{\bf The role of interpolation spaces in neck/collar analysis} : In the present work the  {\it $C^0$-no neck property}  is derived from the main estimates established by the 3rd author and P. Laurain in \cite{LaRi1} (see  the explanations in subsections III.1 and IV.1 below)\footnote{See relevant estimates by T. Lamm and B. Sharp in \cite{LaSh}.}. The {\it $C^0$-no neck property} in our approach is a consequence of the so called {\it $L^{2,1}-$energy quantization}\footnote{The  {\it $L^{2,1}-$energy quantization} (i.e. lemma~\ref{lm-L-2-1quanti-general} below) was first established in \cite{LaRi1} for the ``angular part of the energy'' in the necks.} where $L^{2,1}$ refers to the {\it Lorentz space} pre-dual  to the {\it Marcinkiewicz}  space $L^2-$weak also denoted $L^{2,\infty}$.  
 
  The idea to introduce the Lorentz space $L^{2,1}$ norm in neck analysis for conformally invariant variational problems goes back to the two works of the third authors with F.H.Lin (\cite{LiRi1} and \cite{LiRi2}) while the general derivation of conservation laws for harmonic maps, which is also central in the present work, has been derived by the third author in \cite{Riv1} .

\medskip

There are strong reasons to expect  that the {\it Morse index plus nullity upper-semi-continuity} holds true for a large class of conformally invariant variational problems for which conservation laws opening the way to {\it $C^0$-no neck property} have been derived in the last two decades. We can quote for instance  {\bf half-harmonic map} in one dimension and {\bf free boundary surfaces} for which almost conservation laws have been derived in \cite{DaRi1}, \cite{DaRi2}    by the first and the third authors, in \cite{DaPi}, by   the first author and Pigati and   in   \cite{DMS} by the first author, Mazowiecka and Schikorra, {\bf Bi-harmonic maps} in 4 dimensions where the {\it $C^0$-no neck property} should follow from the analysis in \cite{LaRi2} which it-self is derived 
from the conservation laws established in \cite{LamRi}, {\bf Yang-Mills Fields} in dimension 4 where the $L^{2,1}-$norm estimates of the curvature in neck regions is established by the third author in \cite{Riv0},  or {\bf Willmore surfaces} where the {\it $C^0$-no neck property} for the Gauss map has been proved recently by the third author in collaboration with A. Michelat in \cite{MiRi}.

\bigskip

One of the motivation of the present work is to ultimately establish the upper semi-continuity of the Morse index plus nullity for solution to {\em min-max problems}  via the Ginzburg Landau type relaxation $\frak{L}_\ep$ of the energy $\frak{L}$ of the form
\[
{\frak L}_\ep(u):=\frac{1}{2}\int_\Sigma \lf[|du|^2_h+\frac{1}{\ep^2}F(u)\rg]\ dvol_h+u^\ast\al\ ,
\]
for $u$ in $W^{1,2}(\Sigma,{\R}^m)$ where $F$ is some smoothing of the square of the distance to ${\mathcal{N}}^n$. This question is under investigation by the authors. It would be interesting to see under which assumption the {\it $C^0$-no neck property} holds. We expect the   {\it entropy estimates} of the form
\[
\frac{1}{\ep^2}F(u)=o\lf(\frac{1}{\log\ep^{-1}}\rg)\ ,
\]
to play a role in this context. Such {\it entropy estimate} can be derived from the refinement of Palais-Smale deformation theory due to M.Struwe also known as ``Struwe Monotonicity Trick''.

\subsection{Description of the various steps in the proof of theorem~\ref{th-morse-stability}.} 
\bigskip
We consider a weakly converging sequence of solution to (\ref{0.1}).
We assume for simplicity  of the presentation that there is exactly one bubble of characteristic size  $\delta_k\rightarrow 0$ and centered at $x_k\rightarrow p\in \Sigma$. The so called ``neck region'' connecting
the bubble to the main part of the solution is given by the family of annuli with degenerating conformal class of the form
\[
A(\eta,\delta_k):=B_{\eta}(x_k)\setminus B_{\delta_k/\eta}(x_k)\ .
\]
The absence of bubble in the neck region (by definition) says
\[
\lim_{\eta\rightarrow 0}\limsup_{k\rightarrow +\infty}\ \sup_{ \delta_k/\eta<\rho<2\rho<\eta}\ \int_{B_{2\rho}(x_k)\setminus B_\rho(x_k)}|du_k|^2_h\ dvol_h=0\ .
\]
The main argument in \cite{LaRi1} is implying
\be
\label{0.8}
\lim_{\eta\rightarrow 0}\limsup_{k\rightarrow +\infty}\  \int_{A(\eta,\delta_k)}|du_k|^2_h\ dvol_h=0\ .
\ee
which is also called ``$L^2-$energy quantization''. 

\medskip

One of the main contribution of \cite{Riv1}  is the starting point of the present work. It says that the Euler-Lagrange Equation (\ref{0.1})  can be rewritten in the form
\be
\label{0.7}
\Delta_hu_k=\Om_k\cdot du_k\ ,\quad\mbox{ where }\quad\Om_k\in (u_k)_\ast\, T_x\Sigma\otimes so(m)\quad\mbox{ and }\quad |\Om_k|\le C\, |du_k|\ .
\ee
Thanks to the small Dirichlet energy of $u_k$ in the neck, one of the main result of \cite{Riv1} implies the existence of $A_k\in W^{1,2}(A(\eta,\delta_k),Gl_m({\R}))$
and $B_k\in W^{1,2}(A(\eta,\delta_k),M_m({\R}))$ such that the Euler Lagrange Equation (\ref{0.7}) can be rewritten in the form of a {\bf conservation law}
\be
\label{0.9}
\mbox{div}\lf(A_k\ \nabla u_k+B_k\ \nabla^\perp u_k \rg)=0,
\ee
where in \eqref{0.9} we are using local conformal/complex coordinates $z$ for $(\Sigma,h)$. In \cite{LaRi1}, the conservative form of the equation is used  to ``upgrade'' the {\it $L^2-$energy quantization} (\ref{0.8}) to the following {\it $L^{2,1}-$energy quantization}
\be
\label{0.9-a}
\lim_{\eta\rightarrow 0}\limsup_{k\rightarrow +\infty}\  \int_0^{+\infty}\lf|\lf\{x\in A(\eta,\delta_k)\ ;\ |du_k|_h(x)>s\rg\}\rg|^{1/2}\ ds=0\ .
\ee
This {\it $L^{2,1}-$energy quantization} is a crucial step for proving the the neck contributes in a \underbar{positive} way to the second derivative of the Lagrangian. Indeed we first proceed to a Hodge type decomposition
of $A_k\,\nabla u_k$ in the form
\be
\label{0.10}
A_k\,\nabla u_k=\nabla \varphi_k+\nabla^\perp\psi_k+\nabla \frak{h}_k\quad\mbox{ in }A(\eta,\delta_k)\ ,
\ee
where respectively
\be
\label{0.11}
\Delta\varphi_k=\nabla^\perp B_k\cdot\nabla u_k\quad,\quad\Delta\psi_k=\nabla^\perp A_k\cdot\nabla u_k\quad\mbox{ and }\quad \Delta\frak{h}_k=0\ .
\ee
The harmonic part is decomposed as follows
\[
\frak{h}_k=\frak{h}_k^++\frak{h}_k^-+\frak{h}_k^0\ \ \mbox{ where }\ \ \frak{h}_k^+=\Re\sum_{n>0} h_{n,k}z^n\ \ ,\ \ \frak{h}_k^-=\Re\sum_{n<0} h_{n,k}z^n\ \ \mbox{ and }\ \ \frak{h}_k^0=\frak{h}_{0,k}+C_{\eta,k}^0\,\log|z|\ .
\]
An important fact in the proof of the main result is to prove that the ``$L^{2,1}-$energy quantization'' is implying the following control of the constant $C_{\eta,k}^0$
\be
\label{0.12}
C_{\eta,k}^0=o_{\eta,k}\lf(\frac{1}{\log \frac{\eta^2}{\delta_k}}\rg)
\ee
where we observe that $\log \frac{\eta^2}{\delta_k}$ is nothing but the degenerating class of the annulus $A(\eta,\delta_k)$. The estimate \eqref{).12}  is optimal  for our purposes in the following sense. The contribution to the pointwise estimate of $|du_k|^2(x)$ of $\frak{h}_k^0$ is equal to $(C_{\eta,k}^0)^2/|x|^2$ and the following general result\footnote{We prove that the constant $\pi^2$ is optimal in the first inequality (\ref{0.13}). } holds
\be
\label{0.13}
\forall f\in W^{1,2}_0(A(\eta,\delta_k))\quad \int_{A(\eta,\delta_k)}|\nabla f|^2\ dx^2\ge \frac{\pi^2}{\log^2 \frac{\eta^2}{\delta_k}}\,\int_{A(\eta,\delta_k)}\frac{|f|^2}{|x|^2}\ dx^2>>\ \int_{A(\eta,\delta_k)} |d\frak{h}_k^0|^2\ |f|^2\ dx^2\ .
\ee
For the other contributions to $\nabla u_k$ we use on one hand  classical pointwise estimates for harmonic functions with only positive or negative Fourier frequencies and on the other hand we are using weighted Wente type inequalities (see Lemma \ref{WenteWeight}) for controlling the contributions respectively from $\varphi_k$ and $\psi_k$ whose Laplacians are both given by Jacobians.  Finally we  conclude the {\bf pointwise estimates of the gradient} of $u_k$ by proving the existence of $\beta>0$ independent of $k$ such that
\be
\label{0.14}
|du_k|^2\le\ o_{\eta,k}(1)\ \om_{\eta,k}(x)\ ,
\ee
where $\om_{\eta,k}$ is the {\bf weight function} given in $A(\eta,\delta_k)$ explicitly by
\be
\label{0.15}
\om_{\eta,k}(x)=\frac{1}{|x|^2}\,\lf[\frac{|x|^\beta}{\eta^\beta}+\frac{\delta_k^\beta}{\eta^\beta\,|x|^\beta}+\frac{1}{ \log^2 \frac{\eta^2}{\delta_k}}\rg]\ .
\ee
We generalize (\ref{0.13}) by proving that the first eigenvalue of the Laplacian for Dirichlet boundary conditions on the degenerating annulus $A(\eta,\delta_k)$ with respect to the $L^2-$norm weighted by 
$\om_{\eta,k}$ is bounded from above and below by positive constants \underbar{independent of $k$} in particular we have
\be
\label{0.16}
\exists\ \lambda_0>0\ \quad\forall f\in W^{1,2}_0(A(\eta,\delta_k))\quad \int_{A(\eta,\delta_k)}|\nabla f|^2\ dx^2\ge\lambda_0\ \int_{A(\eta,\delta_k)}\om_{\eta,k}\ | f|^2\ dx^2
\ee
where $\lambda_0$ is independent of $k$. As a consequence of (\ref{0.14}) and (\ref{0.16}) we obtain that the neck are asymptotically not contributing to the negativity of the second derivative of the Lagrangian in the following sense
\be
\label{0.17}
\forall\ w\in \Gamma(u_k^{-1} T{\mathcal{N}}^n)\cap W^{1,2}_0(A(\eta,\delta_k),{\R}^m)\quad\quad D^2{\frak{L}}_{u_k}(w)\ge\lambda_0\ \int_{A(\eta,\delta_k)}\om_{\eta,k}\ |w|^2\ dx^2\ .
\ee
The proof of the main theorem follows from (\ref{0.17}).
\medskip

{\bf The paper is organized as follows} : Section II is devoted to the computation of $D^2{\frak L}$. The proof of the pointwise estimate on the gradient of $u_k$ is the subject of section III. It is relying on several lemmae contained in the appendix (D to H) among which the weighted
Wente estimate (lemma~\ref{WenteWeight}). The section IV.1 is devoted to the introduction of the weight function $\om_{\eta,k}$ and the  proof of  (\ref{0.17}) which is the contribution of the necks to the positivity of the second derivative of ${\frak L}$ is based on the pointwise estimate of the gradient of $u_k$.
The conclusion of section IV.1 is given in lemma~\ref{lm-lower-bound-neck}. The use of the  {\em Sylvester inertia principle} and the diagonalisation of $D^2{\frak L}_{u_k}$ with respect to the weights $\om_{\eta,k}$ is given in subsection IV.2.  Thee proof of the main theorem~\ref{th-morse-stability} is given in section IV.3. In section V it is explained how the main arguments can be modified to cover the case of degenerating underlying Riemann surfaces and  to prove theorem~\ref{th-morse-stability-deg}.

\section{Preliminaries and Notations}
\reset
\subsection{The Hessian to the Dirichlet Energy : $\mathbf{D^2{\frak L}}$.}
The second derivative of the Dirichlet energy at harmonic maps has been computed for several decades already (see for instance \cite{Smi}). More generally, the second derivative of ${\frak{L}}$ at critical points can be found in \cite{ DHS} but for normal perturbations only. These computations  moreover are involving  the intrinsic Riemann tensors  of ${\mathcal{N}}^n$. We present below a general computation which involves the
extrinsic tensor given by the second fundamental form of ${\mathcal{N}}^n$ instead and which has the advantage of being relatively easily computable for a given sub-manifold. 

\medskip

We first consider the case $\om=0$ that is the case of the Dirichlet energy and harmonic maps.

\subsubsection{The harmonic map case.}

Let  $u$ be an harmonic map from $(\Sigma,h)$ into ${\mathcal{N}}^n$. We consider a smooth perturbation of $u$ of the form $u_t=\pi_{N}(u+t\,\phi)$ where $\phi$ is a smooth map in $C^\infty(\Sigma,{\R}^m)$
and $\pi_N$ is the orthogonal projection onto ${\mathcal{N}}^n$ (which is smooth in a tubular neighborhood of ${\mathcal{N}}^n$). We denote $$w:=\frac{du_t}{dt}(0)=d\pi_{u}\phi\ .$$Observe that
$d\pi_z$ for $z\in {\mathcal{N}}^n$ coincide with the matrix orthogonal projection $P_z$ of ${\R}^m$ onto $T_z{\mathcal{N}}^n$. We have moreover in ${\R}^m$ $d\pi_{u+t\,\phi}=d\pi_{\pi(u+t\,\phi)}$. This gives finally
\be
\label{II.1}
\frac{du}{dt}(t)=P_{u_t}\phi\ .
\ee
which gives
\be
\label{II.2}
\frac{d^2u}{dt^2}(0)=d_wP_{u}\phi\ .
\ee
We first compute the 
\be
\label{II.3}
\frac{d E(u_t)}{dt}=\int_\Sigma du_t\cdot_h d\frac{d u_t}{dt}\ dvol_h
\ee
The harmonicity of $u$ is equivalent to
\be
\label{II.4}
\forall\ \phi\in C^\infty(\Sigma,{\R}^m)\quad\quad 0=\frac{d E(u)}{dt}(0)=2\,\int_\Sigma du\cdot_h d(P_u\phi)\ dvol_h=2\,\int_\Sigma \Delta_hu\cdot_h P_u\phi\ dvol_h
\ee
Using the symmetry of the matrix $P_u$ we deduce that this is equivalent to 
\be
\label{II.5}
P_u\lf[\Delta_h u\rg]=0\quad\Longleftrightarrow \Delta_hu=- \, d[P_u]\cdot_h du
\ee
where we have used $P_udu=du$ and where we recall that $\Delta_h$ denotes the positive Laplace Beltrami operator. Combining (\ref{II.2}) and (\ref{II.6}) gives
\be
\label{II.6}
\begin{array}{l}
\ds\frac{d^2 E(u_t)}{dt^2}(0)=\int_\Sigma d\frac{ du_t}{dt}\cdot_h d\frac{d u_t}{dt}\ dvol_h+\int_\Sigma du_t\cdot_h d\frac{d^2 u_t}{dt^2}\ dvol_h\\[5mm]
\ds =\int_\Sigma dw\cdot_h dw\ dvol_h+\int_\Sigma \lf<\Delta_hu, d_wP_{u}\phi \rg>\ dvol_h\\[5mm]
\end{array}
\ee
where $<\cdot,\cdot>$ denotes the scalar multiplication in ${\R}^m$. Since $P_u[P_u\phi]=P_u\phi$ taking the differential along $w$ gives
\be
\label{II.7}
d_wP_u\phi=d_wP_uw+P_u d_wP_u\phi\ .
\ee
Since $P_u[\Delta_hu]=0$ we deduce
\be
\label{II.8}
 \lf<\Delta_hu, d_wP_{u}\phi \rg>=\lf<\Delta_hu  , d_wP_uw\rg>=-\lf<  d[P_u]\cdot_h du , d_wP_u\,w\rg>\ .
\ee
Combining (\ref{II.6}) and (\ref{II.8}) gives finally
\be
\label{II.9}
Q_u(w)=\int_\Sigma |dw|_h^2\ dvol_h-\int_\Sigma \lf<  d[P_u]\cdot_h du , d_wP_u\,w\rg>\ dvol_h
\ee
Recall that for any extensions in an ${\R}^n$ neighborhood of a point $z\in {\mathcal{N}}^n$ of a pair of vectors $X,Y$ in $T_z{\mathcal{N}}^n$ one has
\[
{\mathbb I}_z(X,Y)=d_YX-P_z(d_YX)
\]
Choosing an extension such that $P_z(X)=X$ for all  $z\in {\mathcal{N}}^n$ this gives
\be
\label{II.10-a}
{\mathbb I}_z(X,Y)=d_Y(P_z X)-P_z(d_YX)= d_YP_zX\ .
\ee
We finally obtain
\be
\label{II.9-a}
Q_u(w)=\int_\Sigma |dw|_h^2\ dvol_h-\int_\Sigma \lf<{\mathbb I}_u(du,du)_h   , {\mathbb I}_u(w,w)\rg>\ dvol_h
\ee
We shall denote by $S_u(du)_h$ the map from $\Sigma$ into the space of symmetric matrices real $m\times m$ matrices such that
\be
\label{II.9-b}
\forall x\in \Sigma\quad\forall X,Y\in T_{u(x)}{\mathcal{N}}^n\quad\lf<S_u(du)_h\,X,Y\rg>=\lf<{\mathbb I}_u(du,du)_h   , {\mathbb I}_u(X,Y)\rg>\ .
\ee
Observe that with this definition one has
\be
\label{II.9-c}
|S_u(du)_h|\le C\, |du|^2_h\ .
\ee
where $C>0$ is independent of $u$ and $h$.
\medskip

In the case of the sphere ${\mathcal{N}}^n=S^n$ for instance one has $P_u=I_{n+1}-u\otimes u$ which implies
\be
\label{II.10}
dP_u=-\lf[ du\otimes u-u\otimes du\rg]\quad \mbox{ and }\quad d_wP_u=-w\otimes u-u\otimes w\ .
\ee
This gives
\be
\label{II.11}
d[P_u]\cdot_h du=-\,u\,|du|^2_h\quad\mbox{ and }\quad d_wP_u\,w=- u\,|w|^2\ .
\ee
Combining (\ref{II.9}) and (\ref{II.10}) gives in this case
\be
\label{II.12}
Q_u(w)=\int_\Sigma |dw|_h^2\ dvol_h-\int_\Sigma|du|^2_h\,|w|^2\ dvol_h\ .
\ee
\subsubsection{The computation of $D^2{\frak L}$ in the general case.}

We extend $\al$ in a small neighborhood of ${\mathcal{N}}^n$ and we keep denoting $\al$ for $\pi_N^\ast\al$ in this neighborhood.\footnote{In this work we take the most current (but not exclusively adopted in the literature) convention
\[
dz^i\wedge dz^j(X,Y)=X^i\,Y^j-Y^i\,X^j\ 
\] and 
\[
 dz^k\wedge dz^i\wedge dz^j(X,Y,Z)=X^k\,(Y^i\,Z^j-Z^i\,Y^j)+Y^k\,(Z^i\,X^j-Z^j\,X^i)+Z^k\,(X^i\,Y^j-X^j\,Y^j)\ .
\]}
We write
\[
\al:=\frac{1}{2}\sum_{i,j=1}^m\al_{ij}(z)\ dz^i\wedge dz^j=\sum_{1\le i<j\le m}\ \al_{ij}(z)\ dz^i\wedge dz^j,
\]
in such a way that
\[
d\al=\frac{1}{2}\sum_{i,j,k=1}^m\p_{z_k}\al_{ij}\ dz_k\wedge dz_i\wedge dz_j \ .
\]
which gives 
\[
\begin{array}{l}
\ds d\al(X,Y,Z)=\frac{1}{2}\sum_{i,j,k=1}^m\p_{z_k}\al_{ij} \lf(  X^k\,(Y^i\,Z^j-Z^i\,Y^j)+Y^k\,(Z^i\,X^j-Z^j\,X^i)+Z^k\,(X^i\,Y^j-X^j\,Y^j)  \rg)\\[5mm]
\ds\ =\frac{1}{2} \sum_{k=1}^m X^k \sum_{i,j=1}^m [\p_{z_k}\al_{ij}+\p_{z_j}\al_{ki}+\p_{z_i}\al_{jk}] \ (Y^i\,Z^j-Z^i\,Y^j) 
\end{array}
\]
We introduce then
\[
\begin{array}{l}
\ds H(Y,Z):=\sum_{k=1}^mH^k(Y,Z)\,\p_{z_k}:=\frac{1}{2}\,\sum_{i,j,k=1}^m [\p_{z_k}\al_{ij}+\p_{z_j}\al_{ki}+\p_{z_i}\al_{jk}] \ (Y^i\,Z^j-Z^i\,Y^j)\ \p_{z_k} \\[5mm]
\ds\ =\sum_{k=1}^m\sum_{i<j}^m [\p_{z_k}\al_{ij}+\p_{z_j}\al_{ki}+\p_{z_i}\al_{jk}] \ (Y^i\,Z^j-Z^i\,Y^j)\ \p_{z_k}
\end{array}
\]
With this notation one has
\[
d\al(X,Y,Z)=\lf<X,H(Y,Z)\rg>\ . 
\]
We have
\be
\label{II.13}
\begin{array}{l}
\ds\frac{1}{2}\frac{d}{dt}\int _\Sigma u_t^\ast\al=\frac{1}{2}\frac{d}{dt}\int _\Sigma \sum_{i, j=1}^m\al_{ij}(u_t) du^i_t\wedge d u^j_t\\[5mm]
\ds=\frac{1}{2}\int _\Sigma \sum_{i, j,k=1}^m\p_{z_k}\al_{ij}(u_t) \frac{du^k_t}{dt}\ du^i_t\wedge d u^j_t+ \sum_{i, j=1}^m\al_{ij}(u_t) \ d\frac{d u^i_t}{dt}\wedge d u^j_t+ \sum_{i, j=1}^m\al_{ij}(u_t) \ d u^i_t\wedge d\frac{d u^j_t}{dt}\\[5mm]
\ds=\int_{\Sigma}\sum_{k=1}^m\frac{du^k_t}{dt}\sum_{i<j}[\p_{z_k}\al_{ij}(u_t)+\p_{z_j}\al_{ki}(u_t)+\p_{z_i}\al_{jk}(u_t)] \ du^i_t\wedge d u^j_t\\[5mm]
\ds=\int_{\Sigma}\sum_{k=1}^m\frac{du^k_t}{dt}\sum_{i<j} H^k_{ij}(u_t)\ du^i_t\wedge d u^j_t \end{array}
\ee
Combining (\ref{II.1}), (\ref{II.3}) and (\ref{II.13}) we obtain first that $u$ is a critical point of $\frak{L}$ if and only if
\be
\label{II.13-a}
\forall\ l=1\ldots m\quad 0=\sum_{k=1}^mP^{lk}_u\lf[\Delta_hu^k+\sum_{i<j} H^k_{ij}(u_t)\ du^i\wedge d u^j  \rg].
\ee

Since $P_u(H_u(du\wedge du))=H_u(du\wedge du)$ we deduce
\be
\label{II.13-b}
\Delta_hu+H_u(du\wedge du)=-\,d[P_u]\cdot_h du=-\,{\mathbb I}_u(du,du).
\ee
 
Taking the now second derivative of ${\frak L}-E$ gives
\be
\label{II.14}
\begin{array}{l}
\ds\frac{d^2({\frak L}-E)}{dt^2}=\frac{1}{2}\frac{d^2}{dt^2}\int _\Sigma u_t^\ast\al=\int_{\Sigma}\sum_{k=1}^m\frac{d^2u^k_t}{dt^2}\sum_{i<j} H^k_{ij}(u_t)\ du^i_t\wedge d u^j_t\\[5mm]
\ds+\int_{\Sigma}\sum_{k,l=1}^m\frac{du^k_t}{dt}\,\frac{du^l_t}{dt}\ \sum_{i<j} \p_{z_l}H^k_{ij}(u_t)\ du^i_t\wedge d u^j_t\\[5mm]
\ds+\int_{\Sigma}\sum_{k=1}^m\frac{du^k_t}{dt}\sum_{i<j} H^k_{ij}(u_t)\ \lf[d\frac{du^i_t}{dt}\wedge d u^j_t+du^i_t\wedge d \frac{ du^j_t}{dt}\rg]
\end{array}
\ee
Hence, for $u$ critical point of $\frak{L}$, we have
\be
\label{II.15}
\begin{array}{l}
\ds\frac{d^2{\frak L}(u)}{dt^2}(0)=\int_\Sigma |dw|^2_h\ dvol_h+\int_\Sigma\lf<\lf[\Delta_hu+H_u(du\wedge du)_h\rg], (I-P_u)d_wP_u\phi\rg>\ dvol_h\\[5mm]
\ds+\int_{\Sigma}\sum_{k,l=1}^m \,{w^k}\,{w^l}\ \sum_{i<j} \p_{z_l}H^k_{ij}(u)\ du^i\wedge d u^j+\int_{\Sigma}\sum_{k=1}^m w^k\,\sum_{i<j} H^k_{ij}(u)\ \lf[dw^i\wedge d u^j+du^i\wedge d w^j\rg]\\[5mm]
\end{array}
\ee
and using (\ref{II.7}) and (\ref{II.10-a}) we obtain
\be
\label{II.16-a}
\begin{array}{l}
\ds\frac{d^2{\frak L}(u)}{dt^2}(0)=\int_\Sigma |dw|^2_h\ dvol_h+\int_\Sigma\lf<\lf[\Delta_hu+H_u(du\wedge du)_h\rg], (I-P_u)d_wP_u\phi\rg>\ dvol_h\\[5mm]
\ds+\int_{\Sigma}\sum_{k,l=1}^m \,{w^k}\,{w^l}\ \sum_{i<j} \p_{z_l}H^k_{ij}(u)\ du^i\wedge d u^j+\int_{\Sigma}\sum_{k=1}^m w^k\,\sum_{i<j} H^k_{ij}(u)\ \lf[dw^i\wedge d u^j+du^i\wedge d w^j\rg]\\[5mm]
\end{array}
\ee
Hence, using (\ref{II.13-b}) this gives
\be
\label{II.16-b}
\begin{array}{l}
\ds\frac{d^2{\frak L}(u)}{dt^2}(0)=\int_\Sigma |dw|^2_h\ dvol_h-\int_\Sigma \lf<{\mathbb I}_u(du,du)_h   , {\mathbb I}_u(w,w)\rg>\ dvol_h\\[5mm]
\ds\quad+\int_\Sigma w\cdot \lf[2\, H_u(dw\wedge du)_h+\nabla_wH_u(du,du)_h\rg]\ dvol_h
\end{array}
\ee
where
\[
H_u(dw\wedge du)_h:=\frac{1}{2}\sum_{i,j=1}^mH^k_{ij}(u)\ \lf[dw^i\wedge d u^j+du^i\wedge d w^j\rg]\ \p_{z_k}\ .
\]
We introduce
\be
\label{II.18}
\begin{array}{rcl}
\ds{\mathcal H}_u(w)\ dvol_h:&=&\ds\sum_{i<j}^m H_{ij}(u)\ [dw^i\wedge d u^j+du^i\wedge dw^j]+\sum_{l=1}^m{w^l}\ \sum_{i<j}^m\p_{z_l}H_{ij}(u)\ du^i\wedge d u^j\\[5mm]
\ds\quad&=&\ds 2\, H_u(dw\wedge du)_h+\nabla_wH_(du\wedge du)_h
\end{array}
\ee
and we claim that ${\mathcal H}_u$ is self-adjoint on $V_u=\Gamma(u^{-1} TN)$.
We have
\be
\label{II.16}
\begin{array}{l}
\ds\int_{\Sigma}\sum_{k=1}^m w^k\,\sum_{i<j} H^k_{ij}(u)\ \lf[dv^i\wedge d u^j+du^i\wedge d v^j\rg]\\[5mm]
\ds=-\frac{1}{2}\int_{\Sigma}\sum_{k=1}^m v^i \sum_{i,j}^m H^k_{ij}(u)\ dw^k\wedge d u^j+\frac{1}{2}\int_\Sigma\sum_{i,j}^m v^j  \sum_{k=1}^m H^k_{ij}(u)\ dw^k\wedge  du^i\\[5mm]
\ds-\frac{1}{2}\int_{\Sigma}\sum_{k,l=1}^m v^i \sum_{i,j}^m \p_{z_l}H^k_{ij}(u)\ w^k\ du^l\wedge d u^j+\frac{1}{2}\int_\Sigma\sum_{i,j}^m v^j  \sum_{k,l=1}^m \p_{z_l}H^k_{ij}(u)\ w^k du^l\wedge  du^i\\[5mm]
\ds=\int_{\Sigma}\sum_{k=1}^m v^k \sum_{i<j}^m H^k_{ij}(u)\ [dw^i\wedge d u^j+du^i\wedge dw^j]\\[5mm]
\ds+\frac{1}{2}\int_{\Sigma}\sum_{k,l=1}^m v^k \sum_{i,j=1}^m \p_{z_l}H^k_{ij}(u)\ w^i\ du^l\wedge d u^j-\frac{1}{2}\int_\Sigma\sum_{i,j=1}^m v^k  \sum_{k,l=1}^m \p_{z_l}H^k_{ij}(u)\ w^j du^l\wedge  du^i
\end{array}
\ee
We have
\be
\label{II.17}
\begin{array}{l}
\ds\frac{1}{2}\int_{\Sigma}\sum_{k,l=1}^m \,{w^k}\,{v^l}\ \sum_{i,j=1}^m\p_{z_l}H^k_{ij}(u)\ du^i\wedge d u^j+\frac{1}{2}\int_{\Sigma}\sum_{k,l=1}^m v^k \sum_{i,j=1}^m \p_{z_l}H^k_{ij}(u)\ w^i\ du^l\wedge d u^j\\[5mm]
\ds-\frac{1}{2}\int_\Sigma\sum_{i,j=1}^m v^k  \sum_{k,l=1}^m \p_{z_l}H^k_{ij}(u)\ w^j du^l\wedge  du^i\\[5mm]
\ds=\frac{1}{2}\int_{\Sigma}\sum_{k,l=1}^m \,{w^k}\,{v^l}\ \sum_{i,j=1}^m\lf[\p_{z_l}H^k_{ij}(u)+\p_{z_i}H^l_{kj}(u)-  \p_{z_j} H^l_{ki} -\p_{z_k}H^l_{ij}\rg]\ du^i\wedge d u^j\\[5mm]
\ds+\frac{1}{2}\int_{\Sigma}\sum_{k,l=1}^m \,v^k\,{w^l}\ \sum_{i,j=1}^m\p_{z_l}H^k_{ij}(u)\ du^i\wedge d u^j\\[5mm]
\ds = \frac{1}{2}\int_{\Sigma}\sum_{k,l=1}^m \,v^k\,{w^l}\ \sum_{i,j=1}^m\p_{z_l}H^k_{ij}(u)\ du^i\wedge d u^j
\end{array}
\ee
Combining (\ref{II.16}) and (\ref{II.17}) gives that the operator ${\mathcal H}_u$ is self-adjoint on $V_u=\Gamma(u^{-1} T{\mathcal{N}}^n)$. 

\medskip

We shall denote
\[
{\mathcal L}_uw:=P_u\lf[\Delta_hw-S_u(du)_h\, w+{\mathcal H}_u(w)\rg]\ .
\]
and with this notation
\[
D^2{\frak L}_u(w)=\int_\Sigma w\cdot{\mathcal L}_u w\ dvol_h\ .
\]

\section{Proof of the pointwise estimate of the gradient in the necks.}
 \reset
Let ${\mathcal{N}}^n$ be an arbitrary closed $C^2$ sub-manifold of ${\R}^m$ and $\al$ be an arbitrary $C^2$ 2-form on ${\mathcal{N}}^n$. For a sake of clarity we assume that we have exactly one bubble. The general case is inducing more complicated notations but no change at all in the arguments. 
 
 Hence from now on we shall consider a sequence of critical points  $u_k$ of ${\frak L}$ among maps from $(\Sigma,h)$ into ${\mathcal{N}}^n$ which ``bubble tree converges'' to a pair of  maps $(u_\infty, v_\infty)$ respectively from $\Sigma$ into ${\mathcal{N}}^n$ and from $\C$ into ${\mathcal{N}}^n$. More precisely there exists $p\in \Sigma$, $x_k\rightarrow p$ and $\delta_k\rightarrow 0$ such that 
\begin{itemize}
\item[i)]
\be
\label{I.8}
u_k\longrightarrow u_\infty\quad\mbox{in } C^1_{loc}(\Sigma\setminus\{p\})\ ,
\ee
\item[ii)] in local fixed conformal coordinates for $h$ around $p$
\be
\label{I.9}
u_k(\delta_k\, y+x_k)\longrightarrow v_\infty(y)\quad\mbox{in } C^1_{loc}({\C})\ ,
\ee
\item[iii)]
\be
\label{I.10}
\lim_{\eta\rightarrow 0}\lim_{k\rightarrow +\infty}\ \sup_{\delta_k/\eta<\rho< 2\rho<\eta} \int_{B_{2\rho}(x_k)\setminus B_\rho(x_k)}|du_k|^2_h\ dvol_h=0
\ee
\end{itemize}
The main arguments in \cite{LaRi1} is implying
\be
\label{I.11}
\lim_{\eta\rightarrow 0}\limsup_{k\rightarrow +\infty}\  \int_{A(\eta,\delta_k)}|du_k|^2_h\ dvol_h=0
\ee
where we shall be using the following notation in the rest of the paper
\[
A(\eta,\delta_k):=B_{\eta}(x_k)\setminus B_{\delta_k/\eta}(x_k)\ .
\]
The annulus $A(\eta,\delta_k)$ is  called {\it annular neck region} and the identity (\ref{I.11}) is also known under $L^2-${\it energy quantization}.
The goal of the present section is to convert the $L^2$-{\it energy quantization} given by  (\ref{I.11}) into a pointwise quantization given by the following lemma.
\begin{Lm}
\label{lm-ptwz}
Let $u_k$ be a sequence of critical points of $\frak{L}$ from $(\Sigma,h)$ into ${\mathcal{N}}^n$, a $C^2$ sub-manifold of ${\R}^m$ and $\al$ is a $C^2$ 2-form of ${\mathcal{N}}^n$. Assume $u_k$  ``bubble tree converges''
to a pair of  maps  $(u_\infty, v_\infty)$ in the sense that (\ref{I.8}), (\ref{I.9}) and (\ref{I.10}) are satisfied. Then there exist  $C>0$ and $0<\beta<1$ such that, for any $\eta$ small enough independent of $k$  in local conformal coordinates around $x_k$ such that $(x_1(x_k),x_2(x_k))=(0,0)$ there holds
\be
\label{I.12}
\begin{array}{l}
\ds\forall x\in A(\eta,\delta_k)\quad\quad
\ds |x|^2\ |\nabla u_k|^2(x)\le C\ \lf[\frac{|x|^\beta}{\eta^\beta}+\lf(\frac{\delta_k}{\eta\,|x|}\rg)^\beta\rg]\  \int_{A(2\eta,\delta_k)}|du_k|^2_h\ dvol_h+C_{\eta,k}\ ,
\end{array}
\ee
where
\be
\label{I.13}
\lim_{\eta\rightarrow 0}\limsup_{k\rightarrow +\infty}\ \log^2\lf( \frac{\eta^2}{\delta_k}  \rg)\, C_{\eta,k}=0\ .
\ee
and $C_{\eta,k}$ is uniformly bounded independent of $k$ and of $\eta<\eta_0$ for some $\eta_0>0$. \hfill $\Box$
\end{Lm}
\subsection{The $C^0$-no neck Property and $L^{2,1}-$Energy Quantization.}

The $L^2$ energy quantization (\ref{I.11})  is reinforced in \cite{LaRi1} - estimate (56) - by the following
\be
\label{V.3-a}
\lim_{\eta\rightarrow 0}\lim_{k\rightarrow +\infty}\ \lf\|\frac{1}{r}\frac{\p u_k}{\p\theta}\rg\|_{L^{2,1}(A(\eta,\delta_k))}=0\ .
\ee
Where we recall that the $L^{2,1}$ norm of a function $f$ is equivalent to the following quasi-norm
\[
|f|_{L^{2,1}(A(\eta,\delta_k))}=\int_0^{+\infty}\lf|\lf\{ x\in A(\eta,\delta_k)\ ;\ |f(x)|>s\rg\}\rg|^{1/2}\ ds\ .
\]
Recall that, for any critical point of conformally invariant problems of the form given by $\frak{L}$, the Hopf differentials of $u_k$ given in local conformal coordinates for $(\Sigma,h)$ by
\be
\label{V.3-b}
\frak{H}(u_k):=\p_zu_k\cdot\p_z u_k\ dz\otimes dz
\ee
defines an holomorphic quadratic differential of $(\Sigma,h)$. The space of holomorphic quadratic differentials is finite dimensional and since $\frak{H}(u_k)$ is uniformly bounded in $L^1$, it is pre-compact in any norm.
In particular it is uniformly bounded in $L^2$. We have
\[
\p_{z}u_k=e^{-i\theta}\,(\p_\rho u_k-i\,\rho^{-1}\,\p_\theta u_k) \ .
\]
Hence
\be
\label{V.3-c}
\Re \lf(e^{2\,i\,\theta}\, \p_zu_k\cdot\p_z u_k\rg):=|\p_\rho u_k|^2-\rho^{-2}\,|\p_\theta u_k|^2\ .
\ee
Which gives
\be
\label{V.3-d}
|\p_\rho u_k|\le \, \rho^{-1}\,|\p_\theta u_k|+\,\sqrt{\lf|\Re \lf(e^{2\,i\,\theta}\, \p_zu_k\cdot\p_z u_k\rg)\rg|}\ .
\ee
Generalized H\"older inequality (see \cite{Gra1}) gives
\be
\label{V.3-e}
\begin{array}{rcl}
\ds\lf\|\sqrt{\lf|\Re \lf(e^{2\,i\,\theta}\, \p_zu_k\cdot\p_z u_k\rg)\rg|}\rg\|_{L^{2,1}(A(\eta,\delta_k))}&\le& |A(\eta,\delta_k))|^{1/4}\ \lf\|\sqrt{\lf|\Re \lf(e^{2\,i\,\theta}\, \p_zu_k\cdot\p_z u_k\rg)\rg|}\rg\|_{L^{4}(A(\eta,\delta_k))}\\[5mm]
\ds\quad &\le& \sqrt{\eta}\, \sqrt{\|\frak{H}(u_k)\|_{L^2(\Sigma)}}
\end{array}
\ee
Combining (\ref{V.3-a}) and (\ref{V.3-e}) gives
\be
\label{V.3-f}
\lim_{\eta\rightarrow 0}\limsup_{k\rightarrow +\infty}\ \lf\|\frac{\p u_k}{\p r}\rg\|_{L^{2,1}(A(\eta,\delta_k))}=0\ .
\ee
Hence we have established the following lemma.
\begin{Lm}
\label{lm-L-2-1quanti-general}
Let $u_k$ be a sequence of critical points of $\frak{L}$ satisfying (\ref{I.8}), (\ref{I.9}) and (\ref{I.10}). Then the following holds
\be
\label{V.3-g}
\lim_{\eta\rightarrow 0}\ \limsup_{k\rightarrow+\infty}\ \|du_k\|_{L^{2,1}(A(\eta,\delta_k))}=0\ .
\ee
\hfill $\Box$
\end{Lm}
\subsection{Proof of Lemma~\ref{lm-ptwz}.}
As explained in \cite{Riv1} the Euler-Lagrange equation (\ref{0.1}) can be rewritten in the form
\be
\label{V.1}
\Delta_hu_k=\Om_k\cdot du_k\quad\quad\mbox{ on }\Sigma\ ,
\ee
where there exists $C>0$ such that
\be
\label{V.2}
\Om_k\in \Gamma(T^\ast\Sigma\otimes so(m))\quad\mbox{ and }\quad |\Om_k|(x)\le C\,|du_k|_h(x)\quad\forall x\in \Sigma\ .
\ee
where $so(m)$ denotes the space of real $m\times m$ antisymmetric matrices. 

\medskip

By extending $\Om_k$ outside $A(2\eta,\delta_k)$ by 0  one produces  $\ti{\Om}_k\in L^2(\wedge^1B_1(0),so(m))$.  For $\eta$ small enough and $k$ large enough, using theorem I.4
of \cite{Riv1} we obtain the existence of $A_k\in W^{1,2}(B_1,Gl_m({\R}))$ and $B_k\in W^{1,2}(B_1,M_m({\R}))$ such that
\be
\label{V.4}
\nabla A_k-A_k\,\ti{\Om}_k=\nabla^\perp B_k\quad\quad\mbox{ in }B_1(0)\ ,
\ee
and
\be
\label{V.5}
\begin{array}{l}
\ds\int_{B_1}|\nabla A_k|^2+|\nabla A_k^{-1}|^2\ dx^2+\int_{B_1}|\nabla B_k|^2\ dx^2+\lf\|\mbox{dist}\lf(\lf\{A_k,A_k^{-1}\rg\},SO(m)\rg)\rg\|^2_{L^\infty(B_1(0))}\\[5mm]
\ds \le C\ \int_{A(2\eta,\delta_k)}|\Om_k|^2\ dx^2\le C\ \int_{A(2\eta,\delta_k)}|\nabla u_k|^2\ dx^2\ .
\end{array}
\ee
and we have also obviously
\be
\label{V.6}
\int_{B_1}|\nabla B_k|^2\ dx^2\le C\ \int_{A(2\eta,\delta_k)}|\nabla u_k|^2\ dx^2\ .
\ee
It is proved in \cite{Riv1} that the following conservation law (in conformal coordinates for $h$)  is satisfied
\be
\label{V.7}
\mbox{div}\lf( A_k\,\nabla u_k \rg)=\nabla^\perp B_k\cdot\nabla u_k\quad\mbox{ in }A(\eta,\delta_k)\ .
\ee
We extend  $u_k$, $A_k$ and $B_k$ outside $A(\eta,\delta_k)$ in the whole ${\C}$ using lemma~\ref{lma-extension} and we denote respectively $\ti{A}_k$ and $\ti{B}_k$ these extensions.
Let $\varphi_k$ and $\psi_k$ be the solutions in $\dot{W}^{1,2}({\C})$ respectively of 
\be
\label{V.8}
\Delta\varphi_k=\nabla^\perp \ti{B}_k\cdot\nabla \ti{u}_k\quad\mbox{ in }{\C}\ ,
\ee
and
\be
\label{V.9}
\Delta\psi_k=\nabla^\perp \ti{A}_k\cdot\nabla \ti{u}_k\quad\mbox{ in }{\C}\ .
\ee
We have that $A_k du_k-d\varphi_k-\ast d\psi_k$  realizes an harmonic 1-form in $A(\eta,\delta_k)$. Hence there exists $\frak{h}_k$ harmonic such that
\be
\label{V.10}
A_k\,\nabla u_k-\nabla\varphi_k-\nabla^\perp\psi_k=\nabla \frak{h}_k+ C^1_k \,\nabla^\perp\log r
\ee
We have in particular for any $r\in [\delta_k/\eta,\eta]$ such that
\be
\label{V.11}
 2\pi\, C_k^1 =\int_{\p B_r(0)} A_k\,\frac{1}{r}\frac{\p u_k}{\p\theta}\ \,d\theta-\int_{\p B_r(0)}\p_\nu\psi_k\ \,d\theta=\int_{B_r}\nabla \ti{A}_k\cdot\nabla^\perp u_k-\Delta\psi_k\ dx^2=0
\ee
Hence finally we have
\be
\label{V.12}
A_k\,\nabla u_k=\nabla\varphi_k+\nabla^\perp\psi_k+\nabla \frak{h}_k\ .
\ee
The harmonic part $\frak{h}_k$ is decomposed as follows
\[
\frak{h}_k=\frak{h}_k^++\frak{h}_k^-+\frak{h}_k^0\ \ \mbox{ where }\ \ \frak{h}_k^+=\Re\sum_{n>0} h_{n,k}z^n\ \ ,\ \ \frak{h}_k^-=\Re\sum_{n<0} h_{n,k}z^n\ \ \mbox{ and }\ \ \frak{h}_k^0=\frak{h}_{0,k}+C_{\eta,k}^0\,\log|z|\ .
\]
 Coifman-Lions-Meyer-Semmes estimate (see \cite{Hel}) applies to both $\varphi_k$ and $\psi_k$ and says respectively
\be
\label{V.13}
\|\nabla \varphi_k\|^2_{L^{2,1}({\C})}\le C\ \int_{\C}|\nabla\ti{B}_k|^2\ dx^2\ \int_{\C}|\nabla\ti{u}_k|^2\ dx^2\le C\,\lf[\int_{A(\eta,\delta_k)}|d u_k|^2_h\ dvol_h\rg]^2\ ,
\ee
and
\be
\label{V.14}
\|\nabla \psi_k\|^2_{L^{2,1}({\C})}\le C\ \int_{\C}|\nabla\ti{A}_k|^2\ dx^2\ \int_{\C}|\nabla\ti{u}_k|^2\ dx^2\le C\,\lf[\int_{A(\eta,\delta_k)}|d u_k|^2_h\ dvol_h\rg]^2\ ,
\ee
where we have used (\ref{V.5}) and (\ref{V.6}). Hence we deduce from (\ref{V.3-g})
\be
\label{V.15}
\lim_{\eta\rightarrow 0}\ \limsup_{k\rightarrow+\infty}\ \|d\frak{h}_k\|_{L^{2,1}(A(\eta,\delta_k))}=0\ .
\ee
The following lemma due to A.Michelat and the third author \cite{MiRi} is going to be used for controlling the $L^{2,1}-$norms of $d\frak{h}_k^\pm$ in the neck region.

 \begin{Lm}[\cite{MiRi} Lemma 2.3]\label{MiRi}
 Let $0<4\,r<R$ be fixed radii, $u\colon B_R(0)\setminus \bar{B}_r(0)\to \R$ be a harmonic function such that for some $\rho_0\in (r,R)$
\begin{equation}\label{dernul}
 \int_{\partial B_{\rho_0}(0)}\partial_{\nu} u=0.
\end{equation}
Then for all $\left(\frac{r}{R}\right)^{1/2}<t<1$ we have
\begin{equation}
\|\nabla u\|_{L^{2,1}(B_{t R}(0)\setminus B_{t^{-1} r}(0))}\le C\, \frac{t}{1-t}\|\nabla u\|_{L^{2}(B_R(0)\setminus B_r(0))}\ ,
\end{equation}
for some universal constant $C$ independent of $R$, $r$ and $\al$. \hfill $\Box$
\end{Lm}
We apply lemma~\ref{MiRi} to $\frak{h}^\pm_k$ and this gives for any fixed $t<1$
\be
\label{V.16}
 \|d\frak{h}^\pm_k\|_{L^{2,1}(A(t\,\eta,\delta_k))}\le C\,  \|d\frak{h}^\pm_k\|_{L^{2}(A(\eta,\delta_k))}\ .
\ee
We have
\be
\label{V.16-aa}
\begin{array}{l}
\ds\int_{A(\eta,\delta_k)}|\nabla \frak{h}_k^0|^2\ dx^2+\int_{A(\eta,\delta_k)}|\nabla \frak{h}_k^+|^2\ dx^2+\int_{A(\eta,\delta_k)}|\nabla \frak{h}_k^-|^2\ dx^2\\[5mm]
\ds\le 4\,\int_{A(\eta,\delta_k)}|\nabla u_k|^2\ dx^2+4\,\int_{A(\eta,\delta_k)}|\nabla \varphi_k|^2\ dx^2+4\,\int_{A(\eta,\delta_k)}|\nabla \psi_k|^2\ dx^2\ .
\end{array}
\ee
Combining (\ref{I.11}), (\ref{V.13}), (\ref{V.14}), (\ref{V.15}), \eqref{V.16} and  (\ref{V.16-aa}) we obtain
\be
\label{V.16-a}
\lim_{\eta\rightarrow 0}\ \lim_{k\rightarrow+\infty}\ \|d\frak{h}^0_k\|_{L^{2,1}(A(t\,\eta,\delta_k))}=0\ .
\ee
The following direct  computations are illustrating the ``screening'' each of the $L^{2,q}$ norms for $q=1,2,\infty$ is operating on the neck.
It looks elementary but this is nevertheless central in the present work.
\be
\label{lorentz}
\lf\{
\begin{array}{l}
\ds \|\nabla\log |x|\|_{L^{2,\infty}(A(\eta,\delta_k))}\simeq\ds\sup_{\lambda>0} \la\,\lf|\lf\{x\in A(\eta,\delta_k): |x|^{-1}>\lambda\rg\}\rg|^{1/2} =\sqrt{\pi}\\[5mm]
\ds\|\nabla\log |x|\|_{L^{2}(A(\eta,\delta_k))}=\lf[\int_{A(\eta,\delta_k)}|\nabla \log|x||^2 dx^2   \rg]^{1/2}=\sqrt{2\pi}\ \sqrt{\log\lf(\frac{\eta^2}{\delta_k}\rg)}\\[5mm]
\ds\|\nabla\log |x|\|_{L^{2,1}(A(\eta,\delta_k))}\simeq\ds\int_{0}^{\infty}|x\in A(\eta,\delta_k):~~|x|^{-1}\ge \lambda|^{1/2}d\lambda=\sqrt{2\pi}\,\log\lf[\frac{\eta^2}{\delta_k}\rg]\ .
 \end{array}
 \rg.
 \ee
where the sign $\simeq$ is referring to the fact that we are giving equivalent quasi-norms to the norms $L^{2,\infty}$ and $L^{2,1}$ respectively and the constants behind are universal.

\medskip

We deduce from (\ref{V.16-a}) and (\ref{lorentz}) the following lemma\footnote{In the sphere case one has that $A_k$ is the identity matrix, $\psi_k=0$ and moreover the  conservation law
\[
\mbox{div}(\nabla u_k+B_k\,\nabla^\perp u_k)=0
\] is satisfied throughout the bubble in the whole domain $\Sigma$. As a consequence, in this very special case, one can deduce by an integration by  part that
\be
\label{IV.16-aa}
C_{\eta,k}^0=0\ .
\ee This might not be the case in the general situation where (\ref{V.7}) cannot a-priori be extended throughout $B_{\delta_k/\eta}(x_k)$.}.
\begin{Lm}
\label{lm-estim-c0k}
Under the above notations one has
\be
\label{IV.13}
C_{\eta,k}^0=o_{{\eta},k}\lf( \frac{1}{\log\lf(\frac{\eta^2}{\delta_k}\rg)} \rg), ~~\mbox{as $k\to +\infty,~\eta\to 0$}.
\ee
\hfill $\Box$
\end{Lm}

\medskip
We introduce the following notation
\be
\label{IV.17}
s_1:=[\log\eta^{-1}]\quad,\quad s_2=[\log \lf(\eta\,\delta_k^{-1}\rg)]\quad\mbox{ and }\quad A_k=B_{2^{-k}}(0)\setminus B_{2^{-k-1}}(0)\ .
\ee
We adopt the notations of the previous subsection. 
 \begin{Prop}\label{morreyu}
 There exists $0<\gamma<1$ such that for any $\gamma<\mu<1$ there exists $C>0$ such that for $\eta$ sufficiently small and $k$ large enough the following holds.  For all $j \in [s_1,s_2]$ we have 
 \be\label{morreyu2}
 \begin{array}{l}
 \ds\int_{A_j}|\nabla  u_k|^2\ dx^2
 \le C  \sum_{\ell=s_1}^{s_2}\mu^{|\ell-j|}{\int_{ A_{\ell}}}|\nabla \frak{h}_k|^2+C\,\sum_{\ell=s_1}^{s_2}\mu^{|\ell-j|}\gamma^{\ell}\left(\int_{ B_\eta(0)}|\nabla \tilde u_k|^2\right)^2\\[5mm]
\ds +\,C\ \int_{B_\eta(0)}|\nabla \tilde u_k|^2\ dx^2\,\left(\mu^{j-s_1}\int_{B_{2\eta}(0)\setminus B_\eta(0)}|\nabla \tilde u_k|^2\ dx^2+\mu^{s_2-j}\int_{B_{\frac{\delta_k}{\eta}}(0)\setminus B_{\frac{\delta_k}{2\,\eta}}(0)}|\nabla \tilde{u}_k|^2\ dx^2\right)
\end{array}
 \ee
\end{Prop}
\noindent{\bf Proof of Proposition \ref{morreyu}.} We omit to explicitly write the subscript $k$ during the proof of the proposition. We can apply lemma~\ref{morrey9bis}  to both $\varphi$ and $\psi$ in order to deduce for any $j\in\{s_1,\ldots, s_2\}$
 \begin{eqnarray}\label{V.18}
{\int_{ A_{j}}}|\nabla \varphi|^2\ dx^2 
&\le& \gamma^{j}\,\int_{ B_\eta(0)}|\nabla \varphi|^2\ dx^2+C\,\int_{B_\eta(0)}|\nabla \tilde{B}|^2\ dx^2\ 
\sum_{n=0}^{\infty}\gamma^{|n-j|}\int_{A_n}|\nabla \tilde u|^2\ dx^2\ ,\nonumber\end{eqnarray}
and
 \begin{eqnarray}\label{V.19}
{\int_{ A_{j}}}|\nabla \psi|^2\ dx^2 
&\le& \gamma^{j}\,\int_{ B_\eta(0)}|\nabla \psi|^2\ dx^2+C\,\int_{B_\eta(0)}|\nabla \tilde{A}|^2\ dx^2\ 
\sum_{n=0}^{\infty}\gamma^{|n-j|}\int_{A_n}|\nabla \tilde u|^2\ dx^2\ .\nonumber\end{eqnarray}
Using the fact that 
\be
\label{V.20}
\int_{B_\eta(0)}|\nabla \tilde{A}_k|^2+|\nabla \ti{B}_k|^2\ dx^2\le C\, \int_{A(\eta,\delta)}|\nabla u_k|^2\ dx^2
\ee
we obtain respectively
 \begin{eqnarray}\label{V.18-a}
{\int_{ A_{j}}}|\nabla \varphi|^2\ dx^2 
&\le& \gamma^{j}\,\int_{ B_\eta(0)}|\nabla \varphi|^2\ dx^2+C\,\int_{A(\eta,\delta)}|\nabla u|^2\ dx^2\ 
\sum_{n=0}^{\infty}\gamma^{|n-j|}\int_{A_n}|\nabla \tilde u|^2\ dx^2\ ,\nonumber\end{eqnarray}
and
 \begin{eqnarray}\label{V.19-a}
{\int_{ A_{j}}}|\nabla \psi|^2\ dx^2 
&\le& \gamma^{j}\,\int_{ B_\eta(0)}|\nabla \psi|^2\ dx^2+C\,\int_{A(\eta,\delta)}|\nabla u|^2\ dx^2\ 
\sum_{n=0}^{\infty}\gamma^{|n-j|}\int_{A_n}|\nabla \tilde u|^2\ dx^2\ .\nonumber\end{eqnarray}
By using that $A\,\nabla u=\nabla\varphi+\nabla^\perp\psi+\nabla \frak{h}$ it also holds:
\begin{eqnarray}\label{morrey17}
{\int_{ A_{j}}}|\nabla \tilde u|^2\ dx^2&\le& 4 {\int_{ A_{j}}}|\nabla \frak{h}|^2\ dx^2+4\,{\int_{ A_{j}}}|\nabla \varphi|^2+|\nabla\psi|^2\ dx^2\nonumber\\
&\le&4\, {\int_{ A_{j}}}|\nabla \frak{h}|^2\ dx^2+4\,\gamma^{j}\left(\int_{ B(0,\eta)}|\nabla \tilde u|^2\ dx^2\right)^2\\
&+& C\int_{B(0,\eta)}|\nabla \tilde u|^2\ dx^2\ 
\sum_{n=0}^{\infty}\gamma^{|n-j|}\int_{A_n}|\nabla \tilde u|^2\ dx^2\ .\nonumber\end{eqnarray}
From Lemma \ref{lma-extension} implies $$\int_{B^c_{2\eta}(0)}|\nabla \tilde u|^2\ dx^2=0\quad\mbox{  and }\quad\int_{B_{\frac{\delta}{2\eta}}(0)}|\nabla \tilde u|^2\ dx^2=0\ .$$
If  $\int_{B(0,\eta)}|\nabla \tilde u|^2<\varepsilon_0$  then  by  applying lemma \ref{morrey11} to the sequences   
$$a_k=\int_{A_k}|\nabla \tilde u|^2\ dx^2~~~\mbox{ and}~~~
 b_k=2 {\int_{ A_{k}}}|\nabla \frak{h}|^2\ dx^2+2\,\gamma^{k}\left(\int_{ B(0,\eta)}|\nabla \tilde u|^2\ dx^2\right)^2$$
and using the fact    that $a_k=0$ if $k\ge [\log(2\eta/\delta)]$ and $k\le  [\log((2\eta)^{-1})]$
we obtain for $\eta$ small enough and $k$ large enough, passing $$C\  \int_{B_\eta(0)}|\nabla \tilde u|^2\ dx^2\, \sum_{\ell=s_1}^{s_2}\mu^{|\ell-j|}\int_{A_\ell}|\nabla \tilde u|^2\ dx^2$$ from the r.h.s to the l.h.s.
 of (\ref{morrey12})
\be\label{morrey122}
\begin{array}{l}
\ds\sum_{\ell=s_1}^{s_2}\mu^{|\ell-j|}\int_{A_\ell}|\nabla \tilde u|^2\ dx^2
 \le\, \sum_{\ell=s_1}^{s_2}\mu^{|\ell-j|}b_{\ell}+\,C\  \int_{B_\eta(0)}|\nabla \tilde u|^2\ dx^2 \sum_{n=[\log\left(\frac{1}{2\eta}\right)]}^{s_1-1} \mu^{j-n}\,\int_{A_n}|\nabla \tilde u|^2\ dx^2 \\
\ds+ \,C\ \int_{B_\eta(0)}|\nabla \tilde u|^2\ \sum_{n=s_2+1}^{[\log(2\eta/\delta)]} \mu^{n-j}\,\int_{A_n}|\nabla \tilde u|^2\\
\ds \le\, C\  \sum_{\ell=s_1}^{s_2}\mu^{|\ell-j|}{\int_{ A_{\ell}}}|\nabla \frak{h}|^2+C\sum_{\ell=s_1}^{s_2}\mu^{|\ell-j|}\gamma^{\ell}\left(\int_{ B(0,\eta)}|\nabla \tilde u|^2\right)^2\\\
\ds +\,C\, \left(\int_{B(0,\eta)}|\nabla \tilde u|^2 \right)\left(\mu^{j-s_1}\int_{B(0,2\eta)\setminus B(0,\eta)}|\nabla \tilde u|^2+\mu^{s_2-j}\int_{B(0,\delta/\eta)\setminus B(0,\delta/2\eta)}|\nabla \tilde u|^2\right)\ .
\end{array}
 \ee
In particular we get \eqref{morreyu2} and we conclude the proof of proposition~\ref{morreyu}.\hfill       $\Box$

\medskip

Next we estimate the Dirichlet energy of the harmonic part of $u_k$.
\begin{Lm}\label{estseriesw}
 There exists $4^{-1}<\mu<1$ and $C>0$ independent of $k$ and $\eta$ such that for any $j\in \{s_1,\ldots,s_2\}$:
 \begin{eqnarray}\label{estseriesw2}
   \sum_{\ell=s_1}^{s_2}\mu^{|\ell-j|}{\int_{ A_{\ell}}}|\nabla \frak{h}_k|^2\ dx^2
   &\le&C\left[\left(\frac{2^{-j}}{\eta}\right)^\beta+\left(\frac{\delta}{ 2^{-j}\,\eta}\right)^\beta\right]\int_{A(2\,\eta,\delta_k)}|\nabla \tilde u_k|^2\ dx^2\\
   &+& o_{\eta,k}\left(\frac{1}{\log^2\left(\frac{\eta^2}{\delta_k}\right)}\right),\nonumber\
  \end{eqnarray}
  where $\lim_{\eta\rightarrow 0}\limsup_{k\rightarrow+\infty}o_{\eta,k}(1)=0$ and $\beta:=-\log_2\mu$.\hfill $\Box$
 \end{Lm}
\noindent{\bf Proof of Lemma \ref{estseriesw}.} In this proof again we omit to write the subscript $k$.
 \par
 {\bf 1.} Using (\ref{A-16}) from lemma~\ref{ptw-harm} on $A(2\,\eta,\delta)$, $\frak{h}^+$ satisfies for $\rho\in [2^{-s_2}, 2^{-s_1}]$
 \begin{align} 
    |\nabla \frak{h}^+(\rho,\theta)|^2 
    &\leq \frac{C}{\eta^2}\int_{A(\eta,\delta)}|\nabla \frak{h}^+|^2\ dx^2\le \frac{C}{\eta^2}\ \int_{A(2\,\eta,\delta)}|\nabla \tilde u|^2\ dx^2\ .
\end{align}
Therefore for every $\ell\in \{s_1,\ldots,s_2\}$ we have
$$\int_{A_{\ell}}|\nabla \frak{h}^+(\rho,\theta)|^2\ \rho\,d\rho\,d\theta\le 2^{-2\ell}\frac{C}{\eta^2}\int_{A(2\,\eta,\delta)}|\nabla \tilde u|^2\ dx^2$$
and
\begin{eqnarray}\label{estseriesw+}
\sum_{\ell=s_1}^{s_2}\mu^{|\ell-j|}{\int_{ A_{\ell}}}|\nabla \frak{h}^+|^2 &\le& \frac{C}{\eta^2}\ \left[\sum_{\ell=s_1}^{s_2}\mu^{|\ell-j|}2^{-2\ell}\right]\ \int_{A(2\eta,\delta)}|\nabla \tilde u|^2\ dx^2\nonumber\\
&=& \frac{C}{\eta^2}\ \left[\sum_{\ell=s_1}^{j}\mu^{j-\ell}2^{-2\ell}+\sum_{\ell=j+1}^{s_2}\mu^{\ell-j}2^{-2\ell}\right]\ \int_{A(2\eta,\delta)}|\nabla \tilde u|^2\ dx^2\nonumber\\
&\le& \frac{C}{\eta^2}\ \lf(\mu^{j-s_1}\ 2^{-2s_1}+2^{-2j}\rg)\ \int_{A(2\eta,\delta)}|\nabla \tilde u|^2\ dx^2.
\end{eqnarray}
and we have
\be
\label{exp-1}
\mu^{j-s_1}\ 2^{-2s_1}\le\eta^2\,\mu^{j-s_1}=\eta^2\, 2^{\log_2\mu (j-s_1)}=\eta^2\, \lf(\frac{2^{-j}}{\eta}\rg)^\beta\ .
\ee
Combining (\ref{estseriesw+}) and (\ref{exp-1}) and using the fact that $\lf(\frac{2^{-j}}{\eta}\rg)^2\le\lf(\frac{2^{-j}}{\eta}\rg)^\beta$ we obtain
\begin{eqnarray}\label{fin-estseriesw+}
\sum_{\ell=s_1}^{s_2}\mu^{|\ell-j|}{\int_{ A_{\ell}}}|\nabla \frak{h}^+|^2 &\le& C\, \lf(\frac{2^{-j}}{\eta}\rg)^\beta\  \int_{A(2\,\eta,\delta)}|\nabla \tilde u|^2\ dx^2\ .
\end{eqnarray}
 {\bf 2.}  Using (\ref{A-17}) from lemma~\ref{ptw-harm} on $A(2\,\eta,\delta)$ we deduce that $w^-$ satisfies 
 \begin{align} 
    |\nabla \frak{h}^-(\rho,\theta)|^2& 
    \leq C\, \frac{\delta^2}{\eta^2\rho^4}\ \int_{A(2\eta,\delta)}|\nabla \frak{h}^-|^2\ dx^2 \leq C\,\frac{\delta^2}{\eta^2\rho^4}\ \int_{A(2\eta,\delta)}|\nabla \tilde u|^2\ dx^2.
\end{align}
Therefore for every $\ell\in \{s_1,\ldots, s_2\}$ we have
$$\int_{A_{\ell}}|\nabla \frak{h}^-(\rho,\theta)|^2\ dx^2\le C\,\frac{\delta^2}{\eta^2\,2^{-2\ell}}\int_{A(2\eta,\delta)}|\nabla \tilde u|^2\ dx^2\ .$$
Thus
\begin{eqnarray}\label{estseriesw-}
\sum_{\ell=s_1}^{s_2}\mu^{|\ell-j|}{\int_{ A_{\ell}}}|\nabla \frak{h}^-|^2\ dx^2 &\le& C\,\frac{\delta^2}{\eta^2}\lf[\sum_{\ell=s_1}^{s_2}\mu^{|\ell-j|}2^{2\ell}\rg]\ \int_{A(2\eta,\delta)}|\nabla \tilde u|^2\ dx^2\nonumber\\
&=&  C\,\frac{\delta^2}{\eta^2}\ \lf[\sum_{\ell=s_1}^{j}\mu^{j-\ell}2^{2\ell}+\sum_{\ell=j+1}^{s_2}\mu^{\ell-j}2^{2\ell}\rg]\int_{A(2\eta,\delta)}|\nabla \tilde u|^2\ dx^2\nonumber\\
&\leq& C\ \frac{\delta^2}{\eta^2}\ \lf[2^{2j} +\mu^{s_2-j}2^{2s_2}\rg]\ \int_{A(2\eta,\delta)}|\nabla \tilde u|^2\ dx^2 
\end{eqnarray}
We have
\be
\label{exp-2}
\mu^{s_2-j}2^{2s_2}\le \frac{\eta^2}{\delta^2}\ 2^{(s_2-j)\,\log_2\mu}=\frac{\eta^2}{\delta^2}\ \lf( \frac{\delta}{\eta\,2^{-j}} \rg)^\beta\ .
\ee
Using the fact that $\lf(\frac{\delta}{\eta\, 2^{-j}}\rg)^2\le \lf(\frac{\delta}{\eta\, 2^{-j}}\rg)^\beta$ we obtain
\begin{eqnarray}\label{fin-estseriesw-}
\sum_{\ell=s_1}^{s_2}\mu^{|\ell-j|}{\int_{ A_{\ell}}}|\nabla \frak{h}^-|^2\ dx^2 &\le& C\ \lf(\frac{\delta}{\eta\, 2^{-j}}\rg)^\beta\ \int_{A(2\eta,\delta)}|\nabla \tilde u|^2\ dx^2\ .
\end{eqnarray}
{\bf 3.}  Finally for  $w^0=C^0\log(|x|)$  we use lemma~\ref{lm-estim-c0k}:
\begin{eqnarray}\label{estseriesw0}
\sum_{\ell=s_1}^{s_2}\mu^{|\ell-j|}{\int_{ A_{\ell}}}|\nabla \frak{h}^0|^2\ dx^2 &\lesssim& (C^0)^2\,\sum_{\ell=s_1}^{s_2}\mu^{|\ell-j|}\le  (C^0)^2 =o\left(\frac{1}{\log^2\left(\frac{\eta^2}{\delta}\right)}\right)\ .
\end{eqnarray}
Combining (\ref{fin-estseriesw+}), (\ref{fin-estseriesw-}) and (\ref{estseriesw0}) we obtain (\ref{estseriesw2})  and this concludes the proof of the lemma \ref{estseriesw}. \hfill$\Box$

\medskip

\noindent{\bf Proof of Lemma~\ref{lm-ptwz}.}
We again omit to write explicitly the subscript $k$. Combining  proposition~\ref{morreyu}  and lemma~\ref{estseriesw} we obtain for any $j\in\{s_1,\ldots, s_2\}$
\be
\label{exp-3}
 \begin{array}{l}
 \ds\int_{A_j}|\nabla  \tilde u|^2\ dx^2
 \le C\left[\left(\frac{2^{-j}}{\eta}\right)^\beta+\left(\frac{\delta}{ 2^{-j}\,\eta}\right)^\beta\right]\int_{A(\eta,\delta)}|\nabla \tilde u|^2\ dx^2\\[5mm]
   \ds+ o_{\eta,k}\left(\frac{1}{\log^2\left(\frac{\eta^2}{\delta_k}\right)}\right)\ +C\,\sum_{\ell=s_1}^{s_2}\mu^{|\ell-j|}\gamma^{\ell}\left(\int_{ B_\eta(0)}|\nabla \tilde u|^2\right)^2\\[5mm]
\ds +\,C\ \int_{B_\eta(0)}|\nabla \tilde u|^2\ dx^2\,\left(\mu^{j-s_1}\int_{B_{2\eta}(0)\setminus B_\eta(0)}|\nabla \tilde u|^2\ dx^2+\mu^{s_2-j}\int_{B_{\frac{\delta}{\eta}}(0)\setminus B_{\frac{\delta}{2\,\eta}}(0)}|\nabla \tilde{u}|^2\ dx^2\right).
\end{array}
 \ee
 The $\epsilon$-regularity in \cite{Riv1} implies that for  every $\rho>0$ such that $\frac{\delta}{\eta}<\rho/4<4\rho<\eta$ we have
\begin{equation}\label{espreg}
\rho^2\,\|\nabla u\|^2_{L^{\infty}(B(0,2\rho)\setminus \bar B(0,\rho/2))}\le C\,\|\nabla u\|^2_{L^{2}(B(0,4\rho)\setminus \bar B(0,\rho/4))}.
\end{equation}
Combining (\ref{exp-3}) and (\ref{espreg}) gives for $|x|\simeq 2^{-j}$
\be
\label{exp-4}
 \begin{array}{l}
 \ds |x|^2\ |\nabla u|^2(x)
 \le C\left[\left(\frac{2^{-j}}{\eta}\right)^\beta+\left(\frac{\delta}{ 2^{-j}\,\eta}\right)^\beta\right]\int_{A(\eta,\delta)}|\nabla \tilde u|^2\ dx^2\\[5mm]
   \ds+ o_{\eta,k}\left(\frac{1}{\log^2\left(\frac{\eta^2}{\delta_k}\right)}\right)\ +C\,\sum_{\ell=s_1}^{s_2}\mu^{|\ell-j|}\gamma^{\ell}\left(\int_{ B_\eta(0)}|\nabla \tilde u|^2\ dx^2\right)^2\\[5mm]
\ds +\,C\ \int_{B_\eta(0)}|\nabla \tilde u|^2\ dx^2\,\left(\mu^{j-s_1}\int_{B_{2\eta}(0)\setminus B_\eta(0)}|\nabla \tilde u|^2\ dx^2+\mu^{s_2-j}\int_{B_{\frac{\delta}{\eta}}(0)\setminus B_{\frac{\delta}{2\,\eta}}(0)}|\nabla \tilde{u}|^2\ dx^2\right).
\end{array}
 \ee
 We estimate
 \be
 \label{exp-5}
 \sum_{\ell=s_1}^{s_2}\mu^{|\ell-j|}\,\gamma^{\ell}=\sum_{\ell=s_1}^{j}\mu^{j-\ell}\,\gamma^{\ell}+\sum_{\ell=j}^{s_2}\mu^{\ell-j}\,\gamma^{\ell}\le \mu^j\ \lf[\frac{\gamma}{\mu}\rg]^{s_1}+\gamma^j.
 \ee
 Now we use the fact that
\begin{eqnarray}
\label{exp-6}
 \ds\int_{B(0,2\eta)\setminus B(0,\eta)}|\nabla \tilde u|^2+\int_{B(0,\delta/\eta)\setminus B(0,\delta/2\eta)}|\nabla \tilde u|^2& \le &\int_{A(\eta,\delta)}|\nabla \tilde u|^2\label{exp-6}\\[5mm]
\ds \int_{ B(0,\eta)}|\nabla \tilde u|^2&\le& \int_{A(\eta,\delta)}|\nabla \tilde u|^2.\label{exp-6bis}
 \end{eqnarray}
 Combining (\ref{exp-4}), (\ref{exp-5}) and (\ref{exp-6}) we obtain (\ref{I.12}) and this concludes the proof of Lemma~\ref{lm-ptwz}.\hfill $\Box$

\section{Proof of the main result theorem~\ref{th-morse-stability}.}
\reset
Once the pointwise estimate of the gradient is established the rest of the proof of the main result, theorem~\ref{th-morse-stability}, is formally exactly the same under the general hypothesis of the theorem as for the special case $\al=0$ and ${\mathcal{N}}^n=S^n$ . We will hence be considering this later special case exclusively starting from subsection IV.2 on until the end of section IV in order to have simplified notations. 
\subsection{Introduction of the weight functions and the proof of the positive contribution of the necks to $D^2{\frak L}$}
The following elementary lemma will be central in the proof of the positive contribution of the necks to $D^2{\frak L}$ .
\begin{Lm}
\label{lm-lower-bounds}
There exists $c_0>0$ independent of $k$ and $\eta$, such that, for all $f\in W^{1,2}_0(A(\eta,\delta_k))$ we have respectively
\be
\label{I.14}
\int_{A(\eta,\delta_k)}|\nabla f|^2-\frac{c_0}{\ds \log^2\lf(\frac{\eta^2}{\delta_k}\rg)}\frac{1}{|x|^2}|f|^2\ge 0
\ee
moreover
\be
\label{I.16}
\int_{A(\eta,\delta_k)}|\nabla f|^2-\frac{\delta_k^\beta}{\eta^\beta\,|x|^\beta }\,\frac{c_0}{|x|^2}|f|^2\ dx^2\ge 0
\ee
and finally
\be
\label{I.16-a}
\int_{A(\eta,\delta_k)}|\nabla f|^2-\frac{|x|^\beta}{\eta^\beta}\,\frac{c_0}{|x|^2}|f|^2\ dx^2\ge 0.
\ee
\hfill $\Box$
\end{Lm}
\noindent{\bf Proof of lemma~\ref{lm-lower-bounds}}

{\bf Proof of (\ref{I.14})}
We consider the minimization problem
\be
\label{I.14-a}
\inf\lf\{ \int_{A(\eta,\delta_k)}|\nabla f|^2\ dx^2\ ;\ f\in W^{1,2}_0(A(\eta,\delta_k))\quad\mbox{ and }\quad\int_{A(\eta,\delta_k)}\frac{|f|^2}{|x|^2}\ dx^2=1\rg\}
\ee
Since $|x|^{-2}$ is a $C^\infty$ function on $A(\eta,\delta_k)$, the infimum is achieved. Let $f$ be such a minimum. Then for some $\lambda\in\R$ it solves the following Euler Lagrange equation
\be
\label{I.14-aa}
\lf\{
\begin{array}{l}
\ds-\Delta f=\frac{\la}{|x|^2}\ f\quad\mbox{ in }A(\eta,\delta_k)\\[5mm]
\ds f=0\quad\mbox{ on }\p A(\eta,\delta_k)
\end{array}
\rg.
\ee
 and introduce
\[
f^\ast(r):=\sqrt{\frac{1}{2\pi}\int_0^{2\pi} |f(r,\theta)|^2\ d\theta}\ .
\]
We have clearly
\[
\int_{A_{\eta,\delta}}\frac{|f|^2}{|x|^2}\ dx^2=\int_{A_{\eta,\delta}}\frac{|f^\ast|^2}{|x|^2}\ dx^2
\]
Moreover
\be
\label{I.14-b}
\begin{array}{l}
\ds|\nabla f^\ast|^2=\lf|\frac{d f^\ast}{dr}\rg|^2=\lf|\ds\frac{\frac{1}{2\pi}\ds\int_0^{2\pi}2\,\frac{\p |f|}{\p r}\,|f|(r,\theta)\ d\theta}{\ds 2\,\sqrt{\frac{1}{2\pi}\int_0^{2\pi} |f(r,\theta)|^2\ d\theta}}\rg|^2\le \lf|\sqrt{\frac{1}{2\pi} \int_0^{2\pi}|\p_r|f||^2\ d\theta}\rg|^2\\[10mm]
\ds\quad=\frac{1}{2\pi} \int_0^{2\pi}|\p_r|f||^2\ d\theta\ \le\frac{1}{2\pi} \int_0^{2\pi}|\nabla f|^2(r,\theta)\ d\theta
\end{array}
\ee
This implies
\be
\label{I.14-c}
\int_{A(\eta,\delta_k)}|\nabla f^\ast|^2\ dx^2\le \int_{A(\eta,\delta_k)}|\nabla f|^2\ dx^2
\ee
Hence (\ref{I.14-a}) admits an axially symmetric minimizer which satisfies the ODE
\be
\label{I.14-d}
-{f}''-\frac{{f}'}{r}-\frac{\la}{r^2} f=0\quad\mbox{ and } \quad f(\delta_k/\eta)=f(\eta)=0\ .
\ee
Let $Y(s):=f(e^s)$. One has $$\dot{Y}(s)= e^s\, f'(e^s)=r \, f'(r)\mbox{ and }\ddot{Y}(s)= e^s\, f'(e^s)+ e^{2s}\,f''(e^s)=r \, f'(r)+r^2\,f''(r)$$
Multiplying by $r^2$ (\ref{I.14-d}) and proceeding to the change of variable $r=e^s$, the ODE becomes
\be
\label{I.14-e}
-\ddot{Y}(s)-\la\,Y(s)=0\quad Y(\log (\delta_k/\eta))=0\quad\mbox{and}\quad Y(\log(\eta))=0\ .
\ee
Hence $Y$ has the form
\be
\label{I.14-f}
Y(s)=A\, \cos(\sqrt{\la} s)+ B\,\sin(\sqrt{\la} s)
\ee
The smallest $\la>0$ for which there exists a solution which vanishes both at $\log (\delta/\eta)$ and $\log(\eta)$ is the one for which $$\sqrt{\la}\, [\log\eta-\log(\delta/\eta)]=\pi $$
Therefore the minimizer we are considering satisfies
\be
\label{I-14-g}
\lf\{
\begin{array}{l}
\ds-\Delta f=\frac{\pi^2}{|x|^2\log^2\frac{\eta^2}{\delta_k^2}}\ f\quad\mbox{ in }A(\eta,\delta_k)\\[5mm]
\ds f=0\quad\mbox{ on }\p A(\eta,\delta_k)
\end{array}
\rg.
\ee
Multiplying by $f$ and integrating by parts gives then
\be
\label{I-14-h}
\int_{A(\eta,\delta_k)}|\nabla f|^2\ dx^2=\frac{\pi^2}{\log^2\frac{\eta^2}{\delta_k^2}}\int_{A(\eta,\delta_k)}\frac{| f|^2}{|x|^2}\ dx^2\ .
\ee
Since $f$ is a minimizer we obtain (\ref{I.14}) for $c_0=\pi^2$ in fact.



\medskip

{\bf Proof of (\ref{I.16})}  Let 
\[
y:=\frac{\delta}{\eta}\frac{x}{|x|^2}
\]
and $\ti{f}(y)=f(x)$. We have in particular
\[
|y|=\frac{\delta}{\eta\,|x|}\quad\mbox{ and }\quad c_1\,\lf(\frac{\delta}{\eta}\rg)^2\ \frac{dx^2}{|x|^4} =dy^2
\]
Hence
\be
\label{I.17}
\int_{A_{\eta,\delta_k}}|\nabla f|^2-\frac{\delta_k^\beta}{\eta^\beta\,|x|^\beta }\,\frac{1}{|x|^2}|f|^2\ dx^2=\int_{A_{1,\delta_k/\eta^2}}|\nabla \ti{f}|^2-|y|^\beta\frac{c}{|y|^2}|\ti{f}|^2\ dy^2
\ee
We have that $W^{1,2}_0(B_1(0))$ embeds continuously into $L^{2q}$ for any $1< q<+\infty$ (Poincar\'e Sobolev inequality). For such a $q$, extending $\ti{f}$ by $0$ inside $B_{\delta/\eta^2}(0)$, we have in one hand
\be
\label{I.18}
\int_{B_1(0)}|\nabla\ti{f}|^2(y)\ dy^2\ge C_q \lf[\int_{B_1(0)}|\ti{f}|^{2q}(y)\ dy^2\rg]^{1/q}
\ee
and, choosing $q$ large enough in such a way that $(2-\beta)q' < 2$ where $(q')^{-1}=1-q^{-1}$
\be
\label{I.19}
\int_{B_1(0)}|y|^\beta\frac{c}{|y|^2}|\ti{f}|^2\ dy^2\le \lf[ \int_{0}^1\frac{s}{s^{2q'-\beta q'}} \ ds \rg] \lf[\int_{B_1(0)}|\ti{f}|^{2q}(y)\ dy^2\rg]^{1/q}\le \ti{C}_q\ \int_{B_1(0)}|\nabla\ti{f}|^2(y)\ dy^2
\ee
Hence combining (\ref{I.17}) and (\ref{I.19}) we obtain (\ref{I.16}) for $c_0:=c/\ti{C}_q$.

\medskip

{\bf Proof of (\ref{I.16-a})} We introduce $\ti{f}(y):=f(y\,\eta)$ hence we have
\be
\label{I.20}
\int_{A(1,\delta_k/\eta^2)}|\nabla \ti{f}|^2-|y|^\beta \frac{c_0}{|y|^2}|\ti{f}|^2(y)\ dy^2=\int_{A(\eta,\delta_k)}|\nabla f|^2-\frac{|x|^\beta}{\eta^\beta}\,\frac{c_0}{|x|^2}|f|^2\ dx^2
\ee
and (\ref{I.16-a}) follows from (\ref{I.19}).

This concludes the proof of lemma~\ref{lm-lower-bounds}.\hfill $\Box$

\bigskip

We introduce the following {\bf weight function}\footnote{We could have simply extended $\om_{\eta,k}$ in $B_{\delta_k}{\eta}(0)$ by
\[
\frac{\eta^2}{\delta_k^2}\,\lf[1+\frac{\delta_k^\beta}{\eta^{2\,\beta}}+\frac{1}{ \log^2 \frac{\eta^2}{\delta_k}}\rg]\ .
\]
The reason why we extend $1$ by $ \frac{1}{\eta^4}\frac{(1+\eta^2)^2}{(1+\delta_k^{-2}\,|x|^2)^2}$ starting at $|x|=\delta_k/\eta$ inside $B_{\delta_k}{\eta}(0)$ is because, after the composition with the stereographic projection$\pi$  it will become simply a constant in $\pi^{-1}(B_{\delta_k}{\eta}(0))$ which offers some advantages in later computations. Both choices would lead in any case to the same conclusion.
} 
\be
\label{I.21}
\om_{\eta,k}(x)
\lf\{
\begin{array}{l}
\ds=\frac{1}{|x|^2}\,\lf[\frac{|x|^\beta}{\eta^\beta}+\frac{\delta_k^\beta}{\eta^\beta\,|x|^\beta}+\frac{1}{ \log^2 \frac{\eta^2}{\delta_k}}\rg]\quad\mbox{ in }A(\eta,\delta_k)\\[8mm]
\ds =   \frac{1}{\eta^2}\,\lf[1+\frac{\delta_k^\beta}{\eta^{2\,\beta}}+\frac{1}{ \log^2 \frac{\eta^2}{\delta_k}}\rg]\quad\mbox{ in } \Sigma\setminus B_\eta(0)\\[8mm]
\ds =\frac{\eta^2}{\delta_k^2}\lf[  \frac{1}{\eta^4}\frac{(1+\eta^2)^2}{(1+\delta_k^{-2}\,|x|^2)^2}+\frac{\delta_k^\beta}{\eta^{2\beta}}+\frac{1}{ \log^2 \frac{\eta^2}{\delta_k}}\rg]\quad\mbox{ in }B_{\delta_k/\eta}(0)
\end{array}
\rg.
\ee
Combining Lemma~\ref{lm-ptwz} and lemma~\ref{lm-lower-bounds} we obtain the following lemma
\begin{Lm}
\label{lm-lower-bound-neck}
There exists $\la_0>0$ independent of $k$ and there exists $\eta_0$ such that for any $0<\eta<\eta_0$ and for $k$ large enough the following holds
\be
\label{I.22}
\forall\ {w}\,\in V_{u_k}\quad\mbox{ s. t. }\quad{w}=0\quad \mbox{ in }\quad\Sigma\setminus A(\eta,\delta_k)\quad\quad Q_{u_k}(w)\ge\, \la_0\, \int_{\Sigma}\om_{\eta,k}\ w^2\ dvol_h\ .
\ee
\hfill $\Box$
\end{Lm}
\subsection{ The diagonalization of $D^2{\frak L}$ with respect to the weights $\om_{\eta,k}$.}
We denote by
\[
<w,v>_{\om_{\eta,k}}=\int_\Sigma \ w\,v\ \om_{\eta,k}\, dvol_h\ .
\]
 We consider the diagonalization of the self-adjoint operator ${\mathcal L}_{\eta,k}$ with respect to $<\ ,\ >_{\om_{\eta,k}}$ given by
\[
Q_{u_k}(w)=<{\mathcal L}_{\eta,k}w,w>_{\om_{\eta,k}}=\int_\Sigma \lf[\om_{\eta,k}^{-1}\,\Delta_h w-\om_{\eta,k}^{-1}\,|du_k|_h^2\, w\rg] w\ \om_{\eta,k}\ dvol_h.
\]
We consider the diagonalization of ${\mathcal L}_{\eta,k}$ and we denote by ${\mathcal E}_{k,\eta}(\la)$ the eigenspace for the eigenvalue $\la$. We have the following lemma
\begin{Lm}
\label{lm-ind-uk}
\be
\label{I.23}
\mbox{Ind}_E(u_k)=\mbox{dim}\lf[\bigoplus_{\la<0}{\mathcal E}_{\eta,k}(\la)\rg]\ .
\ee
\hfill $\Box$
\end{Lm}
\noindent{\bf Proof of lemma~\ref{lm-ind-uk}.} On the one hand, the restriction of $Q_{u_k}$ to $W:=\oplus_{\la<0}{\mathcal E}_{\eta,k}(\la)$ is strictly negative hence $\mbox{Ind}_E(u_k)\ge\mbox{dim}\lf[\oplus_{\la<0}{\mathcal E}_{\eta,k}(\la)\rg]$.
Vice versa, let $W$ be a sub-vectorspace of $L^2$ on which $Q$ is strictly negative, we obviously have
\be
\label{I.24}
W\cap \bigoplus_{\la\ge 0}{\mathcal E}_{\eta,k}(\la)=\{0\}\ .
\ee
This implies that dim$W\le \mbox{dim}\lf[\oplus_{\la<0}{\mathcal E}_{\eta,k}(\la)\rg]$. Hence we have proved lemma~\ref{lm-ind-uk}.\hfill $\Box$ 

\medskip

We shall be using the following Lemma
\begin{Lm}
\label{lm-lower-spectr}
There exist two constants $\eta_0>0$ and  $\mu_0>0$ in dependent of $k$ and a family of constants $\mu_{\eta,k}>0$  satisfying 
\be
\label{I.25-b}
\lim_{\eta\rightarrow 0}\limsup_{k\rightarrow +\infty} \mu_{\eta,k}=0\ \quad \mbox{ and }\quad \forall k\in{\N},~~\forall \eta \in(0,\eta_0)\quad\mbox{ one has }\quad 0<\mu_{\eta,k}<\mu_0\ ,
\ee
and such that for any element $\lambda\in{\R}$ with
\be
\label{I.25-a}
\mbox{dim}\lf( {\mathcal E}_{\eta,k}(\la) \rg)>0\quad\Longrightarrow\quad\la\ge -\,\mu_{\eta,k}\ge - \mu_0\ .
\ee
\hfill $\Box$
\end{Lm}
\noindent{\bf Proof of lemma~\ref{lm-lower-spectr}.} 

Observe that, thanks to  lemma~\ref{lm-ptwz} and using respectively (\ref{I.8}) and (\ref{I.9}),  there exists a positive constant $C_0$ independent of $\eta$ small enough ($\eta<\eta_0$) and $k$ such that
\be
\label{I.25}
\forall \ x\in A(\eta,\delta_k) \quad |du_k|^2_h(x)\le C_0\ \om_{\eta,k}(x)\ ,
\ee
Moreover 
\be
\label{I.250-a}
\lim_{\eta\rightarrow 0}\limsup_{k\rightarrow +\infty}\lf\|\frac{|du_k|^2_h}{\om_{\eta,k}}\rg\|_{L^\infty(A(\eta,\delta_k))}=0\ .
\ee
{\bf Claim:}
\be
\label{I.250-b}
\sup_{0<\eta<\eta_0}\sup_{k\in \N}\lf\|\frac{|du_k|^2_h}{\om_{\eta,k}}\rg\|_{L^\infty(\Sigma)}<+\infty.\ee
{\bf Proof of Claim.} We prove (\ref{I.250-b}) by contradiction. Suppose that  there is  sequence $\eta_n\to 0$ as $n\to +\infty$ such that 
\be
\label{I.250-c}
\sup_{k\in \N}\lf\|\frac{|du_k|^2_h}{\om_{\eta_n,k}}\rg\|_{L^\infty(\Sigma)}=+\infty.
\ee
Because of (\ref{I.25}),   if (\ref{I.250-b}) would not be true we would contradict
the strong $C^1$ convergence respectively of $u_k$ towards $u_\infty$ on $\Sigma\setminus B_{\eta_n}(x_k)$ and  of $v_k(y):=u_k(\delta_k\,y+x_k)$ towards $v_\infty$. Hence (\ref{I.250-b}) holds true.
We conclude the proof of the Claim.
\par
\medskip
\underline{We show \eqref{I.25-b}}.\par

We introduce respectively on $\Sigma$
\be
\label{om-eta-infty}
\om_{\eta,\infty}(x)=\lf\{
\begin{array}{l}
\ds\frac{1}{\eta^2}\quad\mbox{ in }\Sigma\setminus B_\eta(p)\\[5mm]
\ds\frac{|x|^\beta}{\eta^\beta}\,\frac{1}{|x|^2}\quad\mbox{ in } B_\eta(p)\ ,
\end{array}
\rg.
\ee
and on $\C$
\be
\label{om-infty}
\hat{\om}_{\eta,\infty}(x)=\lf\{
\begin{array}{l}
\ds\frac{1}{\eta^2}\,\frac{(1+\eta^2)^2}{(1+|y|^2)^2}\quad\mbox{ for }|y|\le \eta^{-1}\\[5mm]
\ds \frac{1}{\eta^\beta}\frac{1}{|y|^{2+\beta}}\quad\mbox{ for }|y|\ge \eta^{-1}\ .
\end{array}
\rg.
\ee
We have on the one hand
\be
\label{I.25-c}
\limsup_{k\rightarrow +\infty}\lf\|\frac{|du_k|^2_h}{\om_{\eta,k}}\rg\|_{L^\infty(\Sigma\setminus B_\eta(x_k))}=\lf\|\frac{|du_\infty|^2_h}{\om_{\eta,\infty}}\rg\|_{L^\infty(\Sigma\setminus B_\eta(x_k))}=O(\eta^2)\ ,
\ee
and on the other hand
\be
\label{I.250-c}
\begin{array}{l}
\ds\limsup_{k\rightarrow +\infty}\lf\|\frac{|du_k|^2_h}{\om_{\eta,k}}\rg\|_{L^\infty(B_{\delta_k\,\eta^{-1}}(x_k))}=\limsup _{k\rightarrow +\infty}\lf\|\frac{\delta_k^2\,|du_k|^2_h}{\delta_k^2\,\om_{\eta,k}}\rg\|_{L^\infty( B_{\delta_k\,\eta^{-1}}(x_k))}\\[8mm]
\ds=\limsup_{k\rightarrow +\infty}  \lf\|\frac{|\nabla v_k|^2}{\delta_k^2\,\om_{\eta,k}(\delta_k\,y)}\rg\|_{L^\infty(B_{\eta^{-1}}(0))} =\lf\|\frac{|\nabla v_\infty|^2}{\hat{\om}_{\eta,\infty}}\rg\|_{L^\infty( B_{\eta^{-1}}(0))}=O(\eta^2)\ .
\end{array}
\ee
where we use the fact that, thanks to the point removability on $S^2$ we have $|\nabla v_\infty|^2\le C\, (1+|y|^2)^{-2}$.\par
Then if we set 
\[
\mu_{\eta,k}:=\lf\|\frac{|du_k|^2_h}{\om_{\eta,k}}\rg\|_{L^\infty(\Sigma)}\ .
\]
the sequence  $\mu_{\eta,k}$ satisfies \eqref{I.25-b}.\par
\medskip
\underline{We next show that \eqref{I.25-a} holds as well}.

Let $w$ be a non zero eigenfunction for a negative eigenvalue $\la$ for ${\mathcal L}_{\eta,k}$. One has
\be
\label{I.26}
\int_\Sigma |dw|^2_h-|du_k|^2_h\,w^2\ dvolg_h=\la\, \int_\Sigma\om_{\eta,k}\,w^2\ dvol_h
\ee
 
We have from (\ref{I.26})
\be
\label{I.27}
\la\, \int_\Sigma\om_{\eta,k}\,w^2\ dvol_h\ge-\int_\Sigma\,|du_k|^2_h\,w^2\ dvolg_h\ge - \,\mu_{\eta,k}\,\int_\Sigma\om_{\eta,k}\,w^2\ dvol_h
\ee
This implies (\ref{I.25-a}) and lemma~\ref{lm-lower-spectr} is proved.\hfill $\Box$

\medskip

We define
\[
<w,v>_{\om_{\eta,\infty}}=\int_\Sigma \ w\,v\ \om_{\eta,\infty}\, dvol_h\ .
\]
Let $w\in V_{u_\infty}$. We consider the diagonalization of the self-adjoint operator ${\mathcal L}_{\eta,\infty}$ with respect to $<\ ,\ >_{\om_{\eta,\infty}}$ given by
\[
Q_{u_\infty}(w)=<{\mathcal L}_{\eta,\infty}w,w>_{\hat{\om}_{\eta,\infty}}=\int_\Sigma \lf[\om_{\eta,\infty}^{-1}\,\Delta_h w-\om_{\eta,\infty}^{-1}\,|du_\infty|_h^2\, w\rg] w\ \om_{\eta,\infty}\ dvol_h.
\]
Let $w\in V_{v_\infty}$. We consider the  self-adjoint operator $\hat{\mathcal L}_{\infty}$ with respect to $<\ ,\ >_{\hat{\om}_{\eta,\infty}}$ given by
\be
\label{opera-infty}
\begin{array}{l}
\ds Q_{v_\infty}(w)=<\hat{\mathcal L}_{\eta,\infty}w,w>_{\hat{\om}_{\eta,\infty}}=\int_{{\C}} \lf[\hat{\om}_{\eta,\infty}^{-1}\,\Delta w-\hat{\om}_{\eta,\infty}^{-1}\,|dv_\infty|^2_{\C}\, w\rg] w\ \hat{\om}_{\eta,\infty}\ dy^2\\[5mm]
\ds=\ti{Q}_{\ti{v}_\infty}(\ti{w}) =<\tilde{\mathcal L}_{\infty}\ti{w},\ti{w}>_{\om_{\infty}}=\int_{S^2} \lf[\ti{\om}_\infty^{-1}\,\Delta_{S^2} \ti{w}-\ti{\om}_\infty^{-1}\,|d\ti{v}_\infty|_{S^2}^2\, \ti{w}\rg] \ti{w}\ \ti{\om}_{\infty}\ dvol_{S^2}\ ,
\end{array}
\ee
where $\ti{v}_\infty:=v_\infty\circ\pi$ and $\ti{w}:=w\circ\pi$ where $\pi$ is the stereographic projection from $S^2$ into $\C$ sending the south pole to infinity and $$\ti{\om}_{\eta,\infty}:= [\hat{\om}_{\eta,\infty}(y)\ (1+|y|^2)^2]\circ \pi$$
Denote by $z=(z_1,z_2,z_3)$ the coordinates on $S^2$, because of (\ref{om-infty}), since
\[
\pi^{-1}(y_1,y_2)=\lf(\frac{y_1}{1+|y|^2},\frac{y_2}{1+|y|^2},\frac{1-|y|^2}{1+|y|^2}\rg)
\]
we have in a neighborhood of the south pole denoting $r^2=z_1^2+z_2^2\simeq\mbox{dist}^2(z,\mbox{South})$
\be
\label{pot-infty}
\mbox{dist}^2(z,\mbox{South})\simeq \frac{|y|^2}{(1+|y|^2)^2}\simeq\frac{1}{|y|^2}\ .
\ee
This gives thanks to (\ref{om-infty})
\be
\label{om-infty-asym}
0< C_\eta^{-1}\le \ti{\om}_{\eta,\infty}(z)\le C_\eta\ \mbox{dist}^{-2+\beta}(z,\mbox{South}).
\ee
Hence we can conclude with the following result:
\begin{Lm}
\label{app-lemma-appendix}
 The conclusion of lemma~\ref{lm-appendix-1}  holds respectively for $(\Sigma, \om_{\eta,\infty},{\mathcal L}_{\eta,\infty})$ and for $(S^2,\ti{\om}_{\eta,\infty}, \ti{\mathcal L}_{\eta,\infty})$.\hfill $\Box$
 \end{Lm}

 We next introduce 
\[
{\mathcal E}^0_{\eta,\infty}:=\lf\{ w\in V_{u_\infty}\ ;\ w\in\bigoplus_{\la\le 0}{\mathcal E}_{\eta,\infty}(\la)\rg\}\ ,
\]
where
\[
{\mathcal E}_{\eta,\infty}(\la):=\lf\{w\in V_{u_\infty}\  ;\ \Delta_hw-|du_\infty|^2_h\,w=\la\ \om_{\eta,\infty}\, w\ \rg\}
\]
 
Because of lemma~\ref{lm-morse-ind}, ${\mathcal E}^0_{\eta,\infty}$ is a finite (hence closed) sub-vector space of $L^2_{\om_{\eta,\infty}}(\Sigma)$ and it is a sub-space  on which $Q_{u_\infty}\le 0$. 
Hence
\be
\label{I.27-a}
\mbox{dim}\lf({\mathcal E}^0_{\eta,\infty}\rg)\le \mbox{Ind}_E(u_\infty)+\mbox{Null}_E(u_\infty)\ ,
\ee
where $\mbox{Null}_E(u_\infty)$ corresponds to the kernel of the bilinear form associated to $Q_{u_\infty}$.

We denote 
\[
V_{v_\infty}=\lf\{ w\in W^{1,2}({\C},S^{n-1})\ ;\ w\cdot v_\infty\equiv 0\quad\mbox{ on }{\C}  \rg\}
\]
and 
\[
\hat{\mathcal E}^0_{\eta,\infty}:=\lf\{ w\in V_{v_\infty}\ ;\ w\in\bigoplus_{\la\le 0}\hat{\mathcal E}_{\eta,\infty}(\la)\rg\}\ ,
\]
where
\[
\hat{\mathcal E}_{\eta,\infty}(\la):=\lf\{\phi\in V_{v_\infty}\  ;\ \Delta_h\phi-|dv_\infty|^2_h\,\phi=\la\ \hat{\om}_{\eta,\infty}\,w\ \rg\}.
\]
We have the following result:

\begin{Prop}\label{eig}
\be
\label{I.33bis}
\mbox{dim}\lf(\hat{\mathcal E}^0_{\eta,\infty}\rg)\le \mbox{Ind}_E(\ti{v}_\infty)+\mbox{Null}_E(\ti{v}_\infty)\ .
\ee
where $\mbox{Null}_E(\ti{v}_\infty)$ corresponds to the kernel of the bilinear form associated to $Q_{\ti{v}_\infty}$ with $\ti{v}_\infty:=v_\infty\circ\pi$.\hfill$\Box$
\end{Prop}

{\bf Proof of proposition \ref{eig}.}
Recall that the composition $\ti{v}_\infty:=v_\infty\circ\pi$, where $\pi$ is the stereographic projection sending the south pole to infinity, extends as a smooth harmonic map over the whole $S^2$. Moreover, for any $w\in V_{v_\infty}$, thanks to conformal invariance the map 
$\ti{w}:=w\circ \pi$ has finite $W^{1,2}$ energy in $S^2\setminus B_\ep(\{South\})$ and
\[
\limsup_{\ep\rightarrow 0}\int_{S^2\setminus B_\ep(\mbox{South})}|d\ti{w}|^2_{S^2}\, dvol_{S^2}=\int_{{\C}}|\nabla w|^2(y)\ dy^2<+\infty
\]
On $S^2$ we define 
\[
\ti{w}_\ep(x):= \lf[1-\chi\lf(\frac{|x|}{\ep}\rg)\rg] \ti{w}- \chi\lf(\frac{|x|}{\ep}\rg)\frac{1}{|B_{2\ep}(0)\setminus B_\ep(0)|}\ \int_{B_{2\ep}(0)\setminus B_\ep(0)} \ti{w}\ dvol_{S^2}
\]
where $|x|$ is the geodesic distance to the south pole (we are taking normal coordinates around the south pole). We have
\be
\label{I.28}
\begin{array}{l}
\ds\int_{S^2}|d\ti{w}_\ep|^2_{S^2}\ dvol_{S^2}\\[5mm]
\ds\le 2\ \int_{S^2\setminus B_\ep(0)}|d\ti{w}|^2_{S^2}\, dvol_{S^2}+2\ \|\chi'\|_\infty\ \frac{1}{\ep^2}\ \int_{B_{2\ep}(0)\setminus B_\ep(0)} |\ti{w}-\ov{w}_\ep|^2\ dvol_{S^2}
\end{array}
\ee
where 
\[
\ov{w}_\ep:=\frac{1}{|B_{2\ep}(0)\setminus B_\ep(0)|}\ \int_{B_{2\ep}(0)\setminus B_\ep(0)} \ti{w}\ dvol_{S^2}
\]
Using Poincar\'e inequality we finally obtain
\be
\label{I.29}
\begin{array}{l}
\ds\int_{S^2}|d\ti{w}_\ep|^2_{S^2}\ dvol_{S^2}\le C\ \int_{{\C}}|\nabla w|^2(y)\ dy^2
\end{array}
\ee
Hence $\ti{w}_\ep$ is uniformly bounded in $W^{1,2}(S^2,{\R}^n)$ and it converges almost everywhere to $\ti{w}$ over $S^2$ hence it weakly converges $\ti{w}$ and this implies the $\ti{w}\in W^{1,2}(S^2,{\R}^n)$.
We have moreover that $\ti{w}\cdot\ti{v}_\infty=0$ almost everywhere on $S^2$ hence $\ti{w}\in V_{\ti{v}_\infty}$. We have proved that the composition with the stereographic projection is realizing
a continuous isomorphism  $\Pi$ between $(V_{v_\infty}, \dot{W}^{1,2}({\C},{\R}^n))$ into $(V_{\ti{v}_\infty}, \dot{W}^{1,2}(S^2,{\R}^n))$. Observe that
\be
\label{I.30}
Q_{v_\infty}(w)=Q_{\ti{v}_\infty}(\ti{w})\ ,
\ee
hence
\be
\label{I.30-aa}
W:=\Pi\lf[\hat{\mathcal E}^0_{\eta,\infty}\rg]
\ee
realizes a sub-vectorspace of $V_{\ti{v}_\infty}$ on which $Q_{\ti{v}_\infty}\le 0$ and one has
\be
\label{I.33}
\mbox{dim}\lf(\hat{\mathcal E}^0_{\eta,\infty}\rg)\le \mbox{Ind}_E(\ti{v}_\infty)+\mbox{Null}_E(\ti{v}_\infty)\ .
\ee
where $\mbox{Null}_E(\ti{v}_\infty)$ corresponds to the kernel of the bilinear form associated to $Q_{\ti{v}_\infty}$.\hfill$\Box$

\subsection{ The proof of the main result theorem~\ref{th-morse-stability}.} 

We prove now the following lemma which implies theorem~\ref{th-morse-stability} in the case ${\mathcal{N}}^n=S^n$ and $\al=0$.
\begin{Lm}
\label{lm-argument}
Under the previous notations for $\eta$ sufficiently small and $k$ large enough  we have
\be
\label{I.35}
\mbox{Ind}_E(u_k)+\mbox{Null}_E(u_k)=\mbox{ dim}\lf(\bigoplus_{\la\le0}{\mathcal E}_{\eta,k}(\la)\rg)\le \mbox{Ind}_E({u}_\infty)+\mbox{Null}_E(u_\infty)+\mbox{Ind}_E(\ti{v}_\infty)+\mbox{Null}_E(\ti{v}_\infty)\ .
\ee
\hfill $\Box$
\end{Lm}

\noindent{\bf Proof of lemma~\ref{lm-argument}.} We consider the finite dimensional sphere given by 
\[
{\mathcal S}_{\eta,k}:=\lf\{ w\in \bigoplus_{\la\le 0}{\mathcal E}_{\eta,k}(\la)\ ;\ <w,w>_{\om_{\eta,k}}=1\rg\}
\]
Let us take $w_k\in {\mathcal S}_{\eta,k}$. We have thanks to lemma~\ref{lm-lower-spectr}
\be
\label{I.36}
\lf\{
\begin{array}{l}
\ds\int_\Sigma\om_{\eta,k}\ w_k^2\ dvol_h=1\ \quad\mbox{ and }\\[5mm]
\ds 0\ge\int_\Sigma| d w_k|^2_h-|du_k|^2_h\, w_k^2\ dvol_h\ge -\mu_0\, \int_\Sigma\om_{\eta,k} \ w_k^2\ dvol_\Sigma=-\mu_0\ .
\end{array}
\rg.
\ee
Because of (\ref{I.250-b})   we have
\be
\label{I.37}
\int_\Sigma|du_k|^2_h\, w_k^2\ dvol_h\le \mu_0\ ,
\ee
hence combining (\ref{I.36}) and (\ref{I.37}) gives
\be
\label{I.37-a}
\int_\Sigma| d w_k|^2_h\  dvol_h\le \mu_0\ .
\ee
 From \eqref{I.37-a} it follows   (up to sequence) that
\be
\label{I.38}
w_k\rightharpoonup w_\infty\quad\mbox{ in }W^{1,2}(\Sigma,h)\quad\mbox{ and }\quad w_k(\delta_k y+x_k)\rightharpoonup \sigma_\infty\quad\mbox{ in }W^{1,2}_{loc}({\C})\ .
\ee

\medskip

\noindent{\bf Claim 1 :}  $w_\infty\ne 0$ or $\sigma_\infty\ne 0$. 

\medskip

\noindent{\bf Proof of the claim 1 :} Let $(\phi_k^j)_{j=1\ldots N_k}$ be an orthonormal basis\footnote{It is very important to develop the argument without assuming that $N_k$ is uniformly bounded because it is not known a-priori and it will come as a direct consequence of the main result.} with respect to $L^2_{\om_{\eta,k}}$ of $\bigoplus_{\la\le 0}{\mathcal E}_{\eta,k}(\la)$ made of eigenfunctions of ${\mathcal L}_{\eta,k}$ :
\be
\label{I.39}
P_{u_k}\lf[\Delta_h\phi_k^j\rg]-|du_k|^2_h\ \phi_k^j=\la_k^j\ \om_{\eta,k}\ \phi_k^j\quad\quad \forall \, j=1\ldots N_k\ .
\ee
where $P_{u_k}$ is the pointwise orthogonal projection onto $T_{u_k} S^{n}$. Since $\phi_k^j$ is section of $u_k^{-1}T S^{n-1}$ we have $P_{u_k}(\phi_k^j)=\phi_k^j$ and since $\Delta_h$ is a positive laplacian
\be
\label{I.39-a}
\begin{array}{l}
\ds P_{u_k}\lf[\Delta_h\phi_k^j\rg]=\Delta_h\phi_k^j+2\,dP_{u_k}\cdot_h d\phi_k^j-\Delta_hP_{u_k}\,\phi_k^j
\end{array}
\ee
Since $w_k\in {\mathcal S}_{\eta,k}$ we have the existence of $c_k^{\,j}$ such that
\be
\label{I.42}
w_k=\sum_{j=1}^{N_k}c_k^{\,j}\ \phi_k^{\,j}\quad\mbox{ where }\quad \sum_{j=1}^{N_k} (c_k^{\,j})^2=1\ .
\ee
We have
\be
\label{I.42-a}
\Delta_hw_k+2\,dP_{u_k}\cdot_h d w_k-\Delta_hP_{u_k}\,w_k-|du_k|^2_h\ w_k=\sum_{j=1}^{N_k}\la_k^{\,j}\ \om_{\eta,k}\ c_k^{\,j}\ \phi_k^{\,j}
\ee
Since $\la_k^{\,j}\in [-C_0/\eta^2,0]$, this gives
\be
\label{I.42-b}
\begin{array}{l}
\ds\int_\Sigma\frac{1}{\om_{\eta,k}}\lf|\Delta_hw_k+2\,dP_{u_k}\cdot_h d w_k-\Delta_hP_{u_k}\,w_k-|du_k|^2_h\ w_k\rg|^2\ dvol_h\\[5mm]
\ds\quad=\int_{\Sigma}\sum_{i,j=1}^{N_k}\la_k^j\ \la_k^i\ \om_{\eta,k}\,c_k^i\ \phi_k^i\ c_k^j\ \phi_k^j\ dvol_h=\sum_{j=1}^n(\la_k^j)^2\ (c_k^j)^2\le \mu_0^2
\end{array}
\ee
Hence there exists $f_k$ uniformly bounded in $L^2$ such that
\be
\label{I.42-c}
\Delta_hw_k+2\,dP_{u_k}\cdot_h d w_k-\Delta_hP_{u_k}\,w_k-|du_k|^2_h\ w_k= \sqrt{\om_{\eta,k}}\ f_k\ 
\ee
Because of (\ref{I.8}), using classical elliptic estimates we obtain
\be
\label{I.43}
\|\nabla^2 w_k\|_{L^2(\Sigma\setminus B_\eta(0))}\le C_\eta  
\ee
Let $v_k(y):=u_k(\delta_k\, y+x_k)$ and $\sigma_k(y):=w_k(\delta_k\, y+x_k)$. We have
\be
\label{I.43-aa}
\Delta_{\C}\sigma_k+2\,dP_{v_k}\cdot d \sigma_k-\Delta_{\C}P_{v_k}\,\sigma_k-|dv_k|^2_h\ \sigma_k= \delta_k^2\ \sqrt{\om_{\eta,k}}\ f_k(\delta_k\, y+x_k)\ .
\ee
We have
\be
\label{I.43-ab}
\int_{B_{1/\eta}(0)}\delta_k^4\ {\om_{\eta,k}}\ f^2_k(\delta_k\, y+x_k)\ dy^2\le C_\eta \ \int_{B_{\delta_k/\eta}(x_k)}|f_k|^2\ dvol_h\ .
 \ee
Hence, thanks to (\ref{I.9}), we have
\be
\label{I.43-a}
\|\nabla^2_{y}\sigma_k\|_{L^2(B_{1/\eta}(0))}=\|\nabla^2_{y}(w_k(\delta_k\,y+x_k))\|_{L^2(B_{1/\eta}(0))}\le C_\eta\ .
\ee
Hence
\be
\label{I.44}
w_k\rightharpoonup w_\infty\quad\mbox{ in }W^{2,2}_{loc}(\Sigma\setminus \{0\})\quad\mbox{ and }\quad \sigma_k(y)=w_k(\delta_k y+x_k)\rightharpoonup \sigma_\infty\quad\mbox{ in }W^{2,2}_{loc}({\C})\ .
\ee
Assume $w_\infty=0$ and $\sigma_\infty=0$. Let
\be
\label{I.45}
\check{w}_k:= w_k\ \chi\lf(2\,\frac{|x-x_k|}{\eta}   \rg)\, \lf(1-\chi\lf( \eta\, \frac{|x-x_k|}{\delta_k}   \rg)\rg)\in W^{1,2}_0(A_{\eta,k},{\R}^n)\ .
\ee
Because of (\ref{I.44}) and assuming $w_\infty=0$ and $\sigma_\infty=0$ we have
\be
\label{I.46}
\lim_{k\rightarrow +\infty}\|\nabla(w_k-\check{w}_k)\|_{L^2(\Sigma)}=0\ .
\ee
Observe that on the one hand
\be
\label{I.47}
\begin{array}{l}
\ds\lf|Q_{u_k}(w_k)-Q_{u_k}(\check{w}_k)\rg|\le\lf|\int_{\Sigma\setminus B_{\eta/2}(x_k)}|\nabla w_k|^2-|\nabla \check{w}_k|^2+|du_k|^2_h\,[w^2_k-\check{w}^2_k]\ dvol_h\rg|\\[5mm]
\ds\quad +\lf|\int_{B_{2\delta_k/\eta}(x_k)}|\nabla w_k|^2-|\nabla \check{w}_k|^2+|du_k|^2_h\,[w^2_k-\check{w}^2_k]\ dvol_h\rg|\ ,
\end{array}
\ee
and on  the other hand
\be
\label{I.48}
\begin{array}{l}
\ds\lf|\int_\Sigma \om_{k,\eta}\,w_k^2\ \ dvol_h-\int_\Sigma \om_{k,\eta}\,\check{w}_k^2\ \ dvol_h\rg|\le \lf|\int_{\Sigma\setminus B_{\eta/2}(x_k)}  \om_{k,\eta}\,[w_k^2-\check{w}_k^2]\ \ dvol_h\rg|\\[5mm]
\ds \quad+ \lf|\int_{B_{2\delta_k/\eta}(x_k)}  \om_{k,\eta}\,[w_k^2-\check{w}_k^2]\ \ dvol_h\rg|
\end{array}
\ee
Since we have supposed that $w_{\infty}=\sigma_{\infty}=0$ then  the last two integrals in \eqref{I.48} go to zero as $k\to +\infty$. Therefore 
combining (\ref{I.46}), (\ref{I.47}) and (\ref{I.48}) gives
\be
\label{I.49}
\lim_{k\rightarrow +\infty}\lf|Q_{u_k}(w_k)-Q_{u_k}(\check{w}_k)\rg|=0\quad\mbox{ and }\quad \lim_{k\rightarrow +\infty}\int_\Sigma \om_{k,\eta}\,\check{w}_k^2\ \ dvol_h=1\ . \ee
Since $\check{w}_k\in W^{1,2}_0(A_{\eta,k},{\R}^m)\cap V_{u_k}$, we have thanks to lemma~\ref{lm-lower-bound-neck}, for any $0<c<c_0$,
\be
\label{I.50}
Q_{u_k}(\check{w}_k)\ge c\ \int_\Sigma \om_{k,\eta}\,\check{w}_k^2\ \ dvol_h\quad.
\ee
Hence, combining (\ref{I.49}) and (\ref{I.50}) gives
\be
\label{I.51}
\liminf_{k\rightarrow+\infty} Q_{u_k}({w}_k)\ge c>0\ ,
\ee
which contradicts the fact that $Q_{u_k}({w}_k)\le 0$. Hence the claim 1 is proved and we have either $w_\infty\ne 0$ or $\sigma_\infty\ne 0$.

\medskip

\noindent{\bf End of the proof of lemma~\ref{lm-argument}.} Let $N>0$ and $(\phi_k^j)_{j=1\ldots N}$ be a free orthonormal family for $L^2_{\om_{\eta,k}}(\Sigma)$ of eigenfunctions of ${\mathcal L}_{\eta,k}$ where we choose a subsequence that we still index by $k$ . Because of the claim 1 we have, modulo extraction of a subsequence, for each $j=1\ldots N$
\be
\label{I.51-a}
\phi_k^j\rightharpoonup \phi_\infty^j \quad\mbox{ weakly in }W^{2,2}_{loc}(\Sigma\setminus\{p\})\quad \mbox{ and }\quad \phi_k^j(\delta_k^j\, y+x_k^j)\rightharpoonup \sigma_\infty^j \quad\mbox{ weakly in }W^{2,2}_{loc}({\C})
\ee
Passing to the limit as $k\rightarrow +\infty$ we obtain that the maps $\phi_\infty^j $ satisfy for some $\la_\infty^j\le 0$
\be
\label{I.51-b}
P_{u_\infty}\lf[\Delta_h\phi^j_\infty\rg]-|du_\infty|^2_{h}\,\phi^j_\infty=\la_\infty^j\ \om_{\eta,\infty}\ \phi^j_\infty\ ,
\ee
and $\sigma_\infty^j$ satisfies
\be
\label{I.51-c}
 P_{v_\infty}\lf[  \Delta_{\C} \sigma^j_\infty \rg]-|\nabla v_\infty|^2\ \sigma^j_\infty =\la_\infty^j \ \hat{\om}_{\eta,\infty}\,\sigma^j_\infty\ .
\ee
Because of claim 1 for all $j=1\ldots N$ we have $(\phi^j_\infty,\sigma^j_\infty)\ne (0,0)$. Assume 
\be
\label{I.51-d}
N>\mbox{dim}\,{\mathcal E}^0_{\eta,\infty}+\mbox{dim}\,\hat{\mathcal E}^0_{\eta,\infty}
\ee
Then thee family $(\phi^j_\infty,\sigma^j_\infty)_{j=1\ldots N}$ is linearly dependent and there exists $(c^1_\infty,\ldots, c^N_\infty)\ne (0,\ldots, 0)$ such that
\be
\label{I.51-e}
\sum_{j=1}^N c^j_\infty\,\phi^j_\infty=0\quad\mbox{ and }\quad \sum_{j=1}^N c^j_\infty\,\sigma^j_\infty=0\ .
\ee
Considering $w_k:=\sum_{j=1}^Nc_\infty^j\, \phi_k^j$ we obtain a contradiction. Hence we have proved $N\le\mbox{dim}\,{\mathcal E}^0_{\eta,\infty}+\mbox{dim}\,\hat{\mathcal E}^0_{\eta,\infty}$ which implies, using (\ref{I.27-a}) and (\ref{I.33}),  lemma~\ref{lm-argument}.\hfill $\Box$
\section{Proof of theorem~\ref{th-morse-stability-deg}}
 \reset
We first prove the following result
\begin{Lm}
\label{lm-collar}
Let $u_k$ be a sequence of critical points of ${\frak L}$ from $A(\eta,\delta_k)=B_\eta(0)\setminus B_{\delta_k/\eta}(0)$ into ${\mathcal{N}}^n$ such that
\be
\label{VI.1}
\lim_{\eta\rightarrow 0}\lim_{k\rightarrow +\infty}\ \sup_{ \delta_k/\eta<\rho<2\rho<\eta}\ \int_{B_{2\rho}(x_k)\setminus B_\rho(x_k)}|du_k|^2_h\ dvol_h=0\ .
\ee
Suppose 
\be
\label{VI.2}
0\le \La:=\lim_{\eta\rightarrow 0}\,\limsup_{k\rightarrow +\infty}\int_{\eta^{-1}\,\delta_k}^\eta\dashint_0^{2\pi}\lf|\frac{d{u}_k}{dr}\rg|\ d\theta\,dr \ .
 \ee
Then there exists $\eta_0>0$ there exits $0<\beta<2$ and $C>0$ independent of $k$ such that for any $\eta<\eta_0$  and any $x\in A(\eta,\delta_k)$
\be
\label{VI.3}
\begin{array}{l}
\ds |x|^2\ |\nabla u_k|^2(x)\le C\ \lf[\frac{|x|^\beta}{\eta^\beta}+\lf(\frac{\delta_k}{\eta\,|x|}\rg)^\beta\rg]\  \int_{A(2\eta,\delta_k)}|du_k|^2_h\ dvol_h+C\, \frac{\La^2}{\log^{2}\lf( \frac{\eta^2}{\delta_k}  \rg)}.
\end{array}
\ee
 \hfill $\Box$
\end{Lm}
\noindent{\bf Proof of lemma~\ref{lm-collar}.} 
\par
{\bf Claim 1: Quantization of $L^2$ norm of $u_k$ in $A(\eta,\delta_k)$: } The assumptions \eqref{VI.1} and \eqref{VI.2} imply that
\begin{equation}\label{qud}
\lim_{\eta\rightarrow 0}\limsup_{k\rightarrow +\infty}\int_{A(\eta,\delta_k)}|\nabla u_k|^2dx =0.
\end{equation}
{\bf Proof of the Claim 1.}  \par
\noindent {\bf Case 1.} We first assume  that 
\begin{equation}
 \|\nabla u_k\|_{L^2(A(\eta,\delta_k))}\le \varepsilon_0
\end{equation}
where $\varepsilon_0$ is the constant appearing in Theorem I.4 in \cite{Riv1}. 
As in subsection III.2, in the collar $A(\eta,\delta_k)$ we decompose
\[
A_k\,\nabla \tilde u_k=\nabla\varphi_k+\nabla^\perp\psi_k+\nabla {\frak h}_k
\]
where $\tilde u_k$ is the Whitney extension of $u_k$ as  in Lemma \ref{lma-extension} and
 \begin{eqnarray}
 \Delta\varphi_k&=&\nabla^\perp B_k\cdot\nabla \tilde u_k~~\mbox{in $\C$}\label{0.11bis}\\
\Delta\psi_k&=&\nabla^\perp A_k\cdot\nabla \tilde u_k~~\mbox{in $\C$}\label{0.11tris}\\
  \Delta\frak{h}_k&=&0~~\mbox{in $A(\eta,\delta_k) $}\label{0.11qris}
 \end{eqnarray}
 The matrix $B_k$ satisfies the condition in \eqref{0.9}.
 The generalized Wente inequality  (see \cite{Beth}) implies that  
 \begin{eqnarray} 
 \|\nabla\varphi_k\|_{L^2(\C)}&\le &C  \|\nabla B_k\|_{L^2(\C)} \|\nabla \tilde{u}_k\|_{L^{2,\infty}(\C)}\label{linftyphi}\\
 \|\nabla \psi_k\|_{L^2(\C)}&\le &C  \|\nabla A_k\|_{L^2(\C)} \|\nabla\tilde{u}_k\|_{L^{2,\infty}(\C)}\label{linftypsi}
 \end{eqnarray}
 Since
 \begin{equation}
  \|\nabla\tilde{u}_k\|_{L^{2,\infty}(\C)}\leq C (\|\nabla{u}_k\|_{L^{2,\infty}(A(\eta/2,\delta_k)}+ \|\nabla\tilde{u}_k\|_{L^{2}(B_{\eta}\setminus B_{\eta/2})}+ \|\nabla\tilde{u}_k\|_{L^{2}(B_{2\delta_k/\eta}\setminus B_{\delta_k/\eta})}),
\ee
then the $\varepsilon$-regularity on the annular domain $A(\eta,\delta_k)$ (see Lemma 13  in \cite{LaRi2}) and the assumption \eqref{VI.1} imply  that
\begin{equation}
\label{VI.4}
\lim_{\eta\rightarrow 0}\limsup_{k\rightarrow+\infty}\|\nabla\varphi_k\|_{L^{2}(A(\eta,\delta_k))}+\|\nabla\psi_k\|_{L^{2}(A(\eta,\delta_k))} =0\ee
Lemma 2.2 in \cite{MiRi} (or Lemma 10 in \cite{LaRi2}) imply that for all $t\in (0,1)$ there exists $C_t$ such that
\begin{equation}
\label{VI.4bis}
\lim_{\eta\rightarrow 0}\limsup_{k\rightarrow+\infty}\|\nabla\frak{h}^{\pm}_k\|_{L^2(t\eta,\delta_k)}\le C_t \lim_{\eta\rightarrow 0}\limsup_{k\rightarrow+\infty}\|\nabla\frak{h}^{\pm}_k\|_{L^{2,\infty}(\eta,\delta_k)}=0.
\end{equation}
  The following estimate holds:
\begin{eqnarray}\label{logu}
\lim_{\eta\rightarrow 0}\limsup_{k\rightarrow+\infty}\left|\int_{A(\eta,\delta_k)}\nabla u_k \cdot \nabla\log(|x|) dx\right|
&\le& \lim_{\eta\rightarrow 0}\limsup_{k\rightarrow+\infty}\int_{\delta_k/\eta}^\eta\int_0^{2\pi}\left|\frac{\partial u}{\partial r}\right|=:\Lambda<+\infty
\end{eqnarray}
Because of   \eqref{VI.4} and \eqref{VI.4bis} we also have that
\begin{eqnarray}
&&\lim_{\eta\rightarrow 0}\limsup_{k\rightarrow+\infty}\left|\int_{A(\eta,\delta_k)}\nabla \varphi_k \cdot \nabla\log(|x|) dx\right|<+\infty\label{logp}\\
&&\lim_{\eta\rightarrow 0}\limsup_{k\rightarrow+\infty}\left|\int_{A(\eta,\delta_k)}\nabla \psi_k \cdot \nabla\log(|x|) dx\right|<+\infty\label{logps}\\
&&\lim_{\eta\rightarrow 0}\limsup_{k\rightarrow+\infty}\left|\int_{A(\eta,\delta_k)}\nabla \frak{h}^{\pm}_k \cdot \nabla\log(|x|) dx\right|<+\infty.\label{loghpm}
\end{eqnarray}
By combining \eqref{logu}-\eqref{loghpm} we finally get
\begin{equation}
\lim_{\eta\rightarrow 0}\limsup_{k\rightarrow+\infty}\int_{A(\eta,\delta_k)}\nabla \frak{h}^0_k \cdot \nabla\log(|x|) dx=
\lim_{\eta\rightarrow 0}\limsup_{k\rightarrow+\infty} C^0_{\eta,k}\log(\frac{\eta^2}{\delta_k})<+\infty.
\end{equation}
In particular  it follows that
\begin{eqnarray}\label{hok}
&&\lim_{\eta\rightarrow 0}\limsup_{k\rightarrow+\infty}\|\nabla \frak{h}^0_k\|_{L^2(A(\eta,\delta_k)}
=\lim_{\eta\rightarrow 0}\limsup_{k\rightarrow+\infty}C^0_{\eta,k}\log^{1/2}(\frac{\eta^2}{\delta_k})=0
\end{eqnarray}
From \eqref{V.4} and \eqref{V.5} it follows that $\|A_k^{-1} \|_{L^{\infty}(A(\eta,\delta_k)}\le C$ for some $C$ independent on $k$ and $\eta.$
Hence by combining \eqref{VI.4},  \eqref{VI.4bis} and \eqref{hok} we get \eqref{qud} and we can conclude the proof of Claim 1.
 \par
 \noindent {\bf Case 2.}  In the general case where
\begin{equation}
 \|du_k\|_{L^2(A(1,\delta_k))}>\varepsilon_0
\end{equation}
   one proceeds as in the proof of Theorem 2 in \cite{LaRi1}.
The rest of the proof of lemma~\ref{lm-collar} follows word by word from the second part of section III.\hfill $\Box$ 
\par
\medskip

For $C\, \La^2<c_0$ where $c_0$ is the constant in (\ref{I.14}) we obtain lemma~\ref{lm-lower-bound-neck}. The rest of {\bf the proof of theorem~\ref{th-morse-stability-deg}}  is deduced from  lemma~\ref{lm-lower-bound-neck} similarly as in the non degenerating case.\hfill $\Box$

\appendix
\addcontentsline{toc}{section}{Appendices}
\section*{Appendices}

\section{The finiteness of the Morse index}
\reset
\begin{Lma}
\label{lm-morse-ind} Let $u$ be a critical point of ${\frak L}$ from a closed Riemann surface $(\Sigma,h)$ into a closed sub-manifold ${\mathcal{N}}^n$ of an euclidian space ${\R}^m$.
Denote by 
\[
Ind_{\frak L}(u):=\sup\lf\{\mbox{dim}(W) \ ;\ W\mbox{ is a sub vector-space of }V\mbox{ s.t. }\lf.Q_u\rg|_{W}<0\rg\}
\]
where $Q_u$ is the second derivative of $\frak L$ at $u$ and is given by
\[
Q_u(w):=D^2{\frak L}_u(w)=\int_\Sigma|dw|^2- <w,S_u(du)_h w>+<w,{\mathcal H}_uw>\, dvol_h\ .
\]
where $S_u(du)_h$ and ${\mathcal H}$ are given respectively by (\ref{II.9-b}) and (\ref{II.18}). This defines a quadratic form defined on $\Gamma(u^{-1}T{\mathcal{N}}^n)$, the sections of the pull-back bundle by $u$ of the tangent bundle to ${\mathcal{N}}^n$. We have 
\be
\label{I.4}
Ind_{\frak L}(u)<+\infty\ .
\ee
Similarly,  denote $\mbox{Null}_{\frak L}(u)$  the dimension of the Kernel of the bilinear form associated to $Q_u$, there holds
\be
\label{I.4-a}
\mbox{Null}_{\frak L}(u)<+\infty\ .
\ee
\hfill $\Box$
\end{Lma}
\noindent{\bf Proof of lemma~\ref{lm-morse-ind}}
Denote by ${\mathcal E}(\la)$ the eigenspaces of the Laplace Beltrami operator defined on $W^{1,2}(\Sigma,{\R})$ :
\[
{\mathcal E}(\la):=\lf\{\phi\in W^{1,2}(\Sigma,{\R})\ ; \ \Delta_h\phi=\la\phi\rg\}\ .
\]
For any $\La\in {\R}$ we shall denote
\[
{\mathcal E}^\La:=\bigoplus_{\la\le\La}{\mathcal E}(\la)
\]
It is well known that, for any $\La\in {\R}$, dim$({\mathcal E}^\La)<+\infty$.\footnote{ We observe that  $\Delta^{-1}_h$ is a compact operator from $L^2(\Sigma)$ to $L^2(\Sigma)$, being $\Sigma$ a compact manifold. Hence for every $\lambda\in\R$, dim$({\mathcal E}(\la))<+\infty$. }  We choose $\La> C\,\||du|_h\|_{L^\infty(\Sigma)}^2$ where $C$ is chosen such that
\[
\lf|\int_\Sigma- <w,S_u(du)_h w>+<w,{\mathcal H}_uw>\, dvol_h\rg|\le \frac{1}{2}\int_\Sigma|dw|^2\ +C \ |du|^2_h\ |w|^2\, dvol_h
\]
 and we introduce
\be
\label{I.5}
V^\La_u:=\lf\{w\in V_u\ ; \ \int_{\Sigma} w_i\, \phi\ dvol_h=0\ ,\ \forall\phi\in {\mathcal E}^\La\ ,\ i=1\ldots n\rg\}\ .
\ee
We have in one hand that the space
\[
\lf\{V_u\ni w\rightarrow \int_{\Sigma} w_i\, \phi\ dvol_h\ ; \ \phi\in {\mathcal E}^\La\ ,\ i=1\ldots n\rg\} \]
is a finite dimensional subspace of $(V_u)^\ast$ and hence $V^\La_u$ is finite co-dimensional. Moreover
\[
\begin{array}{l}
\ds\forall w\in V^\La_u\setminus \{0\}\quad Q_u(w)\ge\frac{1}{2}\sum_{i=1}^n\int_{\Sigma}\lf[{\La}-C\,|du|^2_h\rg]\,|w_i|^2\ dvol_h\ge \frac{1}{2}\, \lf[{\La}-C\,\||du|^2_h\|_\infty\rg]\ \|w\|_{L^2}^2>0
\end{array}
\]
For that reason, for any $W\mbox{ sub vector-space of }V_u\mbox{ s.t. }\lf.Q_u\rg|_{W}<0$ one has and we deduce
\be
\label{I.6}
W\cap V^\La_u=\{0\}\ .
\ee
Since $V^\La_u$ is finite co-dimensional we deduce the projection parallel to $V^\La_u$  of $W$ onto any finite dimensional supplement of $V^\La_u$ in $V_u$ must be injective and we deduce that
\be
\label{I.7}
\mbox{dim}(W)\le\mbox{codim} (V^\La_u)\ .
\ee
A similar argument can be applied to prove (\ref{I.4-a}). This ends the proof of  lemma~\ref{lm-morse-ind}.\hfill $\Box$

\medskip
 
\section{The discreteness of the spectrum of $D^2{\frak L}$ }
\reset
\reset
\begin{Lma}
\label{lm-appendix-1}
Let $u$ be a smooth map from a closed oriented Riemannian 2-dimensional manifold $(\Sigma,h)$ into a closed $C^2$ sub-manifold ${\mathcal{N}}^n$ of ${\R}^m$. Let $\al$ is an arbitrary $C^2$ 2-form of ${\mathcal{N}}^n$. On $\Sigma$ we consider a  positive function ${\om}$, smooth away from one point $p$ and such that there exists $0<\beta<1$ such that
\be
\label{A-1}
0<C^{-1}_0\,\le \om(x)\le \frac{C_0}{\mbox{dist}(x,p)^{2-\beta}}\  \quad\mbox{ and }\quad |\nabla^l\om^{-1}|(x)|\le C_l\, \mbox{dist}(x,p)^{2-l-\beta} \quad\forall\, l=0,1,2\ .
\ee
Denote
\[
L^2_\om(\Sigma):=\lf\{\phi\in L^2(\Sigma)\ ;\ \int_\Sigma\om\, \phi^2\ dvol_h<+\infty\rg\} \ .
\]
Denote for any $z\in {\mathcal{N}}^n$ by $P_z$ the $m\times m$ matrix corresponding to the orthogonal projection onto $T_z{\mathcal{N}}^n$. The composition $P\circ u$ will simply be denoted $P_u$. Let 
\[
W_u:=\lf\{w\in L^2_\om(\Sigma,{R}^m)\ ; \ P_u(w(x))=w(x)\quad\mbox{ for a. e. } x\in \Sigma\rg\}\ .
\]
On $W_u$ we consider the operator 
\[
{\mathcal L}w:= \om^{-1}\,P_u\lf[\Delta_h w+{\mathcal H}_uw-S_u(du)_h w\rg]\ .
\]
where $S_u(du)_h$ and ${\mathcal H}_u$ are given respectively by (\ref{II.9-b}) and (\ref{II.18}). Then there exists an Hilbert basis of $W_u$ for the $L^2_\om$ scalar product made of the eigenfunctions of $w$ and the eigenvalues
statisfy
\[
\la_1<\la_2<\la_3.....\rightarrow +\infty\,.
\]
\hfill $\Box$
\end{Lma}
\noindent{\bf Proof of lemma~\ref{lm-appendix-1}.} Observe first that ${\mathcal L} w$ has a distributional sense for any $w\in W_u$. Indeed
\be
\label{A-2}
{\mathcal L}w=\om^{-1}\,\Delta_hw-2\, \om^{-1}\,dP_{u}\cdot_h dw-\, \om^{-1}\,\Delta_h P_u\, w +\om^{-1}\,P_u\lf[{\mathcal H}_uw\rg]-\om^{-1}\,S_u(du)_h w
\ee
We recall that $P_u\lf[H_u(dw\wedge du)\rg]=H_u(dw\wedge du).$ Therefore we can write    
\begin{eqnarray}
\label{A-2-a}
  P_u\lf[{\mathcal H}w\rg]&=&P_u\lf[2\,H_u(dw\wedge du)+\nabla_w H(du,du)\rg]\\[5mm]
  &=&2\,H_u(dw\wedge du)+P_u\lf[\nabla_w H(du,du)\rg]\\[5mm]
  &=&2d(H_u(w\wedge du))-2\partial_{\xi_{\ell}}H_{i,j}^k(u)du^{\ell}\wedge w^{i}\wedge du^{j}+P_u\lf[\nabla_w H(du,du)\rg]\\[5mm]
 &=&2\,d\lf[H_u(w\wedge du)\rg]+K_u(du)_h w\ ,
\end{eqnarray}
 where 
 $$K_u(du)_h w=-2\partial_{\xi_{\ell}}H_{i,j}^k(u)du^{\ell}\wedge w^{i}\wedge du^{j}+P_u\lf[\nabla_w H(du,du)\rg]$$ 
 with
 \be
 \label{A-2-b}
 |K_u(du)_h|\le C\, |du|^2_h\ .
\ee
Because of the hypothesis~\eqref{A-1} we have that 
\be
\label{A-3}
\lf\{
\begin{array}{l}
\ds|\Delta_h\om^{-1}|(x)\le C_2\, \mbox{dist}(x,p)^{-\beta} \in L^2(\Sigma)\ ,\\[3mm]
d^\ast [\om^{-1}\,dP_{u}]\in L^2(\Sigma)\\[5mm]
\om^{-1}\,\Delta_h P_u\,+\om^{-1}\,[K_u(du)-S_u(du)_h] \in L^2(\Sigma).
\end{array}
\rg.
\ee
 Hence ${\mathcal L}w\in {\mathcal D}'(\Sigma)$.
 
 \medskip
 
 Denote ${\mathcal L}_\la={\mathcal L}+\la\,Id$. 
 
 \medskip
 
 \noindent{\bf Claim 1 :}  {\it There exists $\la_0>0$ such that for any $\la\ge \la_0$ the following holds
 \be
 \label{A-4}
 \lf\{
 \begin{array}{l}
 \ds\forall f\in W_u\quad \exists\,!\ w\in W_u\cap W^{1,2}(\Sigma,{\R}^n)\quad\mbox{ s. t. }\quad{\mathcal L}_\la w=f\\[5mm]
 \ds  \int_\Sigma[|dw|^2_h+\om\,|w|^2]\ dvol_h\le 2 \|f\|^2_{L^2_\om(\Sigma)}
 \end{array}
 \rg.
 \ee}
\noindent{\bf Proof of Claim 1} We start first by establishing a-priori estimates. Multiplying the equation by $w\in W_u\cap W^{1,2}(\Sigma,{\R}^m)$ and integrating by parts gives
\be
\label{A-5}
\int_\Sigma|dw|^2_h\ dvol_h+[\la\, \om\, w^2-<w,S_u(du)_h w>+w{\mathcal H}w]\ dvol_h=\int_\Sigma \om\, f\,w\ dvol_h\le \|w\|_{L^2_\om(\Sigma)}\  \|f\|_{L^2_\om(\Sigma)}\ 
\ee
Because of (\ref{A-2-a}) and (\ref{A-2-b}), for any $\ep>0$ there exists $C_\ep>0$ such that
\be
\label{A-5-b}
\lf|\int_\Sigma w{\mathcal H}w\ dvol_h\rg|\le \ep\int_\Sigma |dw|^2_h\ dvol_h+C_\ep\,\int_\Sigma |du|^2_h\,|w|^2\ dvol_h\ .
\ee
Hence for $\la\ge C\,\||du|^2_h\|_\infty$ for $C>0$ large enough independent of $w$, we deduce
 \be
 \label{A-6}
 \int_\Sigma[|dw|^2_h+\om\,|w|^2]\ dvol_h\le 2\,\|w\|_{L^2_\om(\Sigma)}\  \|f\|_{L^2_\om(\Sigma)}\
 \ee
 which implies
 \be
 \label{A-7}
  \int_\Sigma[|dw|^2_h+\om\,|w|^2]\ dvol_h\le 4\, \|f\|^2_{L^2_\om(\Sigma)}
 \ee
 The previous a-priori estimate is giving the uniqueness of $w$ in $W_u\cap W^{1,2}(\Sigma,{\R}^m)$. 
 
 \medskip
 
 Regarding the proof of the existence, we consider on  $W_u\cap W^{1,2}(\Sigma,{\R}^m)$ the minimization of
 \be
 \label{A-8}
 \begin{array}{l}
\ds L_\la(w):=\int_\Sigma[|dw|^2_h+\la\, \om\, w^2-<w,S_u(du)_h w>+w{\mathcal H}w-2\,\om\, f\cdot w]\ dvol_h\\[5mm]
\ds\quad\ge\frac{1}{2} \int_\Sigma[|dw|^2_h+ 2\,\om\,w^2]\  dvol_h-2\,\|f\|_{L^2_\om(\Sigma)} \, \|w\|_{L^2_\om(\Sigma)}\\[5mm]
\ds\quad\ge \frac{1}{2} \|dw\|_{L^2(\Sigma)}^2+\lf[\|w\|_{L^2_\om(\Sigma)}-\|f\|_{L^2_\om(\Sigma)}\rg]^2- \|f\|^2_{L^2_\om(\Sigma)}
 \end{array}
 \ee
 for $\la\ge C\,\||du|^2_h\|_\infty$. For a minimizing sequence $w_k$ the energy $ \|w_k\|_{W^{1,2}(\Sigma)}^2$ remains uniformly bounded 
  we can extract a subsequence such that $w_k\rightharpoonup w_\infty$ in $W^{1,2}$ and, thanks to Rellich Kondrachov theorem, the weak $W^{1,2}$ convergence  implies that the pointwise condition
  $P_{u}(w_k)=0$ passes to the limit and we obtain that $w_\infty\in W_u\cap W^{1,2}(\Sigma,{\R}^m)$.

  Moreover, the Euler-Lagrange variational theory implies that  the following holds
 \be
 \label{A-9}
 \forall\,\varphi\in W_u\cap W^{1,2}(\Sigma,{\R}^n)\quad\int_\Sigma d\varphi\cdot_h dw_\infty+[\la\, \om\, w_\infty- S(du)_h w_\infty +{\mathcal H}w_\infty]\cdot\varphi-\om\,f\cdot \varphi\ dvol_h=0\ .
 \ee
 and we deduce that ${\mathcal L}_\la w_\infty=f$. This concludes the proof of the Claim 1.
 
 \medskip
 
 \noindent{\bf Claim 2.}{\it
 \be
 \label{A-10}
{ \mathcal K}:=\lf\{ w\in W_u\cap W^{1,2}(\Sigma,{\R}^m)\ ;\ \int_\Sigma[|dw|^2_h+\om\,|w|^2]\ dvol_h\le 1\rg\} 
 \ee
 is a compact subset of $W_u$}. 
 
 \medskip
 
 \noindent{\bf Proof of Claim 2.} We consider a sequence $w_k\in {\mathcal K}$. It is weakly pre-compact  in $W^{1,2}(\Sigma,{\R}^n)$ and we can then extract a subsequence such that 
 $w_k\rightharpoonup w_\infty$ weakly in $W^{1,2}$. Clearly, by lower semi-continuity of the norm $L^2_\om(\Sigma)$ for such a convergence and since the pointwise constraint $P_u(w_k)$ passes to the limit 
 we have $w_\infty$. We have for any $1<p<2/(2-\beta)$
 \be
 \label{A-11}
 \int_\Sigma\om\,|w_k-w_\infty|^2\ dvol_h\le \|\om\|_{L^p(\Sigma)}\ \|w_k-w_\infty\|^2_{L^{2p'}(\Sigma)}\ .
 \ee
 Rellich Kondrachov's theorem gives the compactness of the continuous embedding $W^{1,2}(\Sigma)\hookrightarrow L^q(\Sigma)$ for any $q<+\infty$ since $\Sigma$ is compact. Because of (\ref{A-1}) $\|\om\|_{L^p}<+\infty$ for any $p<2/(2-\beta)$. Thus we deduce that
 \be
 \label{A-12}
 \lim_{k\rightarrow +\infty}\|w_k-w_\infty\|_{L^2_\om(\Sigma)}=0\ .
 \ee
 Any sequence of ${ \mathcal K}$ posses a subsequence which strongly converges for the metric topology induced by $L^2_\om(\Sigma)$ this implies that ${\mathcal K}$ is a compact subset of $W_u$ and {\bf Claim 2 is proved}.
 
 \medskip
 
 We have proved so far that the self adjoint operator ${\mathcal L}^{-1}_\la$ for the $L^2_\om$ norm is compact. $L^2_\om$ defines a separable Hilbert space. Theorem 6.11 of \cite{Bre} is implying the conclusion of lemma~\ref{lm-appendix-1} and this is ending the proof.\hfill $\Box$

 \section{A Whitney type extension lemma on annuli.}
 \reset
 \begin{Lma}
 \label{lma-extension}
 Let $0< 2r<R$ and $a\in W^{1,2}(B_R(0)\setminus B_r(0))$ then there exists $\ti{a}\in W^{1,2}({\C})$ with $$\mbox{Supp}(\nabla\ti{a})\subset B_{2R}(0)\setminus B_{r/2}(0)\ ,$$ such that
 \be
 \label{A-13}
 \int_{\C}|\nabla \ti{a}|^2\ dx^2\le C\ \int_{B_R(0)\setminus B_r(0)}|\nabla a|^2\ dx^2\ ,
 \ee
  where $C$ is independent of $r$ and $R$. We can also impose moreover
  \be
  \label{A-14}
   \int_{B_{2R}(0)\setminus B_R(0)}|\nabla \ti{a}|^2\ dx^2\le C\ \int_{B_R(0)\setminus B_{R/2}(0)}|\nabla a|^2\ dx^2\ ,
  \ee
  and
  \be
  \label{A-15}
   \int_{B_{r}(0)\setminus B_{r/2}(0)}|\nabla \ti{a}|^2\ dx^2\le C\ \int_{B_{2r}(0)\setminus B_{r}(0)}|\nabla a|^2\ dx^2.
  \ee
  \hfill $\Box$
  \end{Lma}
\noindent{\bf Proof of lemma~\ref{lma-extension}.} Let $\chi$ be a function in $C^\infty({\R}_+)$ such that $\chi\equiv 1$ on $[0,1]$ and $\chi\equiv 0$ on $[2,+\infty)$. We denote for any $t>0$
\[
\forall \,x\in {\C}\quad\quad\chi_t(x):=\chi\lf(\frac{|x|}{t}\rg)\ .
\]
For $|x|\ge r$ we define
\[
\hat{a}(x):=\chi_r(x)\,a(x)+(1-\chi_r(x))\,\ov{a}_r\quad\mbox{ where }\quad\ov{a}_r:=\frac{1}{|B_{2r}\setminus B_r|}\int_{B_{2r}\setminus B_r} a(x)\ dx
\]
Because of Poincar\'e inequality there holds
\be
\label{wh1}
\int_{{\C}\setminus B_r}|\nabla \hat{a}|^2\ dx^2\le C\  \int_{B_{2r}(0)\setminus B_{r}(0)}|\nabla a|^2\ dx^2\ .
\ee
Now, we extend $a$ inside $B_r$ by taking
\[
\ti{a}(x):=\hat{a}\lf(r^2\,\frac{x}{|x|^2}\rg)
\]
Since the inversion is conformally invariant we have
\be
\label{wh2}
\int_{B_r}|\nabla \ti{a}|^2\ dx^2=\int_{{\C}\setminus B_r}|\nabla \hat{a}|^2\ dx^2\ .
\ee
We do something similar in order to extend $a$ outside $B_R$. For $|x|<R$ we define
\[
\breve{a}(x):=(1-\chi_R(x))\,a(x)+\chi_R(x)\,\ov{a}_R\quad\mbox{ where }\quad\ov{a}_R:=\frac{1}{|B_{R}\setminus B_{R/2}|}\int_{B_{R}\setminus B_{R/2}} a(x)\ dx
\]
Because of Poincar\'e inequality there holds
\be
\label{wh3}
\int_{B_R}|\nabla \breve{a}|^2\ dx^2\le C\  \int_{B_{R}(0)\setminus B_{R/2}(0)}|\nabla a|^2\ dx^2\ .
\ee
Now, we extend $a$ outside $B_R$ by taking
\[
\ti{a}(x):=\hat{a}\lf(R^2\,\frac{x}{|x|^2}\rg)
\]
Since the inversion is conformally invariant we have
\be
\label{wh4}
\int_{{\C}\setminus B_R}|\nabla \ti{a}|^2\ dx^2=\int_{ B_R}|\nabla \breve{a}|^2\ dx^2\ .
\ee
The proof of lemma~\ref{lma-extension} is complete.\hfill $\Box$

\section{Pointwise gradient estimates of harmonic functions on annuli}
 \reset
\begin{Lma}
\label{ptw-harm}
Let $\eta^2> 4\,\delta$. Let $\frak{h}$ be a real harmonic functions defined on $A(\eta,\delta):=B_\eta(0)\setminus B_{\delta/\eta}(0)$ given in Fourier by 
\[
\frak{h}=\frak{h}^++\frak{h}^-\quad\mbox{ where }\ \frak{h}^+=\Re\lf[\sum_{n>0}\frak{h}_nz^n\rg]\quad\mbox{and}\quad \frak{h}^-=\Re\lf[\sum_{n<0}\frak{h}_nz^n\rg]
\]
then for any $\rho\in [2\delta/\eta,\eta/2]$ and $\theta\in [0,2\pi]$ there holds
\be
\label{A-16}
|\nabla \frak{h}^+(\rho,\theta)|^2\le \frac{C}{\eta^2}\ \int_{A(\eta,\delta)}|\nabla \frak{h}^+|^2 dx^2\ ,
\ee
and 
\be
\label{A-17}
|\nabla \frak{h}^-(\rho,\theta)|^2\le \frac{C\,\delta^2}{\eta^2}\ \frac{1}{\rho^4}\,\int_{A(\eta,\delta)}|\nabla \frak{h}^-|^2 dx^2\ .
\ee
\hfill $\Box$
\end{Lma}
\noindent{\bf Proof of lemma~\ref{ptw-harm}} Let $\frak{f}^+(z):=\sum_{n>0}\frak{h}_nx^n$ and $\frak{f}^-(z):=\sum_{n<0}\frak{h}_nx^n$. We have
\be
\label{A-18}
2\, |\nabla \frak{h}^\pm|^2=2\,|\nabla\Re(\frak{f}^\pm)|^2=2\,|\nabla\Im(\frak{f}^\pm)|^2=|\nabla \frak{f}^\pm|^2\ .
\ee
We have
\begin{align}
    \int_{A(\eta,\delta)}|\nabla \frak{f}^+|^2&=4\pi\sum_{n>0}|\frak{h}_n|^2n^2\int_{\delta/\eta}^\eta\rho^{2n-1}d\rho\nonumber\\
    &=2\pi\sum_{n>0} |\frak{h}_n|^2n\bigg((\eta)^{2n}-\lf(\frac{\delta}{\eta}\rg)^{2n}\bigg)\nonumber\\
    &\geq \pi\sum_{n>0}|\frak{h}_n|^2\,n\ \eta^{2n}
\end{align}
where we have used in the last inequality that $\delta\leq \frac{1}{2}\eta^2$. For $\rho<\eta/2$ we have
\begin{align}\label{pestwplus}
    |\nabla \frak{f}^+(\rho,\theta)|^2&\leq\left(2\sum_{n>0}|\frak{h}_n|n\rho^{n-1}\right)^2\nonumber
    \\&\leq \sum_{n>0}|\frak{h}_n|^2\,n\,\eta^{2n}\sum_{n>0}\frac{\rho^{2n-2}}{\eta^{2n}}n\nonumber\\
    &=\left( \sum_{n>0}|\frak{h}_n|^2\,n\,\eta^{2n}\right)\frac{1}{\eta^2}\sum_{n>0}n\left(\frac{\rho^2}{\eta^2}\right)^{n-1}\nonumber\\
    &\leq \frac{C}{\eta^2}\int_{A(\eta,\delta)}|\nabla \frak{f}^+|^2
\end{align}
Combining (\ref{A-18}) and (\ref{pestwplus}) gives (\ref{A-16}).
\begin{align}
    \int_{A(\eta,\delta)}|\nabla \frak{f}^-|^2&=4\pi\sum_{n<0}|\frak{h}_n|^2n^2\int_{\delta/\eta}^\eta\rho^{2n-1}d\rho\\
    &=2\pi\sum_{n<0} |\frak{h}_n|^2\,n\ \left(\eta^{2n}-\lf(\frac{\delta}{\eta}\rg)^{2n}\right)\\
    &\geq 2\,\pi\sum_{n<0}|\frak{h}_n|^2\ |n|\ \lf(\lf(\frac{\eta}{\delta}\rg)^{2|n|}-\eta^{2n} \rg)\ge \pi\sum_{n<0}|\frak{h}_n|^2\ |n|\ \lf(\frac{\eta}{\delta}\rg)^{2|n|}
\end{align}
where we assume in the last inequality that $\delta\leq \frac{1}{2}\eta^2$. For $\rho>2\,\delta/\eta$ we have
\begin{align}\label{pestwminus}
    |\nabla \frak{f}^-(\rho,\theta)|^2&\leq\left(2\sum_{n<0}|\frak{h}_n|\,n\,\rho^{n-1}\right)^2\nonumber
    \\&\leq \sum_{n<0}|\frak{h}_n|^2|n|\,\left(\frac{\delta}{\eta}\right)^{2n}\sum_{n>0}\frac{\rho^{-2n-2}}{\left(\frac{\delta}{\eta}\right)^{-2n}}n\nonumber\\
    &=\left( \sum_{n<0}|\frak{h}_n|^2\ |n|\ \lf(\frac{\delta}{\eta}\rg)^{2n}\right)\frac{\delta^2}{\rho^4\eta^2}\sum_{n>0}n\left(\frac{\delta^2}{\eta^2\rho^2}\right)^{n-1}\nonumber\\
    &\leq C\,\frac{\delta^2}{\rho^4\eta^2}\int_{A(\eta,\delta)}|\nabla \frak{f}^-|^2
\end{align}
Combining (\ref{A-18}) and (\ref{pestwminus}) gives (\ref{A-17}). This concludes the proof of lemma~\ref{ptw-harm}.\hfill $\Box$

 \section{Weighted Wente Estimates} 
 \reset
 The following lemma is giving a Wente type estimate with weight. The weight $f$ in the r.h.s. of \eqref{wente-weight}might not be optimal but is sufficient for later purposes.
 
 \begin{Lma}\label{WenteWeight}
 Let $\varphi$ be a solution of
 \begin{equation}\label{dectildeu1}
\left\{ \begin{array}{l}
-\Delta\varphi=\p_{x_1}a\ \p_{x_2}b-\p_{x_2}a\ \p_{x_1}b\quad\quad ~~\mbox{in $B_1(0)$} \\[5mm]
\varphi=0\quad\quad ~~\mbox{on $\quad\p B_1(0)$}
\end{array}\right.\end{equation}
 Then
 \begin{equation}
 \label{wente-weight}
 \int_{B(0,1)}|x|^2|\nabla \varphi|^2  \le C\  \int_{B_1(0)}f(|x|)|\nabla b|^2\ dx^2\ \int_{B_1(0))}|\nabla a|^2\ dx^2\ ,
 \end{equation}
 where 
 \be
 \label{weight-formula}
 f(r)=r^2\log^2\left(1+\frac{1}{r}\right)\log\left(1+\log\frac{1}{r}\right)\ ,
 \ee
 and $C>0$ is a universal constant.
 \hfill $\Box$
 \end{Lma}
In order to prove lemma~\ref{WenteWeight} we need to decompose the Jacobian $\nabla^\perp b\cdot\nabla a$ into a countable sum of Jacobians supported each on a dyadic annulus of the form $A_k:=B_{2^{-k}}(0)\setminus B_{2^{-k-1}}(0)$. We will thus first prove the following lemma
\begin{Lma}
\label{lm-dyadic-jac}
Let $b$ be in $W^{1,2}(B_1(0),{\R})$. There exists $b_k$ such that 
\begin{itemize}
\item[i)]
\be
\label{A-19}
\mbox{supp}(\nabla b_k)\ \subset A_k\cup A_{k+1}\
\ee
\item[ii)]
\be
\label{A-20}
\nabla b=\sum_{k=0}^\infty\nabla b_k
\ee
where the convergence of the partial sum to the series has to be understood in $L^2$
\item[iii)]
\be
\label{A-21}
\int_{B_1(0)}|\nabla b_k|^2\ dx^2\le C\, \int_{ A_k\cup A_{k+1}}|\nabla b|^2\ dx^2
\ee
where $C$ is a universal constant.
\end{itemize}\hfill $\Box$
\end{Lma}
\noindent{\bf Proof of lemma~\ref{lm-dyadic-jac}.} Let $\chi$ be a smooth compactly supported function in ${\R}_+$ such that
\[
\begin{array}{l}
\ds\chi\equiv 1\quad\mbox{ in }[0,1/2]\quad\mbox{ and }\quad \chi\equiv 0\quad\mbox{ in }[1,+\infty)
\end{array}
\]
For any $\rho<1$ we write $\chi_\rho(x):=\chi(|x|/\rho)$. With this notation we have in particular
\[
\mbox{Supp}(\nabla \chi_\rho)\subset B_\rho\setminus B_{\rho/2}\ .
\]
For any $\rho<1/2$ and any $c\in W^{1,2}(B_1)$ such that supp$(\nabla c)\subset B_{2\rho}$ we proceed to the following decomposition
\be
\label{decomp}
\nabla c:=\underbrace{\nabla\bigg(\chi_\rho\ ({c}-\overline{c}_\rho)\bigg)}_{\textbf{(I) }\text{supported in }B_\rho}+\underbrace{\nabla\bigg((1-\chi_\rho)(c-\overline{c}_\rho)\bigg)}_{\textbf{(II) }\text{supported in }B_{2\rho}\setminus B_{\rho/2}}.
\ee
where we take
\[
\ov{c}_\rho:=\frac{1}{|B_\rho\setminus B_{\rho/2}|}\int_{B_\rho\setminus B_{\rho/2}} c(x)\ dx
\]
We denote
\be
\label{op-dec}
T_{\rho}\ :\  c\in W^{1,2}(B_1)\ \longrightarrow T_\rho(c):=\chi_\rho\ ({c}-\overline{c}_\rho)
\ee
Because of Poincar\'e inequality we have
\be
\label{op-a}
\int_{B_\rho\setminus B_{\rho/2}}|\nabla T_\rho(c)|^2\ dx^2\le C\ \int_{B_\rho\setminus B_{\rho/2}}|\nabla c|^2\ dx^2
\ee

We define inductively 
\be
\label{op-b}
b^0:=b\quad\mbox{ , }\quad b^k:=T_{2^{-k}}(b^{k-1})\quad\mbox{ for }k\ge 1\quad\mbox{ and }\quad b_k:=b^k-b^{k+1}\quad\mbox{ for }k\ge 0\ .
\ee
By definition we have
\be
\label{op-c}
\mbox{Supp}(\nabla b^k)\subset B_{2^{-k}}\quad\mbox{ and }\quad\nabla b^k=\nabla b^{k-1}=\nabla b\quad\mbox{ in }\quad B_{2^{-k-1}}\quad .
\ee
This gives for $k\ge 0$
\be
\label{op-d}
\mbox{Supp}(\nabla b_k)\subset B_{2^{-k}}\setminus B_{2^{-k-2}}
\ee
Combining (\ref{op-a}) and (\ref{op-c}) gives moreover
\be
\label{op-e}
\begin{array}{l}
\ds\int_{B_{2^{-k}}\setminus B_{2^{-k-1}}}|\nabla b_k|^2\ dx^2\le 2\,\int_{B_{2^{-k}}\setminus B_{2^{-k-1}}}|\nabla b^k|^2+|\nabla b^{k+1}|^2\ dx^2\\[5mm]
\ds=2\,\int_{B_{2^{-k}}\setminus B_{2^{-k-1}}}|\nabla T_{2^{-k}}(b^{k-1})|^2\ dx^2\le C\ \int_{B_{2^{-k}}\setminus B_{2^{-k-1}}}|\nabla b^{k-1}|^2\ dx^2\\[5mm]
\ds= C\ \int_{B_{2^{-k}}\setminus B_{2^{-k-1}}}|\nabla b|^2\ dx^2
\end{array}
\ee
and
\be
\label{op-f}
\begin{array}{l}
\ds\int_{B_{2^{-k-1}}\setminus B_{2^{-k-2}}}|\nabla b_k|^2\ dx^2\le 2\,\int_{B_{2^{-k-1}}\setminus B_{2^{-k-2}}}|\nabla b^k|^2+|\nabla b^{k+1}|^2\ dx^2\\[5mm]
\ds=2\,\int_{B_{2^{-k-1}}\setminus B_{2^{-k-2}}}|\nabla b|^2+|\nabla T_{2^{-k}}(b^{k})|^2\ dx^2\le C\, \int_{B_{2^{-k-1}}\setminus B_{2^{-k-2}}}|\nabla b|^2+|\nabla b^k|^2\ dx^2\\[5mm]
\ds= 2\,C\ \int_{B_{2^{-k-1}}\setminus B_{2^{-k-2}}}|\nabla b|^2\ dx^2
\end{array}
\ee
Combining (\ref{op-d}), (\ref{op-e}) and (\ref{op-f}) give (\ref{A-19}), (\ref{A-20}) and (\ref{A-21}) and lemma~\ref{lm-dyadic-jac} is proved.\hfill $\Box$

\medskip
Before going to the proof of the weighted Wente estimate a further intermediate lemma is needed.

\begin{Lma}
\label{lm-wente-weight}
Let $j\in {\N}$ and $b\in W^{1,2}(B_1(0))$ such that Supp$(b)\subset B_{2^{-j}}(0)$. Let $a\in W^{1,2}(B_1(0))$ and $\varphi$ be the solution of
\be
\label{A-22}
\left\{ \begin{array}{l}
-\Delta\varphi=\p_{x_1}a\ \p_{x_2}b-\p_{x_2}a\ \p_{x_1}b\quad\quad ~~\mbox{in $B_1(0)$} \\[5mm]
\varphi=0\quad\quad ~~\mbox{on $\quad\p B_1(0)$}
\end{array}\right.
\ee
Then for any $0\le k<j$ there holds
\be
\label{A-23}
\int_{B_{2^{-k}}(0)\setminus B_{2^{-k-1}}(0) }|\nabla\varphi|^2\ dx^2\le C\ 2^{2k-2j}\ \int_{B_1(0)}|\nabla a|^2\ dx^2\ \int_{B_1(0)}|\nabla b|^2\ dx^2
\ee
where $C>0$ is a universal constant. 

Assume now that Supp$(b)\subset B_{2^{-j}}(0)\setminus B_{2^{-j-1}}(0)$ then for any $k>j$ there holds
\be
\label{A-24}
\int_{B_{2^{-k}}(0) }|\nabla\varphi|^2\ dx^2\le C\ 2^{2j-2k}\ \int_{B_1(0)}|\nabla a|^2\ dx^2\ \int_{B_1(0)}|\nabla b|^2\ dx^2
\ee
\hfill $\Box$
\end{Lma}
\noindent{\bf Proof of lemma~\ref{lm-wente-weight}.}
We first prove (\ref{A-23}). Under the hypothesis $\varphi$ is harmonic in the annulus $B_1(0)\setminus B_{2^{-j}}(0)$. In this annulus  $\varphi$ takes the form
\[
\varphi=\varphi_0+ C^0\, \log |x| +\sum_{n\in {\Z}^\ast}\varphi_n x^n+\ov{\varphi_n}\ \ov{x}^n\ .
\]
Since $\varphi=0$ on $\p B_1(0)$ this gives $\varphi_0=0$ and $\varphi_n=-\,\ov{\varphi_{-n}}$. Integrating (\ref{A-22}) over $B_r(0)$ for $r>2^{-j}$ gives also
\be
\label{A-25}
\int_{\p B_r(0)}\p_\nu\varphi\ dl=0\ ,
\ee
which implies that $C^0=0$. Hence we obtain
\be
\label{A-26}
\varphi=\varphi^++\varphi^-\quad\mbox{ where } \quad \varphi^+=\sum_{n>0}\varphi_n x^n+\ov{\varphi_n}\ \ov{x}^n\quad\mbox{ and }\quad\varphi^-=-\sum_{n>0}\ov{\varphi_{n}}\ x^{-n}+{\varphi_{n}}\ \ov{x}^{-n}
\ee
Using lemma~\ref{ptw-harm} we deduce that
\be
\label{A-27}
\begin{array}{l}
\ds\int_{B_{2^{-k}}(0)\setminus B_{2^{-k-1}}(0) }|\nabla\varphi^-|^2\ dx^2\le C\ 2^{2k-2j}\ \int_{B_1}|\nabla\varphi^-|^2\ dx^2\\[5mm]
\ds\quad\le C\ 2^{2k-2j}\ \int_{B_1(0)}|\nabla a|^2\ dx^2\ \int_{B_1(0)}|\nabla b|^2\ dx^2
\end{array}
\ee
Observe that thanks to (\ref{A-26}) one has for $2^{-j}<r<1$ since $|\varphi_n|^2=|\varphi_n|^2$
\be
\label{A-28}
\int_{\p B_r(0)}|\nabla \varphi^+|^2\ dl=4\pi\ \sum_{n>0}n^2\, |\varphi_n|^2 r^{2n}\le 4\pi\ \sum_{n>0}n^2\, |\varphi_{-n}|^2 r^{-2n}=\int_{\p B_r(0)}|\nabla \varphi^-|^2\ dl
\ee
Then (\ref{A-23}) follows from (\ref{A-27}) and (\ref{A-28}).

\medskip

Assuming now Supp$(b)\subset B_{2^{-j}}(0)\setminus B_{2^{-j-1}}(0)$, we have that $\varphi$ is harmonic on $B_{2^{-j-1}}(0)$. The monotonicity formula for harmonic map implies that
\be
\label{A-29}
\begin{array}{l}
\ds\forall r<2^{-j-1}\quad\quad \frac{1}{r^2}\int_{B_r(0)}|\nabla \varphi|^2\ dx^2\le 2^{2j+2}\ \int_{B_{2^{-j-1}}(0)}|\nabla \varphi|^2\ dx^2\\[5mm]
\ds\quad\le C\ 2^{2j+2}\ \int_{B_1(0)}|\nabla a|^2\ dx^2\ \int_{B_1(0)}|\nabla b|^2\ dx^2
\end{array}
\ee
and (\ref{A-24}) follows by taking $r=2^{-k}$. This concludes the proof of lemma~\ref{lm-wente-weight}.\hfill $\Box$

\medskip

\noindent{\bf Proof of lemma~\ref{WenteWeight}} 
Let us write $$\varphi=\sum_{k=0}^{\infty}\varphi_k\ ,$$ where
\begin{equation}
\label{decompo}
\lf\{
    \begin{array}{l}
    -\Delta\varphi_k=\nabla^\perp a\cdot\nabla b_k \text{ in }B_1\\[5mm]
    \varphi_k=0\text{ in }\partial B_1.
    \end{array}
    \rg.
\end{equation}
where $b_k$ are given by  lemma~\ref{lm-dyadic-jac}. We keep denoting $A_k=B_{2^{-k}}(0)\setminus B_{2^{-(k+1)}}(0)$ and we also introduce the following notation for slightly larger anuli
$\ti{A}_k=B_{2^{-k+1}}(0)\setminus B_{2^{-(k+1)}}(0)$.

\medskip
 We then have by the Dominated Convergence Theorem,
\begin{align}
    \int_{A_k}|\nabla \varphi|^2&=\int_{A_k}\sum_{i,j}\nabla\varphi_j\cdot\nabla \varphi_j=\sum_{i,j}\int_{A_k}\nabla\varphi_i\cdot\nabla\varphi_j,
\end{align}
and
\begin{equation}
    \sum_{k=0}^\infty\int_{A_k}|x|^2|\nabla\varphi|^2\ dx^2=\sum_{k=0}^\infty 2^{-2k}\,\int_{A_k}|\nabla\varphi|^2\ dx^2\ \cong\ \sum_{k=0}^\infty
2^{-2k}\sum_{i,j}\int_{A_k}\nabla\varphi_i\cdot\nabla\varphi_j.
\end{equation}
We now consider the case of each possible pairing of frequencies separately.
\be
\label{freq}
\begin{array}{l}
  \ds  \sum_{k=0}^\infty2^{-2k}\int_{A_k}|\nabla\varphi|^2=\sum_{k=0}^{\infty}2^{-2k}\bigg(\int_{A_k}\sum_{|i-k|<2,|j-k|<2}\nabla\varphi_i\cdot\nabla\varphi_j
 \ds   +2\int_{A_k}\sum_{|i-k|<2,j>k+2}\nabla\varphi_i\cdot\nabla\varphi_j\\[5mm]
 \ds   +2\int_{A_k}\sum_{|i-k|<2,j<k-2}\nabla\varphi_i\cdot\nabla\varphi_j
    +2\int_{A_k}\sum_{|i-k|>2,|j-k|>2}\nabla\varphi_i\cdot\nabla\varphi_j\bigg)\ .
    \end{array}
\ee
We study each of the 4 cases separately

\medskip

\textbf{Case 1:} Using (\ref{A-21}) we obtain
\be
\label{case1}
\begin{array}{l}
 \ds   \sum_{k=0}^{\infty}2^{-2k}\sum_{\substack{|i-k|<2\\|j-k|<2}}\int_{A_k}\nabla\varphi_i\cdot\nabla\varphi_j\leq \sum_{k=0}^{\infty}2^{-2k}\sum_{\substack{|i-k|<2\\|j-k|<2}}\bigg(\int_{A_k}|\nabla\varphi_i|^2\bigg)^{1/2}\bigg(\int_{A_k}|\nabla\varphi_j|^2\bigg)^{1/2}\\[5mm]
 \ds      \leq C\sum_{k=0}^{\infty}2^{-2k}\sum_{|i-k|<2}\int_{A_k}|\nabla\varphi_i|^2
    \leq C \sum_{k=0}^{\infty}2^{-2k}\sum_{|i-k|<2}\int_{\tilde{A}_i}|\nabla b_i|^2\int_{D^2}|\nabla a|^2\\[5mm]
\ds    \leq C\sum_{k=0}^\infty\int_{\tilde{A}_k}|\nabla b_k|^22^{-2k}\int_{D^2}|\nabla a|^2
    \leq C\int_{D^2}|\nabla b|^2|x|^2\int_{D^2}|\nabla a|^2.
    \end{array}
\ee
Notice that the weight for these frequencies can be taken to be $f(r)=r^2.$\\
\textbf{Case 2:} Lemma~\ref{lm-wente-weight} is implying
\begin{align}
    &2\sum_{k=0}^{\infty}2^{-2k}\sum_{\substack{|i-k|\leq2\\j>k+2}}\int_{A_k}\nabla\varphi_i\cdot\nabla\varphi_j\leq   2\sum_{k=0}^{\infty}2^{-2k}\sum_{|i-k|\leq2}\bigg(\int_{A_k}|\nabla\varphi_i|^2\bigg)^{1/2}\sum_{j>k+2}\bigg(\int_{A_k}|\nabla\varphi_j|^2\bigg)^{1/2}\\
    &\leq \int_{D^2}|\nabla a|^2\sum_{k=0}^\infty 2^{-2k}\bigg(\sum_{|i-k|\leq 2}\bigg(\int_{\tilde{A}_i}|\nabla b_i|^2\bigg)^{1/2}\bigg)\bigg(\sum_{j>k+2}2^{-(j-k)}\bigg(\int_{\tilde{A}_j}|\nabla b_j|^2\bigg)^{1/2}\bigg)\\
    &\leq \int_{D^2}|\nabla a|^2\sum_{k=0}^\infty 2^k 2^{-2k}\bigg(\int_{\tilde{A}_k}|\nabla b_k|^2\bigg)^{1/2}\sum_{j>k+2}\bigg(\int_{\tilde{A}_j}|\nabla b_j|^22^{-2j}\bigg)^{1/2}\\
    &\leq \int_{D^2}|\nabla a|^2\sum_{k=0}^\infty \bigg(\int_{\tilde{A}_k}2^{-2k}|\nabla b_k|^2\bigg)^{1/2}\sum_{j>k+2}\frac{1}{j}\bigg(\int_{\tilde{A}_j}|\nabla b_j|^22^{-2j}\log^2(1+2^j)\bigg)^{1/2}\\
    &\leq \int_{D^2}|\nabla a|^2\sum_{k\geq 1}\frac{1}{k}\bigg(\int_{\tilde{A}_k}|\nabla b_k|^22^{-2k}\log^2(1+2^k)\bigg)^{1/2}\bigg(\sum_{j=1}^\infty \int_{\tilde{A}_j}2^{-2j}|\nabla b_j|^2\log^2(1+2^j)\bigg)^{1/2}\\
    &\leq \int_{D^2}|\nabla a|^2\lf[\sum_{k\geq 1}\frac{1}{k^2}\rg]\ \sum_{k=1}^\infty \int_{\tilde{A}_k}2^{-2k}|\nabla b_k|^2\log^2(1+2^k)\\
    &\leq \int_{D^2}|\nabla a|^2\int_{D^2}|\nabla b|^2\log^2\lf(1+\frac{1}{|x|}\rg)\ |x|^2\\
\end{align}
\textbf{Case 3:} Lemma~\ref{lm-wente-weight} again is implying
\begin{align}
    &\sum_{k=0}^{\infty}2^{-2k}\sum_{\substack{|i-k|<2\\j<k-2}}\int_{A_k}\nabla\varphi_i\cdot\nabla\varphi_j\leq \sum_{k=0}^{\infty}2^{-2k}\sum_{|i-k|<2}\bigg(\int_{A_k}|\nabla\varphi_i|^2\bigg)^{1/2}\sum_{0\leq j\leq k-2}\bigg(\int_{A_k}|\nabla\varphi_j|^2\bigg)^{1/2}\\
    &\leq C\sum_{k=0}^{\infty}2^{-2k}\bigg(\int_{A_k}|\nabla\varphi_k|^2\bigg)^{1/2}\sum_{0\leq j\leq k-2}\bigg(\int_{A_k}|\nabla\varphi_j|^2\bigg)^{1/2}\\
    &\leq C\sum_{k=0}^{\infty}2^{-2k}\bigg(\int_{A_k}|\nabla\varphi_k|^2\bigg)^{1/2}\sum_{0\leq j\leq k-2}2^{-k}2^j\bigg(\int_{\tilde{A}_j}|\nabla b_j|^2\int_{D^2}|\nabla a|^2\bigg)^{1/2}\\
    &\leq C\sum_{k=0}^{\infty}2^{-k}\bigg(\int_{A_k}|\nabla b_k|^22^{-2k}\bigg)^{1/2}\sum_{0\leq j\leq k-2}2^{-k}2^{2j}\bigg(\int_{\tilde{A}_j}|\nabla b_j|^22^{-2j}\bigg)^{1/2}\ \int_{D^2}|\nabla a|^2
 \end{align}   
  and we have      
  \begin{align}  
    &\sum_{k=0}^{\infty}2^{-k}\bigg(\int_{A_k}|\nabla b_k|^22^{-2k}\bigg)^{1/2}\sum_{0\leq j\leq k-2}2^{-k}2^{2j}\bigg(\int_{\tilde{A}_j}|\nabla b_j|^22^{-2j}\bigg)^{1/2}\\
    &\leq C\bigg(\sum_{k=0}^{\infty}\int_{A_k}|\nabla b_k|^22^{-2k}\bigg)^{1/2}\bigg(\sum_{k=0}^\infty 2^{-4k}\bigg(\sum_{0\leq j\leq k-2}2^{-k}2^{2j}\bigg(\int_{\tilde{A}_j}|\nabla b_j|^22^{-2j}\bigg)^{1/2}\bigg)^2\bigg)^{1/2}\\
    &\leq C\bigg(\sum_{k=0}^{\infty}\int_{A_k}|\nabla b_k|^22^{-2k}\bigg)^{1/2}\lf(\sum_{k=0}^\infty 2^{-4k}\lf[\sum_{j\leq k}\frac{2^{4j}}{j^2}\rg]\ \sum_{0\leq j\leq k-1}\int_{\tilde{A}_j}|\nabla b_j|^22^{-2j}\log^2(1+2^j)\rg)^{1/2}\\
    &\leq C\bigg(\sum_{k=0}^{\infty}\int_{A_k}|\nabla b_k|^22^{-2k}\log^2(1+2^k)\bigg)\bigg(\sum_{k=0}^\infty 2^{-4k}\sum_{j\leq k}\frac{2^{4j}}{j^2}\bigg)^{1/2}\\
    &\leq C\int_{D^2}|\nabla b|^2|x|^2\log^2\lf(1+\frac{1}{|x|}\rg)\\
\end{align}
\textbf{Case 4:}
\begin{equation}
    \sum_{|i-k|>2,|j-k|>2}\int_{A_k}\nabla\varphi_i\cdot\nabla\varphi_j
\end{equation}
We treat each of the four sub-cases separately and we make again an intensive use of lemma~\ref{lm-wente-weight}.\\
\textbf{Case 4.1:} $(i>k+2,j>k+2)$
\begin{align}
    &\sum_{k=0}^\infty2^{-2k}\sum_{i=k+2}^\infty\sum_{j=k+2}^\infty\int_{A_k}\nabla\varphi_i\cdot\nabla\varphi_j\\
    \leq& \sum_{k=0}^\infty2^{-2k}\sum_{i=k+2}^\infty\sum_{j=k+2}^\infty\bigg(\int_{A_k}|\nabla\varphi_i|^2\bigg)^{1/2}\bigg(\int_{A_k}|\nabla\varphi_j|^2\bigg)^{1/2}\\
    \leq& \sum_{k=0}^\infty2^{-2k}\sum_{i=k+2}^\infty\sum_{j=k+2}^\infty 2^{-(i-k)}2^{-(j-k)}\bigg(\int_{\tilde{A}_i}|\nabla b_i|^2\bigg)^{1/2}\bigg(\int_{\tilde{A}_j}|\nabla b_j|^2\bigg)^{1/2}\int_{D^2}|\nabla a|^2\\
    \leq& \sum_{k=0}^\infty\bigg(\sum_{i=k+2}^\infty 2^{-i}\bigg(\int_{\tilde{A}_i}|\nabla b_i|^2\bigg)^{1/2}\bigg)^2\int_{D^2}|\nabla a|^2\\
    \leq& \sum_{k=0}^\infty\bigg(\sum_{i=k+2}^\infty\frac{1}{i^2}\bigg)\bigg(\sum_{i=k+2}^\infty \int_{\tilde{A}_i}|\nabla b_i|^22^{-2i}\ i^2\bigg)\int_{D^2}|\nabla a|^2\\   
    \leq& \sum_{k=0}^\infty\frac{1}{k+2}\sum_{i=k+2}^\infty \int_{\tilde{A}_i}|\nabla b_i|^22^{-2i}\log^2(1+2^i)\int_{D^2}|\nabla a|^2\\  
    =&\bigg(\sum_{i=2}^\infty \int_{\tilde{A}_i}|\nabla b_i|^22^{-2i}\log^2(1+2^i)\sum_{k=0}^{i-2}\frac{1}{k+2}\bigg)\int_{D^2}|\nabla a|^2\\ 
    \leq&\bigg(\sum_{i=2}^\infty \int_{\tilde{A}_i}|\nabla b_i|^22^{-2i}\log^2(1+2^i)\log(1+i)\bigg)\int_{D^2}|\nabla a|^2\\
    \leq&C \int_{D^2}|\nabla b|^2|x|^2\log^2\lf(1+\frac{1}{|x|}\rg)\log(1+\log(|x|^{-1}))\int_{D^2}|\nabla a|^2\\
\end{align}
\textbf{Case 4.2:} $(i<k-2,j<k-2)$
Recall that since $\Delta\varphi_i=0$ in $B_{2^{-k}}(0)$ we have
\begin{align}
    \int_{A_k}|\nabla\varphi_i|^2\ dx^2&\leq C2^{-2k}2^{2i}\int_{A_i}|\nabla\varphi_i|^2\ dx^2\\
    &\leq 2^{-2k}2^{4i}\int_{\tilde{A}_i}|\nabla b_i|^2 \,2^{-2i}\, dx^2\ \int_{D^2}|\nabla a|^2\ dx^2.
\end{align}
We can now estimate
\begin{align}
   & \sum_{k=0}^\infty2^{-2k}\sum_{i,j<k-2}\int_{A_k}\nabla\varphi_i\cdot\nabla\varphi_j\\
   \leq &\sum_{k=0}^\infty2^{-2k}\bigg(\sum_{i<k-2}\bigg(\int_{A_k}|\nabla\varphi_i|^2\bigg)^{1/2}\bigg)^2\\
   \leq &\sum_{k=0}^\infty2^{-2k}\bigg(\sum_{i<k-2}2^{-k}2^{2i}\bigg(\int_{\tilde{A}_i}|\nabla b_i|^22^{-2i}\bigg)^{1/2}\bigg)^2\int_{D^2}|\nabla a|^2\\
    \leq & C\sum_{k=0}^\infty2^{-4k}\bigg(\sum_{i<k-2}\frac{2^{4i}}{i^2}\bigg)\bigg(\sum_{i<k-2}\int_{\tilde{A}_i}|\nabla b_i|^22^{-2i}\log^2(1+2^i)\bigg)\int_{D^2}|\nabla a|^2\\
    \leq & C\sum_{k=0}^\infty2^{-4k}\sum_{0<i<k-2}\frac{2^{4i}}{i^2}\bigg(\sum_{i=0}^\infty\int_{\tilde{A}_i}|\nabla b_i|^22^{-2i}\log^2(1+2^i)\bigg)\int_{D^2}|\nabla a|^2\\
    \leq & C\int_{D^2}|\nabla b|^2|x|^2\log^2\lf(1+\frac{1}{|x|}\rg)\ \int_{D^2}|\nabla a|^2\\
\end{align}
\textbf{Case 4.3:} $(i<k-1,j>k+2)$
\begin{align}
    &\sum_{k=0}^\infty 2^{-2k}\sum_{\substack{j>k+2\\ i<k-1}}\int_{A_k}\nabla\varphi_i\cdot\nabla\varphi_j\\
    \leq& \sum_{k=0}^\infty 2^{-2k}\sum_{i<k-2}\bigg(\int_{A_k}\nabla\varphi_i|^2\bigg)^{1/2}\bigg(\sum_{j>k+2}\bigg(\int_{A_k}|\nabla\varphi_j|^2\bigg)^{1/2}\bigg)\\
    \leq& \sum_{k=0}^\infty 2^{-2k}\bigg(\sum_{i<k-2}2^{-k}2^{2i}\bigg(\int_{\tilde{A}_i}|\nabla b_i|^2 2^{-2i}\bigg)^{1/2}\bigg)\sum_{j>k+2}2^k\bigg(\int_{\tilde{A}_j}|\nabla b_j|^22^{-2j}\bigg)^{1/2}\int_{D^2}|\nabla a|^2\\
    \leq& \sum_{k=0}^\infty 2^{-2k}\bigg(\sum_{i<k-2}\frac{2^{4i}}{i^2}\bigg)^{1/2}\bigg(\sum_{i<k-2}\int_{\tilde{A}_i}|\nabla b_i|^22^{-2i}\log^2(1+2^i)\bigg)^{1/2}\\
    &\bigg(\sum_{j>k+2}\frac{1}{j^2}\bigg)^{1/2}\bigg(\sum_{j>k+2}\int_{\tilde{A}_j}|\nabla b_j|^22^{-2j}\log^2(1+2^j)\bigg)^{1/2}\int_{D^2}|\nabla a|^2\\
    \leq& \bigg(\sum_{k=0}^\infty 2^{-4k}\sum_{i<k-2}\frac{2^{4i}}{i^2}\bigg)^{1/2}\bigg(\sum_{k=0}^\infty\bigg(\sum_{j>k+2}\frac{1}{j^2}\bigg)\sum_{j>k+2}\int_{\tilde{A}_j}|\nabla b_j|^22^{-2j}\log^2(1+2^j)\bigg)^{1/2}\\
    &\bigg(\sum_{i<\infty}\int_{\tilde{A}_i}|\nabla b_i|^22^{-2i}\log^2(1+2^i)\bigg)^{1/2}\int_{D^2}|\nabla a|^2\\
\end{align}
With the estimate
\begin{align}
    &\sum_{k=0}^\infty\bigg(\sum_{j>k+2}\frac{1}{j^2}\bigg)\sum_{j>k+2}\int_{\tilde{A}_j}|\nabla b_j|^22^{-2j}\log^2(1+2^j)\\
    \leq &\sum_{k=0}^\infty\frac{1}{k+1}\sum_{j>k+2}\int_{\tilde{A}_j}|\nabla b_j|^22^{-2j}\log^2(1+2^j)\\
    =&\sum_{j=1}^\infty\int_{\tilde{A}_j}|\nabla b_j|^22^{-2j}\log^2(1+2^j)\sum_{k=1}^j\frac{1}{k+1}\\
    \leq&\sum_{j=1}^\infty\int_{\tilde{A}_j}|\nabla b_j|^22^{-2j}\log^2(1+2^j)\log(1+j)
    \end{align}
we obtain the desired result. Collecting all the estimates above give the weighted Wente inequality~\ref{wente-weight} and lemma~\ref{WenteWeight} is proved.\hfill $\Box$

\section{Morrey decrease of Annulus Energies for Solutions to Wente type Equation.}
\reset
We start with some preliminary results.
\begin{Lma}\label{morrey1}
Let  $a,b\in W^{1,2}(B_1(0))$ and  $\varphi$ be a solution of
\begin{equation}
-\Delta\varphi=\p_{x_1}a\ \p_{x_2}b-\p_{x_2}a\ \p_{x_1}b\quad\mbox{in $B_1(0)$}
\end{equation}
Then the following estimate hold.
\begin{eqnarray}\label{morrey2}
\int_{B_{\frac{1}{2}}(0)\setminus B_{\frac{1}{4}}(0)}|\nabla \varphi|^2\ dx^2&\le &\frac{2}{3}\int_{B_1(0)\setminus B_{\frac{1}{2}}(0)}|\nabla \varphi|^2+C\ \int_{B_1(0)}|\nabla a|^2\ dx^2\ \int_{B_1(0)} f(|x|)\ |\nabla b|^2\ dx^2
\end{eqnarray}
where $f$ is given by (\ref{weight-formula}) and $C>0$ is a universal constant. 
\hfill $\Box$
\end{Lma}
\noindent{\bf Proof of Lemma \ref{morrey1}.}
We decompose $\varphi=\sigma+\xi$ with $\sigma,$  $\xi$
 satisfying respectively
 \begin{equation}\label{systsigmaxi}
 \left\{\begin{array}{cc}
\ds \Delta\sigma =0&~~\mbox{in $B_1(0) $}\\[5mm]
\ds\sigma =\varphi & ~~\mbox{on $\partial B_1(0) $}\end{array}\right. 
  ~~\quad\mbox{and}\quad~~~ \left\{\begin{array}{cc}
\ds \Delta\xi =\Delta \varphi &~~\mbox{in $B_1(0) $}\\[5mm]
\xi =0 &~~\mbox{on $\partial B_1(0) $}\end{array}\right.\end{equation}
Since $\sigma $ is harmonic the monotonicity formula for harmonic functions gives
\begin{eqnarray}
\int_{B_{\frac{1}{2}}(0)}|\nabla \sigma|^2\ dx^2&\le& \frac{1}{4}\int_{B_1(0)}|\nabla \sigma|^2\ dx^2\nonumber\\&=&
\frac{1}{4}\int_{B_1(0)\setminus B_{\frac{1}{2}}(0)}|\nabla \sigma|^2\ dx^2+\frac{1}{4}\int_{B_{\frac{1}{2}}(0)}|\nabla \sigma|^2\ dx^2\ .\label{harmestring}
\end{eqnarray}
From \eqref{harmestring} it follows that
\be
\label{harm-decrease}
\int_{B_{\frac{1}{2}}(0)}|\nabla \sigma|^2\ dx^2 \le \frac{1}{3}\int_{B_1(0)\setminus B_{\frac{1}{2}}(0)}|\nabla \sigma|^2\ dx^2
\ee
We also have by the Dirichlet principle that
\begin{equation}\label{extensionharm}
\int_{B_1(0)}|\nabla \sigma|^2\ dx^2 \le \int_{B_1(0)}|\nabla \varphi |^2\ dx^2
\end{equation}
As far as the function $\xi$ is concerned we can apply the weighted Wente Lemma~\ref{WenteWeight} and deduce that 
\begin{equation}
\int_{B_{1}(0)\setminus B_{\frac{1}{4}}(0)}|\nabla \xi|^2\ dx^2\le C \int_{B_1(0)}|\nabla a|^2\ dx^2\ \int_{B_1(0)}\ f(|x|)\ |\nabla b|^2\ dx^2\ .
\end{equation}
Next we have combining the previous inequalities
\be
\label{morrey4}
\begin{array}{l}
\ds\int_{B_{\frac{1}{2}}(0)\setminus B_{\frac{1}{4}}(0)}|\nabla \varphi|^2\ dx^2=\ds\int_{B_{\frac{1}{2}}(0)\setminus B_{\frac{1}{4}}(0)}|\nabla \xi+\nabla \sigma|^2\ dx^2\\[5mm]
\ds\quad\le\ds \frac{3}{2} \ \int_{B_{\frac{1}{2}}(0)\setminus B_{\frac{1}{4}}(0)}|\nabla \sigma|^2\ dx^2+3\,
 \int_{B_{\frac{1}{2}}(0)\setminus B_{\frac{1}{4}}(0)}|\nabla \xi|^2\ dx^2 \\[5mm]
 \ds\quad\le \ds \frac{1}{2} \int_{B_1(0)\setminus B_{\frac{1}{2}}(0)}|\nabla \sigma|^2+3\,
 \int_{B_{\frac{1}{2}}(0)\setminus B_{\frac{1}{4}}(0)}|\nabla \xi|^2\ dx^2 \\[5mm]
 \ds\quad\le \ds \frac{2}{3} \int_{B_1(0)\setminus B_{\frac{1}{2}}(0)}|\nabla \varphi|^2+9\,
 \int_{B_{1}(0)\setminus B_{\frac{1}{4}}(0)}|\nabla \xi|^2\ dx^2 \\[5mm]
\ds\quad \le \ds \frac{2}{3} \int_{B_1(0)\setminus B_{\frac{1}{2}}(0)}|\nabla \varphi|^2+C\ \int_{B_1(0)}|\nabla a|^2\ dx^2\ \int_{B_1(0)}\ f(|x|)\ |\nabla b|^2\ dx^2\ . 
\end{array}
\ee
This concludes the proof of  lemma \ref{morrey1}.\hfill $\Box$

\begin{Lma}\label{morrey9bis}
Under the hypothesis of lemma \ref{morrey1}, for any $\al\in(0,2)$ there exists $C_\al>0$ such that for any $k\in {\N}$
 \begin{eqnarray}\label{morrey9}
{\int_{ A_{k+1}}}|\nabla \varphi|^2\ dx^2&\le&\gamma^{k+1}{\int_{ A_{0}}}|\nabla \varphi|^2\ dx^2+C_\al\ \int_{B_1(0)}|\nabla a|^2\ dx^2\
\sum_{n=0}^{\infty}\gamma^{|n-k|}\int_{A_n}|\nabla b|^2\ dx^2\ .
\end{eqnarray}
where $\gamma=\max(2^{-\alpha},\frac{2}{3})$ and we are making use of the notation $A_k:=B_{2^{-k}}(0)\setminus B_{2^{-k-1}}(0)$.\hfill $\Box$
\end{Lma}
\noindent{\bf Proof of Lemma \ref{morrey9bis}.} By rescaling and iterating the estimate \eqref{morrey4} we get for any $k\ge 0$
 \begin{eqnarray}\label{morrey5}
{\int_{ A_{k+1}}}|\nabla \varphi|^2&\le&\frac{2}{3}{\int_{ A_{k}}}|\nabla \varphi|^2+C\ \int_{B_1(0)}|\nabla a|^2\ dx^2 \int_{B_{2^{-k}}(0)}f(2^{k}|x|)\,|\nabla b|^2\ dx^2\nonumber\\
&\le&\frac{2}{3}\left[\frac{2}{3}{\int_{ A_{k-1}}}|\nabla \varphi|^2+C\int_{B_1(0)}|\nabla a|^2\ dx^2\ \int_{B_{2^{-k+1}}(0)}f(2^{k-1}|x|)\ |\nabla b|^2\ dx^2\right]\nonumber\\\
&+&C\int_{B_1(0)}|\nabla a|^2\ dx^2\  \int_{B_{2^{-k}}(0)}f(2^{k}|x|)\ |\nabla b|^2\ dx^2\\
&&\ldots\nonumber\\
&\le&\left(\frac{2}{3}\right)^{k+1}{\int_{ A_{0}}}|\nabla \varphi|^2+C\ \int_{B_1(0)}|\nabla a|^2\ dx^2\
\sum_{j=0}^{k}\left(\frac{2}{3}\right)^j\int_{B_{2^{-k+j}}(0)}\ f(2^{k-j}|x|)\ |\nabla b|^2\ dx^2\ .\nonumber
\end{eqnarray}
 We observe that $f(r)=o(r^{\alpha})$ as $r\to 0$ for every $\alpha\in (0,2).$ We have
 \begin{eqnarray}\label{morrey10}
&&\sum_{j=0}^{k}\left(\frac{2}{3}\right)^j\int_{B_{2^{-k+j}}(0)}\,2^{\alpha(k-j)}|x|^{\alpha}|\nabla b|^2\ dx^2\nonumber\\
&\lesssim&
\sum_{j=0}^{k}\left(\frac{2}{3}\right)^j\sum_{n=k-j}^{\infty}\int_{A_n} 2^{\alpha(k-j-n)} |\nabla b|^2\ dx^2
\nonumber\\
&=&\sum_{j=0}^{k}\left(\frac{2}{3}\right)^j2^{\alpha(k-j)}\sum_{n=k-j}^{\infty}2^{-n\alpha} \int_{A_n}  |\nabla b|^2\ dx^2\nonumber\\
&=&\sum_{i=0}^{k}\left(\frac{2}{3}\right)^{k-i}\ 2^{\alpha i}\ \sum_{n=i}^{\infty}2^{-n\alpha} \int_{A_n}  |\nabla b|^2\ dx^2 \nonumber\\
&=&\left(\frac{2}{3}\right)^{k}\sum_{n=0}^{\infty}2^{-n\alpha}\,\sum_{i=0}^{\inf(n,k)}\left(2^{\alpha}\frac{3}{2}\right)^{i}\ \int_{A_n}  |\nabla b|^2\ dx^2\nonumber\\\
&\le &\left(\frac{2}{3}\right)^{k}\sum_{n=0}^{\infty}2^{-n\alpha}\left(2^{\alpha}\frac{3}{2}\right)^{\inf(n,k)} \int_{A_n}  |\nabla b|^2\ dx^2 \nonumber\\
&=& \left(\frac{2}{3}\right)^{k}\sum_{n=0}^{k}2^{-n\alpha}\left(2^{\alpha}\frac{3}{2}\right)^{n} \int_{A_n}  |\nabla b|^2\ dx^2 \nonumber\\
&+&\left(\frac{2}{3}\right)^{k}\sum_{n=k+1}^{\infty}2^{-n\alpha}\left(2^{\alpha}\frac{3}{2}\right)^{k} \int_{A_n}  |\nabla b|^2\ dx^2 \nonumber\\
&=&\sum_{n=0}^{k} \left(\frac{3}{2}\right)^{n-k} \int_{A_n}  |\nabla b|^2\ dx^2+ \sum_{n=k+1}^{\infty}2^{(k-n)\alpha}  \int_{A_n}  |\nabla b|^2\ dx^2\nonumber\\
&\le& \sum_{n=0}^{\infty}\gamma^{|n-k|}\int_{A_n}|\nabla b|^2\ dx^2\ .
\end{eqnarray}
Then the estimate \eqref{morrey9} holds. This concludes the proof of the Lemma \ref{morrey9bis}.\hfill$\Box$
\section{A  lemma on weighted series}
\reset
 \begin{Lma}\label{morrey11}
Let $\gamma\in (0,1)$. For any $s_1<s_2$ in ${\N}$ for any pair of sequences of non negative numbers $(a_n)_{n\in{\N}}$ and $(b_n)_{n\in{\N}}$ satisfying for $s_1\le k\le s_2$
\begin{equation}\label{morrey12bis}
a_k\le b_k+\varepsilon_0\,\sum_{n=0}^{\infty}\gamma^{|n-k|} a_n\ .
\end{equation}
for some $\varepsilon_0>0$. Then for any $\gamma<\mu<1$  there exists $C_{\mu,\gamma}>0$ such that
\begin{eqnarray}\label{morrey12}
\sum_{\ell=s_1}^{s_2}\mu^{|\ell-k|}\,a_{\ell} &\le&  \sum_{\ell=s_1}^{s_2}\mu^{|\ell-k|}b_{\ell}+C_{\mu,\gamma}\ \varepsilon_0\  \sum_{\ell=0}^{\infty}  \mu^{|\ell-k|}\,a_\ell  
 \end{eqnarray}
\hfill $\Box$
\end{Lma}
\noindent{\bf Proof of lemma~\ref{morrey11}.}
From \eqref{morrey12bis} we get:
\begin{equation}\label{morrey12tris}
\sum_{\ell=s_1}^{s_2}\mu^{|\ell-k|}a_{\ell} \le  \sum_{\ell=s_1}^{s_2}\mu^{|\ell-k|}b_{\ell}+\varepsilon_0\sum_{n=0}^{\infty}a_n\lf[\sum_{\ell=s_1}^{s_2}\gamma^{|n-\ell|}\mu^{|\ell-k|}\rg]
\end{equation}
\noindent{\bf Case 1.}  $s_1\le n<k$\par
\begin{eqnarray}
&&\sum_{\ell=s_1}^{s_2}\gamma^{|n-\ell|}\mu^{|\ell-k|}=\sum_{\ell=s_1}^{n}\gamma^{|n-\ell |}\mu^{|\ell-k|}\nonumber\\
&+&\sum_{\ell=n}^{k}\gamma^{|n-\ell|}\mu^{|\ell-k|}+\sum_{\ell=k+1}^{s_2}\gamma^{|n-\ell|}\mu^{|\ell-k|}\nonumber\\
&=&\sum_{\ell=s_1}^{n}\gamma^{n-\ell}\mu^{k-\ell}+\sum_{\ell=n}^{k}\gamma^{\ell-n}\mu^{k-\ell}+\sum_{\ell=k+1}^{s_2}\gamma^{\ell-n}\mu^{\ell-k}\nonumber\\
&\le&\gamma^n\mu^k \sum_{\ell=s_1}^{n}(\mu\gamma)^{-\ell}+
\gamma^{-n}\mu^k \sum_{\ell=n}^{k}(\mu^{-1}\gamma)^{\ell}+\gamma^{-n}\mu^{-k}\sum_{\ell=k+1}^{\infty}(\mu\gamma)^\ell\nonumber\\
&\le&\lf[\frac{(\mu\gamma)^{-s_1+n}}{1-\mu\gamma}+\frac{1}{1-\mu^{-1}\gamma}\rg]\,\mu^{k-n}+\frac{\mu\gamma}{1-\mu\gamma}\,\gamma^{k-n}\le C_{\mu,\gamma}\,  \mu^{k-n}.
\end{eqnarray}
\noindent{\bf Case 2.}  $n\le s_1$\par
\begin{eqnarray}
&&\sum_{\ell=s_1}^{s_2}\gamma^{|n-\ell|}\mu^{|\ell-k|}=\sum_{\ell=s_1}^{k}\gamma^{|n-\ell |}\mu^{|\ell-k|}+\sum_{\ell=k+1}^{s_2}\gamma^{|n-\ell|}\mu^{|\ell-k|}\nonumber\\
&=&\sum_{\ell=s_1}^{k}\gamma^{\ell-n}\mu^{k-\ell}+\sum_{\ell=k+1}^{s_2}\gamma^{\ell-n}\mu^{\ell-k}\nonumber\\
&\le&
\gamma^{-n}\mu^k \sum_{\ell=n}^{k}(\mu^{-1}\gamma)^{\ell}+\gamma^{-n}\mu^{-k}\sum_{\ell=k+1}^{\infty}(\mu\gamma)^\ell\nonumber\\
&\le&\frac{1}{1-\mu^{-1}\gamma}\,\mu^{k-n}+\frac{\mu\gamma}{1-\mu\gamma}\,\gamma^{k-n}\le C_{\mu,\gamma}\,  \mu^{k-n}.
\end{eqnarray}
{\bf Case 3.}  $ k\le n\le s_2$\par
\begin{eqnarray}\label{morrey13}
&&\sum_{\ell=s_1}^{s_2}\gamma^{|n-\ell|}\mu^{|\ell-k|}=\sum_{\ell=s_1}^{k}\gamma^{n-\ell}\mu^{k-\ell}+\sum_{\ell=k+1}^{n}\gamma^{n-\ell}\mu^{\ell-k} +\sum_{\ell=n+1}^{s_2}\gamma^{\ell-n}\mu^{\ell-k} 
\nonumber\\
&\le& \frac{\gamma^{n-k}}{1-\gamma\,\mu}+\frac{\mu^{n-k}}{1-\gamma\mu^{-1}}+\frac{\mu^{n-k}}{1-\gamma\,\mu}\le C_{\mu,\gamma}\,  \mu^{n-k}.
\end{eqnarray}
{\bf Case 4.}  $ s_2\le n$\par
\begin{eqnarray}\label{morrey14}
&&\sum_{\ell=s_1}^{s_2}\gamma^{|n-\ell|}\mu^{|\ell-k|}=\sum_{\ell=s_1}^{k}\gamma^{n-\ell}\mu^{k-\ell}+\sum_{\ell=k+1}^{s_2}\gamma^{n-\ell}\mu^{\ell-k} 
\nonumber\\
&\le& \frac{\gamma^{n-k}}{1-\gamma\,\mu}+\frac{\gamma^{n-s_2}\,\mu^{-k+s_2}}{1-\gamma\mu^{-1}}\le C_{\mu,\gamma}\,  \mu^{n-k}.
\end{eqnarray}
Combining all the above cases and (\ref{morrey12tris})  gives (\ref{morrey12}) and lemma~\ref{morrey11} is proved. \hfill $\Box$

\end{document}